\documentclass[a4paper,12pt]{article}
\pdfoutput=1
\usepackage[paper=letterpaper,margin=1in]{geometry}
\usepackage{amsmath,amssymb,amsfonts,amsthm,mathtools,epsfig,cite,setspace,bigstrut,longtable,array,breqn,url,color}
\usepackage{todonotes}
\usepackage{hyperref,caption,subcaption,multirow,setspace,wrapfig,footmisc}
\usepackage{cancel}
\usepackage{pdfpages,lipsum,comment}
\usepackage[toc,page]{appendix}

\begin{document}


\vskip 0.25in

\newcommand{\sref}[1]{\S~\ref{#1}}
\newcommand{\nn}{\nonumber}
\newcommand{\tr}{\mathop{\rm Tr}}
\newcommand{\cM}{{\cal M}}
\newcommand{\cW}{{\cal W}}
\newcommand{\cN}{{\cal N}}
\newcommand{\cH}{{\cal H}}
\newcommand{\cK}{{\cal K}}
\newcommand{\cZ}{{\cal Z}}
\newcommand{\cO}{{\cal O}}
\newcommand{\cA}{{\cal A}}
\newcommand{\cB}{{\cal B}}
\newcommand{\cC}{{\cal C}}
\newcommand{\cD}{{\cal D}}
\newcommand{\cE}{{\cal E}}
\newcommand{\cF}{{\cal F}}
\newcommand{\cX}{{\cal X}}
\newcommand{\IA}{\mathbb{A}}
\newcommand{\IP}{\mathbb{P}}
\newcommand{\IQ}{\mathbb{Q}}
\newcommand{\IH}{\mathbb{H}}
\newcommand{\IR}{\mathbb{R}}
\newcommand{\IC}{\mathbb{C}}
\newcommand{\IF}{\mathbb{F}}
\newcommand{\IV}{\mathbb{V}}
\newcommand{\II}{\mathbb{I}}
\newcommand{\IZ}{\mathbb{Z}}
\newcommand{\re}{{\rm~Re}}
\newcommand{\im}{{\rm~Im}}
\newcommand{\FC}[1]{{\textcolor{red}{#1}}}

\newcommand{\ol}{\overline}
\newcommand{\diff}{\partial}
\newcommand{\dbar}{\ol{\partial}}
\newcommand{\tmat}[1]{{\tiny \left(\begin{matrix} #1 \end{matrix}\right)}}
\newcommand{\mat}[1]{\left(\begin{matrix} #1 \end{matrix}\right)}

\newtheorem*{theorem}{\bf THEOREM}
\def\thetheorem{\thesection.\arabic{theorem}}
\newtheorem{conjecture}{\bf CONJECTURE}
\newtheorem{proposition}{\bf PROPOSITION}
\def\thetheorem{\thesection.\arabic{proposition}}
\newtheorem{observation}{\bf OBSERVATION}
\def\thetheorem{\thesection.\arabic{observation}}

\def\theequation{\thesection.\arabic{equation}}
\newcommand{\setall}{\setcounter{equation}{0}
        \setcounter{theorem}{0}}
\newcommand{\setequation}{\setcounter{equation}{0}}
\renewcommand{\thefootnote}{\arabic{footnote}}

~\\
\vskip 1cm

\begin{center}
\Large \textbf{Learning 3-Manifold Triangulations:} \\ \large Pachner Graph Network Analysis \\ \& Machine Learning IsoSigs
\\
\medskip
Francesco Costantino$^{1}$,
Yang-Hui He$^{2,3,4,5}$,\\
Elli Heyes$^{3,2}$,
\renewcommand{\thefootnote}{\fnsymbol{footnote}}
Edward Hirst\footnotetext[1]{corresponding author}$^{*6}$
\renewcommand{\thefootnote}{\arabic{footnote}}

\renewcommand{\arraystretch}{0.5} 
{\small
{\it
\begin{tabular}{rl}
  ${}^{1}$ &
  Institut de Math\'ematiques de Toulouse, 118 route de Narbonne, F-31062 Toulouse, France\\
  ${}^{2}$ &
  London Institute for Mathematical Sciences, Royal Institution, London W1S 4BS, UK\\
  ${}^{3}$ &
  Department of Mathematics, City, University of London, EC1V 0HB, UK\\
  ${}^{4}$ &
  Merton College, University of Oxford, OX14JD, UK\\
  ${}^{5}$ &
  School of Physics, NanKai University, Tianjin, 300071, P.R.~China\\
  ${}^{6}$ &
  Centre for Theoretical Physics, Queen Mary University of London, E1 4NS, UK
\end{tabular}
}
~\\
~\\
~\\
francesco.costantino@math.univ-toulouse.fr; \ 
hey@maths.ox.ac.uk; \\
elli.heyes@city.ac.uk; \ 
e.hirst@qmul.ac.uk \\
\footnotesize{Report Number: QMUL-PH-23-16}
}
\end{center}

\renewcommand{\arraystretch}{1.5} 

\vspace{10mm}

\begin{abstract}
Real 3-manifold triangulations can be uniquely represented by isomorphism signatures. Databases of these isomorphism signatures are generated for a variety of 3-manifolds and knot complements, using SnapPy and Regina, then these language-like inputs are used to train various machine learning architectures to differentiate the manifolds, as well as their Dehn surgeries, via their triangulations. Gradient saliency analysis then extracts key parts of this language-like encoding scheme from the trained models. The isomorphism signature databases are taken from the 3-manifolds' Pachner graphs, which are also generated in bulk for some selected manifolds of focus and for the subset of the SnapPy orientable cusped census with $<8$ initial tetrahedra. These Pachner graphs are further analysed through the lens of network science to identify new structure in the triangulation representation; in particular for the hyperbolic case, a relation between the length of the shortest geodesic (systole) and the size of the Pachner graph's ball is observed. 
\end{abstract}

\newpage
\tableofcontents

\section{Introduction}

Low-dimensional topology, as the study of manifolds of dimension $2,3$, and $4$, has yielded countless invaluable results and insights. Initially it was driven by the quest of the Poincar\'e conjecture (stating that a $3$-manifold homotopy equivalent to the $3$-sphere is homeomorphic to it) and W. Thurston's geometrisation conjecture \cite{thurston1982three} which were finally proven by G. Perelman through Ricci flow regularisation techniques \cite{Perelman:2006un,Perelman:2006up}. But during the eighties the discovery of the first relevant examples of topological quantum field theories by Witten \cite{Witten:1988hf} established the field as a central testing ground for most physical theories. 
Therefore, beyond their intrinsic interest within mathematics, this field has long had large influence within physics \cite{TopologyPhysics}, and more recently, its influence has been cemented in the development of the younger field of computer science \cite{Zomorodian_2005}, where in particular concepts such as triangulation and meshing are essential for rendering and visualisation.

As computational power has exponentially developed in recent decades, so too has the potential for implementation of resource-expensive computational methods. 
Many problems completely intractable before have now become rather accessible, largely through methods of `big data' and `big compute'.
However mathematical research presents a unique challenge, where datasets can often be generated infinitely and sampled cheaply, exhaustive computation becomes infeasible, demanding statistical practises. 
It is here, from a presently prominent subfield of computational science that a wealth of new techniques have been made available for approaching mathematical research; that field being: \textit{machine learning} (ML).

Machine learning covers a broad range of techniques for computational statistics, and the subfield itself can be split into three core subsubfields, which each have already seen compelling successes in mathematics and mathematical physics. 
These are: supervised learning \cite{Brodie:2019dfx,Douglas:2020hpv,Krippendorf:2022hzj,aslan2023group,Berman:2023rqb,Halverson:2023ndu,Alawadhi:2023gxa,Butbaia:2024tje,He:2020eva,Bao:2020nbi,Bao:2021auj,Bao:2021olg,Bao:2021ofk,Berman:2021mcw,Arias-Tamargo:2022qgb,Dechant:2022ccf,Aggarwal:2023swe,Choi:2023rqg,He:2023wwt,Hirst:2023kdl} (introduced in \cite{nasteski2017overview}), unsupervised learning \cite{Seong:2023njx,Chen:2023whk} (introduced in \cite{ghahramani2003unsupervised}), and reinforcement learning \cite{Halverson:2019tkf,Abel:2021rrj,Cole:2021nnt,Niarchos:2023lot,Berglund:2023ztk,Gukov:2024buj} (introduced in \cite{ernst2024introduction}); with notable introduction to string theory in \cite{He:2017aed,Carifio:2017bov,Krefl:2017yox,Ruehle:2017mzq}.
This work focuses on techniques from the first of these subsubfields, in particular using supervised learning methods to classify and study $3$-manifolds.

The direct application of machine learning methods within low-dimensional topology is increasingly widespread. Initiated in \cite{hughes2016neural}, where neural networks (NNs) were first used to predict knot invariants, after the surprisingly accurate prediction of the hyperbolic volume of a knot complement from its Jones polynomial \cite{Jejjala:2019kio}, the study carried on in \cite{Gukov:2020qaj} where it was shown that NNs could effectively distinguish high complexity knot presentations of the trivial knot from non-trivial ones. In \cite{Craven:2020bdz, Craven:2021ckk}, notable knot invariants (the Jones polynomial, Rasmussen's $s$-invariant, and Khovanov homology) where mutually predicted by NNs with high accuracies, indicating new possible hidden relations between these invariants; and further knot invariant interrelations were identified in \cite{Davies2021}. An especially useful application of these ideas has been developed in \cite{gukov2023searching}, where the recognition of ribbon knots has been pushed beyond the bounds of human capabilities through these techniques. 
Besides knots, ML methods have been used for the study of a particular kind of closed $3$-manifolds \cite{putrov2023graph}, where graph $3$-manifolds ($3$-manifolds whose Jaco-Shalen-Johansson decomposition contain no hyperbolic piece) where studied via their presentation through surgeries on plumbing graphs and various ML algorithms where implemented to reduce presentations of such $3$-manifolds to easier ones via so-called Neumann's moves. 

In this paper a different approach to representation of $3$-manifolds is considered, and one that has many direct practical applications: \textit{triangulations}. 
As proved by Casler \cite{Casler}, every $3$-manifold admits a triangulation, i.e. a decomposition into $3$-dimensional simplexes (tetrahedra). By compactness only a finite number of these is sufficient and the (finite) datum of how to glue them can be summarised most succinctly in a single string of characters: the isomorphism signature (IsoSig). IsoSigs, defined in \cite{Burton_2007} and reexplained in Appendix \ref{sec:isosigencoding}, are especially compact encodings of $3$-manifold triangulations, making bulk storage and study of triangulations far more accessible. Our first interest in this paper is to examine how amenable ML methods, and in particular NNs, are to extracting the highly-convoluted information about the triangulations from this highly-compressed representation.

Eight $3$-manifolds were selected for focus in this study, chosen due to their importance within the field and variety of properties: $S^3, S^2\times S^1$, $\mathbb{RP}^3$, $L(7,1)$ and $L(7,2)$ (homotopy equivalent but not homeomorphic lens spaces), the three torus $T^3$,  the Poincar\'e homology sphere $PHS$, and the `smallest closed hyperbolic' $3$-manifold with lowest volume $H_{SC}$ (a.k.a. the Weeks manifold). After producing thousands of different triangulations of each of these manifolds we trained NNs (and transformers) to distinguish them pairwise, observing some striking performances, which gradient saliency techniques allowed interpretation of in terms of the IsoSig components (see §\ref{sec:diffman}).

A single $3$-manifold can be encoded via infinitely many distinct triangulations, where each IsoSig completely and uniquely encodes one of these triangulations. It is a classical theorem \cite{Lickorish} that two such triangulations of the same $3$-manifold can always be connected to each other by a finite sequence of local modifications, known as \textit{Pachner moves}, which are of two types (the $2-3$ move with its inverse, and the $1-4$ move with its inverse). 
Therefore one can define an abstract graph associated to each $3$-manifold, known as the \textit{Pachner graph}, whose nodes are triangulations and edges are associated to Pachner moves between them. 

In this work, Pachner graphs were generated up to depth $5$ for the eight $3$-manifolds of focus, starting from a minimal triangulation.
As well as these Pachner graphs providing the IsoSig data used in ML, their graph structure is of its own independent interest.
Through a range of network analysis techniques these Pachner graphs were comparatively studied, with results in §\ref{sec:Pachner}.
Techniques related to network clustering, shortest paths, centrality, and cycles were employed; notably identifying $L(7,2)$'s most central triangulation was not the initial (minimal) one but one which was at distance two from it, containing $9$ tetrahedra. Also, the Weeks manifold (and other hyperbolic manifolds not of focus here) exhibited significantly faster Pachner graph growth with depth, reaching computational limits and leading us to conjecture that Pachner graphs of closed hyperbolic $3$-manifolds grow much faster than those of non-hyperbolic $3$-manifolds.

Further investigations extended network analysis to the Pachner graphs of the \texttt{snappy} \cite{SnapPy} census of cusped orientable hyperbolic $3$-manifolds\footnote{A word of caution here is needed: cusped hyperbolic $3$-manifolds are encoded in \texttt{snappy} via their \textit{ideal} triangulations, namely the (non-compact) manifolds are homeomorphic to the complement of the vertices in the triangulations encoded by the given IsoSigs. So our Pachner graphs describe such ideal triangulations for these cusped hyperbolic $3$-manifolds.}, focusing on those with up to $7$ tetrahedra (there are 4815 of them), see §\ref{sec:PachnerCusped}.
Our analyses evidenced an unexpected behavior which led us to Conjecture \ref{conj:systole} that the growth of the Pachner graph for these manifolds is inversely related to the length of the shortest closed geodesic (a.k.a. the systole). 
Additionally in §\ref{sec:knots}, knot complements where learnt via their cusped triangulation IsoSigs; focusing on eight knot complement $3$-manifolds: the unknot, the trefoil knot, the figure eight knot, the knots $5_1,5_2,6_1$, and $8_2$. The results are comparable to those for closed $3$-manifolds and in interestingly show no learning in distinguishing the figure eight knot from the knot $5_2$, which are the two knots with smallest hyperbolic volume.
We also supplemented the analysis for these knot complements by a similar analysis for their $0$ and $1$-surgeries: in §\ref{sub:dehnsurg} we showed how the unknot and the trefoil knot are very accurately distinguished by these surgeries.

Computational work was completed in \texttt{python}, with use of the \texttt{regina} \cite{regina} and \texttt{snappy} \cite{SnapPy} libraries for manipulating IsoSigs, as inspired by work in \cite{Burton_2007,Burton_2011}.
Pachner graphs were constructed and analysed using \texttt{networkx} \cite{Hagberg2008ExploringNS}, whilst machine learning used the \texttt{sci-kit learn} \cite{sklearn} and \texttt{tensorflow} \cite{tensorflow2015whitepaper} libraries.
Coding scripts and data are made available at this work's respective \href{https://github.com/edhirst/IsoSigPGML.git}{GitHub}\footnote{\url{https://github.com/edhirst/IsoSigPGML.git}} repository, including an introductory notebook with useful functions for ones own further investigations manipulating and visualising Pachner graphs.

\section{Preliminaries}
\subsection{Triangulation IsoSigs}
Any compact connected $3$-manifold $M$ with boundary can be realised by gluing a finite set of $T$ tetrahedra ($3$-simplices) along their faces via affine maps and then by deleting a small open neighborhood of the vertices of the tetrahedra. Such a structure is called an ideal triangulation and it induces a triangulation of $\partial M$ by means of $4T$ triangles (as each tetrahedron contributes $4$ triangles to the triangulation of $\partial M$)\footnote{Since we allow self-gluing of tetrahedra, in some literature this is referred to as a `$\Delta$-triangulation'.}.
In the special case of a closed manifold one can consider triangulations of $M$ minus a finite set of disjoint open balls: indeed one can recover $M$ from $M\setminus B^3$ in a unique way by gluing back a $B^3$ along its boundary $S^2$.

Each triangulation can be encoded by numbering the tetrahedra along with their vertices, indicating which pairs of faces of two tetrahedra of the triangulation are glued and how (there are $6$ possible identifications for each pair of faces). 
All these data can be encoded in a single string of characters defined in \cite{Burton_2011}, called the isomorphism signature ({\bf IsoSig}) of the triangulation.

\begin{figure}[!t]
    \centering
    \begin{subfigure}{0.45\textwidth}
        \centering
        \includegraphics[width=0.89\textwidth]{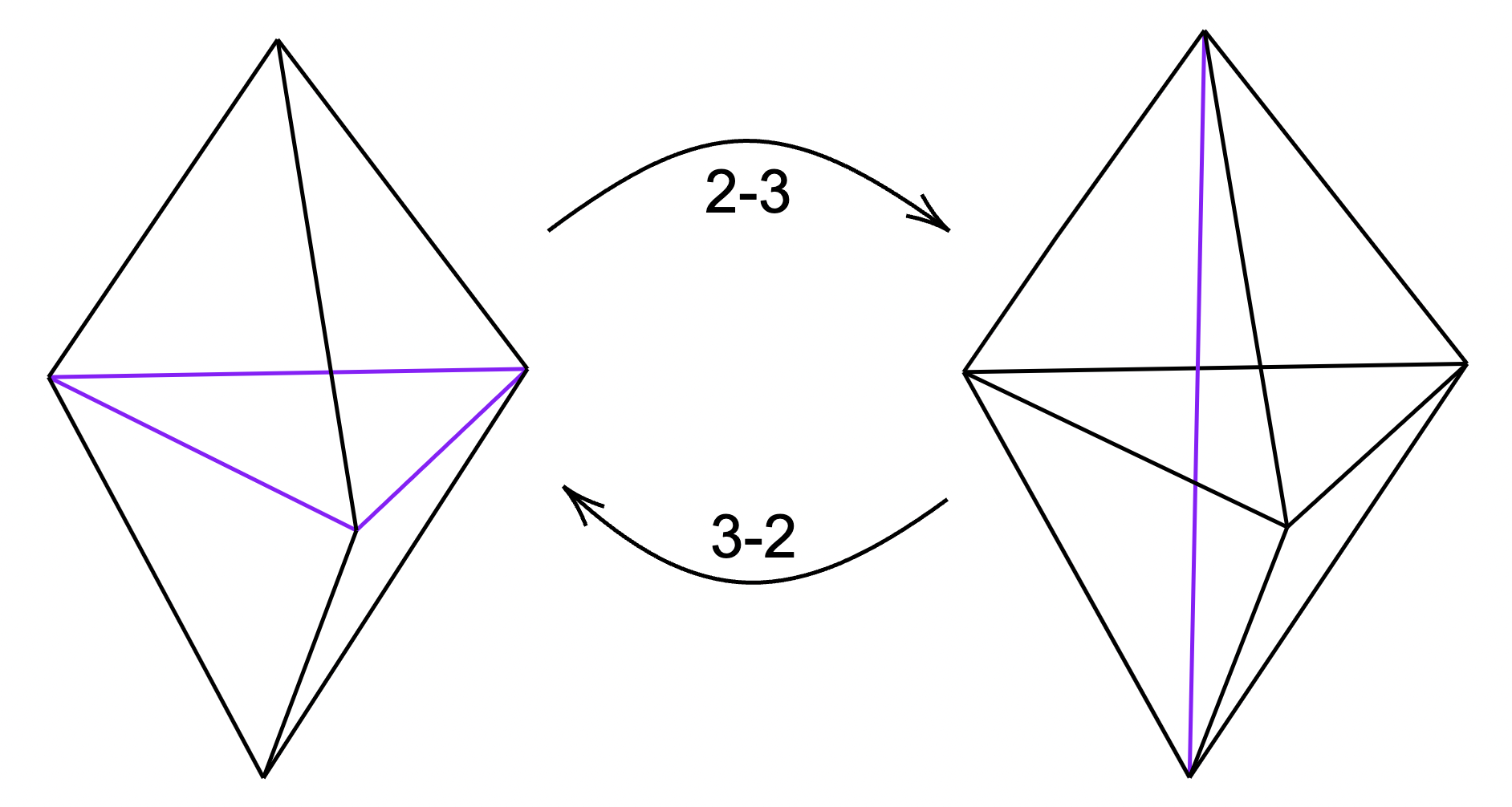}
        \caption{}
    \end{subfigure} 
    \begin{subfigure}{0.45\textwidth}
        \centering
        \includegraphics[width=0.99\textwidth]{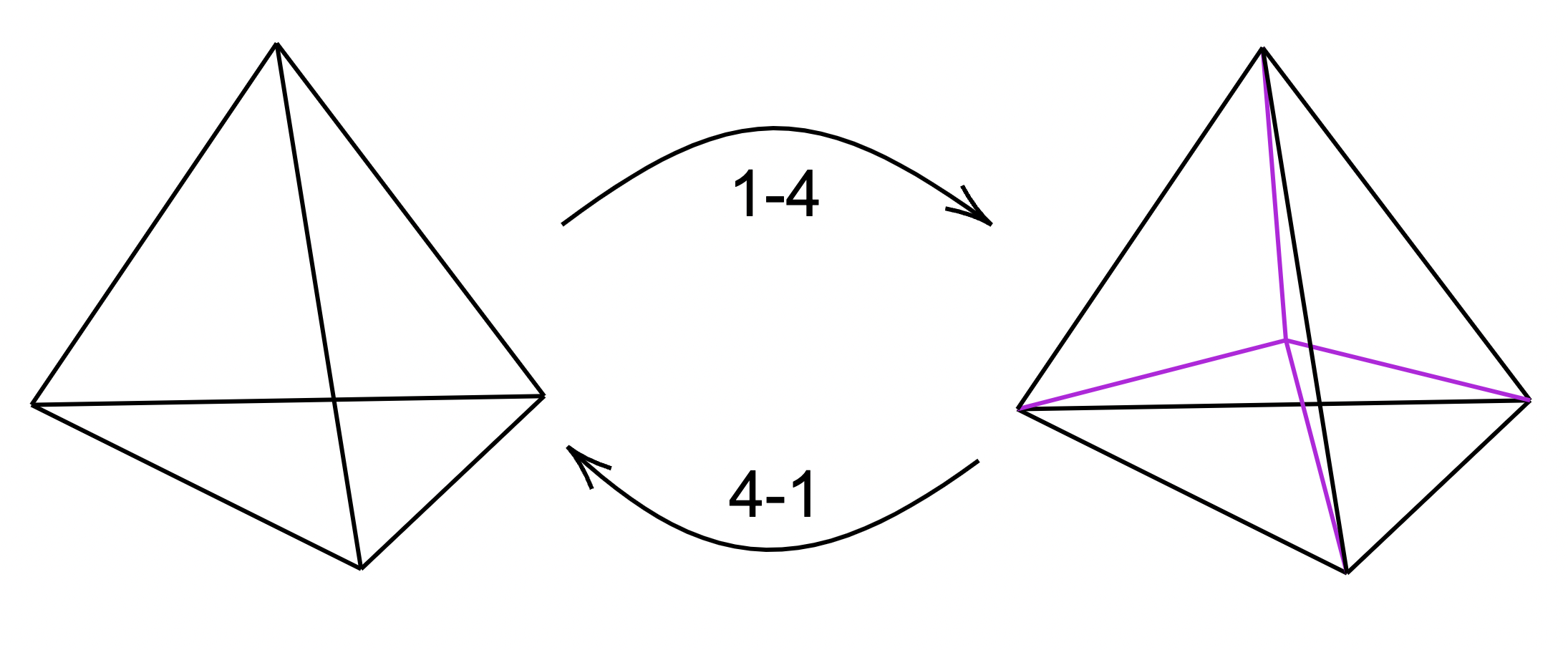}
        \caption{}
    \end{subfigure}
    \caption{Diagrammatics of the four Pachner moves, including (a) the $2-3$ Pachner move with its inverse $3-2$, and (b) the $1-4$ Pachner move, with its inverse $4-1$. For the $2-3$ move the coloured lines bound a face which is replaced by an edge which then bounds the three tetrahedra inside. For the $1-4$ move a vertex (coloured) is introduced into the centre and then connected to the remaining vertices, splitting the single tetrahedra into four.}\label{PachnerDiagrams}
\end{figure}

IsoSigs, are an exceptionally compact representation method for 3-manifold triangulations, introduced in detail with examples in Appendix \ref{sec:isosigencoding}, they are formed of a sequence of letters from the alphabet \{\texttt{a,b,...,z,A,B,...Z,0,1,...,9,+,-}\}, encoding numbers in base 64 which enumerate in turn: the number of tetrahedra in the triangulation, the destination sequence of face gluings, the type sequence listing boundaries, and the permutation sequence dictating the vertex matchings of the face gluings.
The encoding relies on many tricks for removing redundancies to make the IsoSig encoding especially efficient compared to traditional methods of listing each of these components separately.
Importantly for the IsoSig representation, the following holds: 
\begin{theorem}[Theorem 6 in \cite{Burton_2011}]
If two IsoSig strings $S_1$ and $S_2$ encode two triangulations then the triangulations are combinatorially isomorphic iff $S_1=S_2$.
\end{theorem}
Where two triangulations are isomorphic if they are combinatorically equivalent, such that some relabelling of tetrahedra, faces, edges, and vertices makes them the same.
It should be noted though, that a $3$-manifold admits infinitely many different (non-isomorphic) triangulations, which are all related to each other via the so-called bistellar or Pachner moves \cite{Lickorish}.
When applying a {\bf Pachner} move, in general the IsoSig of a triangulation changes (for most moves even the number of tetrahedra changes). 
One can then build a graph $\mathcal{P}$ whose nodes are the isomorphism types of triangulations of a fixed manifold and whose edges represent that the two triangulations can be transformed into one another via a possible Pachner move\footnote{Note that there may be (likely rare) scenarios where two different moves can transform an input triangulation into the same output triangulation, in these cases we still only connect the respective nodes by a single edge.}. 
By the above cited fundamental result this graph is connected \cite{Lickorish}. 
Before proceeding we remind the reader that these moves are of four types: $2-3$, $3-2$, $1-4$ and $4-1$ moves where for instance a $2-3$ move consists in exchanging two tetrahedra with a suitable combination of three tetrahedra and a $1-4$-move consists in subdividing a tetrahedron into four tetrahedra by coning from its barycenter to its vertices. The moves $3-2$ and $4-1$ are the inverses of respectively the $2-3$ and the $1-4$-moves. Note that the number of vertices of the triangulation does not change under $2-3$ and $3-2$ moves, but does change under $1-4$ and $4-1$ moves, exactly by one unit. Diagrammatics for these moves are shown in Figure \ref{PachnerDiagrams}.

The `Pachner graph' $\mathcal{P}$ described above is naturally endowed with two functions which we will denote $v$ and $t$ respectively given by the number of vertices and of tetrahedra of a triangulation. We will denote $\mathcal{P}_n$ the full subgraph of $\mathcal{P}$ of triangulations whose number of vertices is $n$: a particularly interesting subgraph is $\mathcal{P}_1$. At first sight one might be puzzled by this but we remind the reader that in our triangulations the tetrahedra are not embedded and self-identifications are possible. In fact it was proved independently by Matveev \cite{Matveev_spines} and Piergallini \cite{Piergallini} that for every 3-manifold the subgraph of $\mathcal{P}_1$ formed by nodes corresponding to triangulations with at least $2$ tetrahedra is non-empty and connected \cite{Matveev_spines,Piergallini}; there are only $3$ triangulations with only 1 tetrahedra and they correspond to manifolds $L(3,1),L(5,1)$ and $S^3$, so excluding them is not a problem.  

\subsection{Existing Results Surrounding Pachner Graph Paths}\label{algorithms}
In this subsection we review some main results concerning the existence of algorithms for distinguishing $3$-manifolds and their complexity. Our account is by no means complete nor exhaustive and we redirect the reader to the excellent survey by M. Lackenby \cite{Lackenby} and to S. Matveev's book \cite{Matveev}.
The rough overall message is that although algorithms exist for basically all interesting tasks in $3$-dimensional topology, when their theoretical complexity bounds are known, they are huge. Still, extensive experiments have been run for `small' triangulations and they contrast sharply with the theoretical bounds.

One of the main facts to keep in mind is that a theoretical algorithm exists to tell whether two triangulations belong to the same $3$-manifold or not. There are various proofs of this fact \cite{Kuperberg,Scott-Short}, but they all use the proof of the Geometrisation Conjecture by Perelman. However, the best `known' upper bound for the complexity of these algorithms, due to G.~Kuperberg, is huge and is of the form $2^{2^{\cdots^t}}$ where $t$ is the number of tetrahedra in the triangulations and the length of the tower of powers is a universal yet unknown constant \cite{Kuperberg}. 

The special case of recognising the $3$-sphere is known to be NP and has a more explicit algorithm, due to Rubinstein \cite{Rubinstein}. The number of Pachner moves needed to connect a $n$-tetrahedra triangulation of $S^3$ to the standard one has been bound by Mijatovi\'c \cite{Mijatovic} to be $6 10^6 n^2 2^{20000 n^2}$.

The fact that an algorithm exists to tell whether two triangulations belong to the same manifold implies that there is also a theoretical computable function which bounds the distance of two triangulations in the Pachner graph of a given $3$-manifold (see \cite{Lackenby} Theorem 4.3) based only on their numbers of tetrahedra. But for general  
$3$-manifolds this is purely theoretical at present\footnote{All the previous results can be largely improved when the $3$-manifolds are link exteriors in the sphere, but we will not review the main results here as in this paper we are mainly concerned with closed orientable manifolds.}.

We also mention that, although the theoretical bounds, when they exist, are huge, concrete experiments have been carried out and show quite a different behavior. 
In \cite{Burton_2011} B.~Burton showed that if one restricts to $1$-vertex triangulations with less than $10$ tetrahedra (there are about $82 10^6$ of them!) then these can be connected to simpler ones (i.e. triangulations with less tetrahedra) by a sequence of Pachner moves whose length does not exceed $17$ (of course this is possible only if one does not start from a minimal triangulation). In particular if one restricts to triangulations of $S^3$ only $9$ moves are needed. 
Furthermore in the same paper it is observed that the sequence of simplifying moves can be found so that the number of tetrahedra during the sequence never increases by more than $3$. 

Since a theoretical algorithm exists, but is likely exceptionally slow and computationally expensive, there is strong motivation for discovery of other methods for identifying equivalence of triangulations (that they represent the same 3-manifold and are in the same Pachner graph).
Existing results show promise for identifying paths to triangulations with fewer tetrahedra, but not between triangulations with general numbers of tetrahedra.
It is hence well motivated to tackle this broader general equivalence problem, where this work initiates the application of ML methods for this.

\section{Data Generation \& Analysis}\label{data}
To examine and learn the 3-manifold triangulation isomorphism via IsoSigs, we select a variety of 3-manifolds to focus on, in particular those which span a range of interesting properties.
Eight 3-manifolds were selected for these investigations, they include: the 3-sphere $S^3$; the product of the 2-sphere and the 1-sphere $S^2 \times S^1$; 3-dimensional real projective space $\mathbb{RP}^3$, the lens space $L(7,1)$, the lens space $L(7,2)$; the 3-torus $T^3$; the Poincaré homology sphere $PHS$; and the smallest closed hyperbolic manifold $H_{SC}$.
The 3-manifold abbreviations, along with the initial IsoSigs used to generate the datasets of equivalent triangulation IsoSigs are given in Table \ref{tab:manifoldisosigs}.
The IsoSig encoding is of particular interest due to its exceptionally compressed nature, and a thorough introduction to its structure is given in Appendix \ref{sec:isosigencoding}.
Of these manifolds $S^2 \times S^1$ is the only one not admitting a constant curvature isotropic complete Riemannian metric (we shall say for short a `non-isotropic manifold'), whilst the others are split amongst their curvatures, with $(S^3, \mathbb{RP}^3, L(7,1), L(7,2),PHS)$ positive curvature, $T^3$ zero curvature, and $H_{SC}$ negative curvature.
In addition they have properties as listed in Table \ref{mfld_properties}.

This selection of manifolds is somewhat arbitrary, but includes a list of the prototypical manifolds often considered in other works.
Specifically $S^3$, $T^3$, represent the heuristically `simplest' manifolds with $\{+,0\}$ curvatures, while $H_{SC}$ is selected as a representative with $-$ curvature.
In fact, this choice is the hyperbolic 3-manifold with minimum volume (given in Table \ref{mfld_properties}, contrasting to the other selected manifolds which all have 0 volume).
$S^2 \times S^1$ is a common choice for considering composite manifolds which may be considered as a direct product of lower dimensional manifolds.
$\mathbb{RP}^3$ is the simplest 3d projective space, a particular important construction with a broad range of uses.
Then, $PHS$ is the only homology 3-sphere with finite fundamental group (other than $S^3$ itself), making this an interesting comparison case for differentiating with $S^3$.
Finally the two Lens spaces were the first constructions of 3-manifolds which were not differentiable by their homology and fundamental group alone (Table \ref{mfld_properties} shows that these properties match).
Under special (more trivial) conditions of the Lens space constructions both $S^3$ and $S^2 \times S^1$ can be produced, making these four a particularly interesting and potentially difficult set to differentiate between.

\begin{table}[!t]
\centering
\begin{tabular}{|c|c|c|}
\hline
Manifold Name    & Initial IsoSig &  Number of Tetrahedra               \\ \hline
$S^3$            & \texttt{cMcabbgqs}   & 2                \\ \hline
$S^2 \times S^1$ & \texttt{cMcabbjaj}  & 2                 \\ \hline
$\mathbb{RP}^3$  & \texttt{cMcabbgqw}  & 2                 \\ \hline
$L(7,1)$         & \texttt{eLAkbcbddhhjhk}   & 4           \\ \hline
$L(7,2)$         & \texttt{cMcabbjqw}  & 2                 \\ \hline
$T^3$            & \texttt{gvLQQedfedffrwawrhh}  & 6       \\ \hline
PHS              & \texttt{fvPQcdecedekrsnrs}   & 5        \\ \hline
$H_{SC}=Weeks$         & \texttt{jLvAzQQcfeghighiiuquanobwwr} & 9 \\ \hline
\end{tabular}
\caption{The eight 3-manifolds of focus in this investigation, listed with the initial 1-vertex triangulation IsoSigs used to generate their Pachner graphs, and the respective number of tetrahedra in these initial triangulations.}
\label{tab:manifoldisosigs}
\end{table}

\begin{table}[!t]
\centering
\addtolength{\leftskip}{-1.5cm}
\addtolength{\rightskip}{-1.5cm}
\begin{tabular}{|cc|cccccccc|}
\hline
\multicolumn{2}{|c|}{\multirow{2}{*}{Property}}                                                             & \multicolumn{8}{c|}{Manifold} \\ \cline{3-10} 
\multicolumn{2}{|c|}{}                                                                                      & \multicolumn{1}{c|}{$S^3$} & \multicolumn{1}{c|}{$S^2 \times S^1$} & \multicolumn{1}{c|}{$\mathbb{RP}^3$} & \multicolumn{1}{c|}{$L(7,1)$}       & \multicolumn{1}{c|}{$L(7,2)$}       & \multicolumn{1}{c|}{$T^3$}                                                          & \multicolumn{1}{c|}{$PHS$}                                                     & \multicolumn{1}{c|}{$H_{SC}$}      \\ \hline
\multicolumn{2}{|c|}{Curvature}                                                                             & \multicolumn{1}{c|}{+}     & \multicolumn{1}{c|}{$\cdot$}                & \multicolumn{1}{c|}{+}               & \multicolumn{1}{c|}{+}              & \multicolumn{1}{c|}{+}              & \multicolumn{1}{c|}{0}                                                              & \multicolumn{1}{c|}{$+$}                                                         & \multicolumn{1}{c|}{$-$}           \\ \hline
\multicolumn{2}{|c|}{Volume}                                                                                & \multicolumn{1}{c|}{0}     & \multicolumn{1}{c|}{0}                & \multicolumn{1}{c|}{0}               & \multicolumn{1}{c|}{0}           & \multicolumn{1}{c|}{0}          & \multicolumn{1}{c|}{0}                                                              & \multicolumn{1}{c|}{0}                                                         & \multicolumn{1}{c|}{0.94}       \\ \hline
\multicolumn{2}{|c|}{Homology}                                                                              & \multicolumn{1}{c|}{0}     & \multicolumn{1}{c|}{$\mathbb{Z}$}     & \multicolumn{1}{c|}{$\mathbb{Z}/2$}  & \multicolumn{1}{c|}{$\mathbb{Z}/7$} & \multicolumn{1}{c|}{$\mathbb{Z}/7$} & \multicolumn{1}{c|}{$\mathbb{Z}+\mathbb{Z}+\mathbb{Z}$}                             & \multicolumn{1}{c|}{0}                                                         & \multicolumn{1}{c|}{$\mathbb{Z}/5+\mathbb{Z}/5$}   \\ \hline
\multicolumn{1}{|c|}{\multirow{2}{*}{{\begin{tabular}[c]{@{}c@{}}\\ Fund\\ Group\end{tabular}}}} & G & \multicolumn{1}{c|}{\footnotesize{x}}     & \multicolumn{1}{c|}{\footnotesize{a}}                & \multicolumn{1}{c|}{\footnotesize{a}}               & \multicolumn{1}{c|}{\footnotesize{a}}              & \multicolumn{1}{c|}{\footnotesize{a}}              & \multicolumn{1}{c|}{\footnotesize{a, b, c}}                                                        & \multicolumn{1}{c|}{\footnotesize{a, b}}                                                      & \multicolumn{1}{c|}{\footnotesize{a, b}}    \\ \cline{2-10} 
\multicolumn{1}{|c|}{}                                                                         & R   & \multicolumn{1}{c|}{\footnotesize{x}}     & \multicolumn{1}{c|}{\footnotesize{x}}                & \multicolumn{1}{c|}{\footnotesize{aa}}              & \multicolumn{1}{c|}{\footnotesize{aaaaaaa}}        & \multicolumn{1}{c|}{\footnotesize{aaaaaaa}}        & \multicolumn{1}{c|}{\begin{tabular}[c]{@{}c@{}}\footnotesize{bcBC}, \\ \footnotesize{abAB}, \\ \footnotesize{aCAc}\end{tabular}} & \multicolumn{1}{c|}{\begin{tabular}[c]{@{}c@{}}\footnotesize{abbaB}, \\ \footnotesize{aaaabAb}\end{tabular}} & \multicolumn{1}{c|}{\begin{tabular}[c]{@{}c@{}}\footnotesize{ababAbbAb}, \\ \footnotesize{abaBaaBab}\end{tabular}}  \\ \hline
\end{tabular}
\caption{Table of selected manifold properties. Properties include the sign of the curvature (note $S^2 \times S^1$ is non-isotropic so curvature is not applicable); Euler characteristic, hyperbolic volume, homology group, and fundamental group (given as generators and relators).}
\label{mfld_properties}
\end{table}

For each 3-manifold's initial triangulation and IsoSig, a dataset was generated of equivalent triangulations (each represented with their respective IsoSig) from the manifold's Pachner graph computed exhaustively up to as high a \textit{depth} (i.e. integer number of moves away from the initial IsoSig) as feasible within memory limits.
The Pachner graphs were generated with moves $\{2-3, 3-2\}$  (with further functionality in the code to consider $\{1-4, 4-1\}$ moves).
Preliminary investigations by these authors for this work using supervised neural networks to classify 3-manifolds from IsoSig triangulations showed promising performances with accuracies $>0.99$ for near-perfect learning.
However, these architectures relied heavily on the IsoSig length for distinguishing manifolds; where gradient saliency methods revealed that the only used input was the first padded entry indicating that change of IsoSig input length was the dominantly used feature in the classification.

\begin{figure}[!t]
    \centering    
    \includegraphics[width=0.7\textwidth]{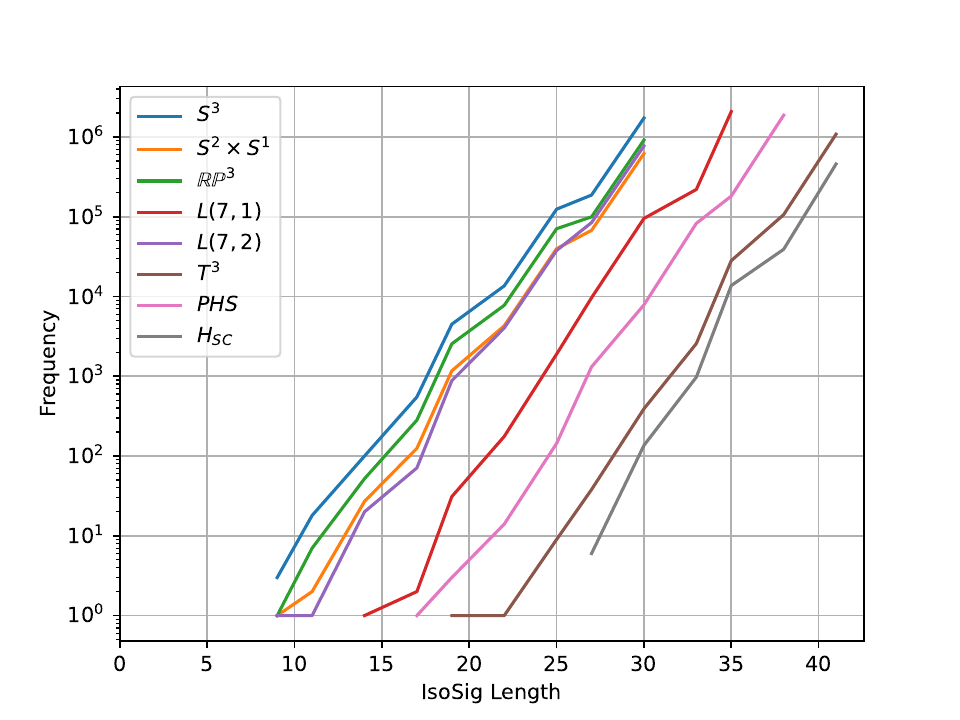}
    \caption{IsoSig length frequency distributions across the deep Pachner graphs generated for the 8 manifolds considered.}
    \label{isosiglengths}
\end{figure}

Since the IsoSig length is heavily correlated with the number of tetrahedra, when we initiate generation using triangulations with a similar (e.g. small) number of tetrahedra then the graphs will display the same IsoSig lengths throughout, however the frequency distribution of the observed lengths varies, and this is what the architectures can leverage to provide efficient estimates of the considered manifold (note the final depth often all have the same length and dominate these distributions).
For the manifolds considered the IsoSig lengths observed in the deep Pachner graph generation were: $(9, 11, 14, 17, 19, 22, 25, 27, 30,  33, 35, 38, 41)$.
The respective distributions are shown in Figure \ref{isosiglengths}.

Thus the ML focus was shifted to IsoSigs of the same length, such that the architectures could not shortcut to good performance and the true IsoSig language representation would have to be learnt.
The datasets were created by randomly sampling all those of length 30 characters (the most evenly populous across the Pachner graphs) from the generated Pachner graphs (often exhibiting millions of IsoSigs).
The sampled datasets of length 30 IsoSigs were all of size 2000, except for $T^3$ with 391 and $H_{SC}$ with 136, where this was the number of length 30 IsoSigs in the generated Pachner graphs.
Therefore investigations with these two manifolds had bias in their datasets from these different sizes, which is mitigated by considering the MCC performance measure (over accuracy) which accounts for imbalanced classes \cite{Chicco2020}.

The code scripts were written in \texttt{python}, using the \texttt{regina} library \cite{regina} for manipulating triangulations and applying Pachner moves, whilst also using the orientable cusped census hosted on \texttt{snappy} \cite{SnapPy} to source some manifold IsoSigs.
Scripts created Pachner graphs, as \texttt{networkx} objects \cite{Hagberg2008ExploringNS}, by exhaustively searching all vertices, edges, faces, and tetrahedra of each triangulation that Pachner moves could be applied to.
Where performing a move generated a new triangulation the respective IsoSig was listed as a new node in the graph, connected by an edge (labelled with the move type) to the IsoSig the move was applied to. 
Conversely, where performing a move lead to a previously generated IsoSig an edge was introduced connecting the respective nodes already in the graph (ignoring the trivial inversion of moves).
This Pachner graph generator function started from an initial input IsoSig and was adaptable to perform (input-specified) Pachner \textit{moves} up to an input-specified \textit{depth} from the initial IsoSig; outputting the Pachner graph and the list of IsoSigs, one for each triangulation of the manifold.

\subsection{Pachner Graph Network Analysis}\label{sec:Pachner}
The Pachner graphs used for IsoSig data generation are interesting combinatorial objects in their own right.
As one performs Pachner moves with the hope of simplifying a triangulation, or looking to seek some specific triangulation property, one is effectively performing a walk on this graph.
The graph network structure hence determines the potential success of these search algorithms, and more information about their network properties can guide the refinement of these search algorithms to make them as optimal as possible.

Therefore, with this in mind, we study these 1-vertex Pachner graphs $\mathcal{P}_1$ for the considered 3-manifolds, using only $\{2-3, 3-2\}$ moves, generated to as large a depth as can be performed in 240 core hours on a high-performance compute cluster \cite{apocrita}.
Then with these large 1-vertex Pachner graphs we perform extensive analysis of their network properties: degrees, clustering, shortest paths, centrality, cycles.
The distribution of number of tetrahedra in the triangulations was also tracked, since seeking a 3-manifold triangulation with the minimum number of tetrahedra is highly sought-after in the field (for simplicity in computing topological properties and constructing composite 3-manifolds).

Example Pachner graphs for the $S^3$ manifold, generated to depths 3 and 4 are shown in Figure \ref{S3_PGs}. 
The equivalent Pachner graphs for the remaining 7 considered manifolds are shown in Appendix \ref{sec:extra_PGs}.
Each of these Pachner graphs were generated from the initial IsoSig given in Table \ref{tab:manifoldisosigs}, using moves $\{2-3,3-2\}$; within these each Pachner graph node is an independent IsoSig / triangulation, with edges representing the moves connecting them (undirected since each $3-2$ move changes the IsoSig in the opposite direction to the equivalent $2-3$ move).
The Pachner graph nodes are labelled with the number of tetrahedra in the respective triangulation.
Since each of the initial IsoSigs represent 1-vertex triangulations, and only $\{2-3,3-2\}$ moves are considered which preserves the number of vertices in the triangulation, these graphs are (part of) the respective unique 1-vertex restricted Pachner graph $\mathcal{P}_1$ for each manifold \cite{Burton_2011}. 

\begin{figure}[!t]
    \centering
    \begin{subfigure}{0.47\textwidth}
        \centering
        \includegraphics[width=0.98\textwidth]{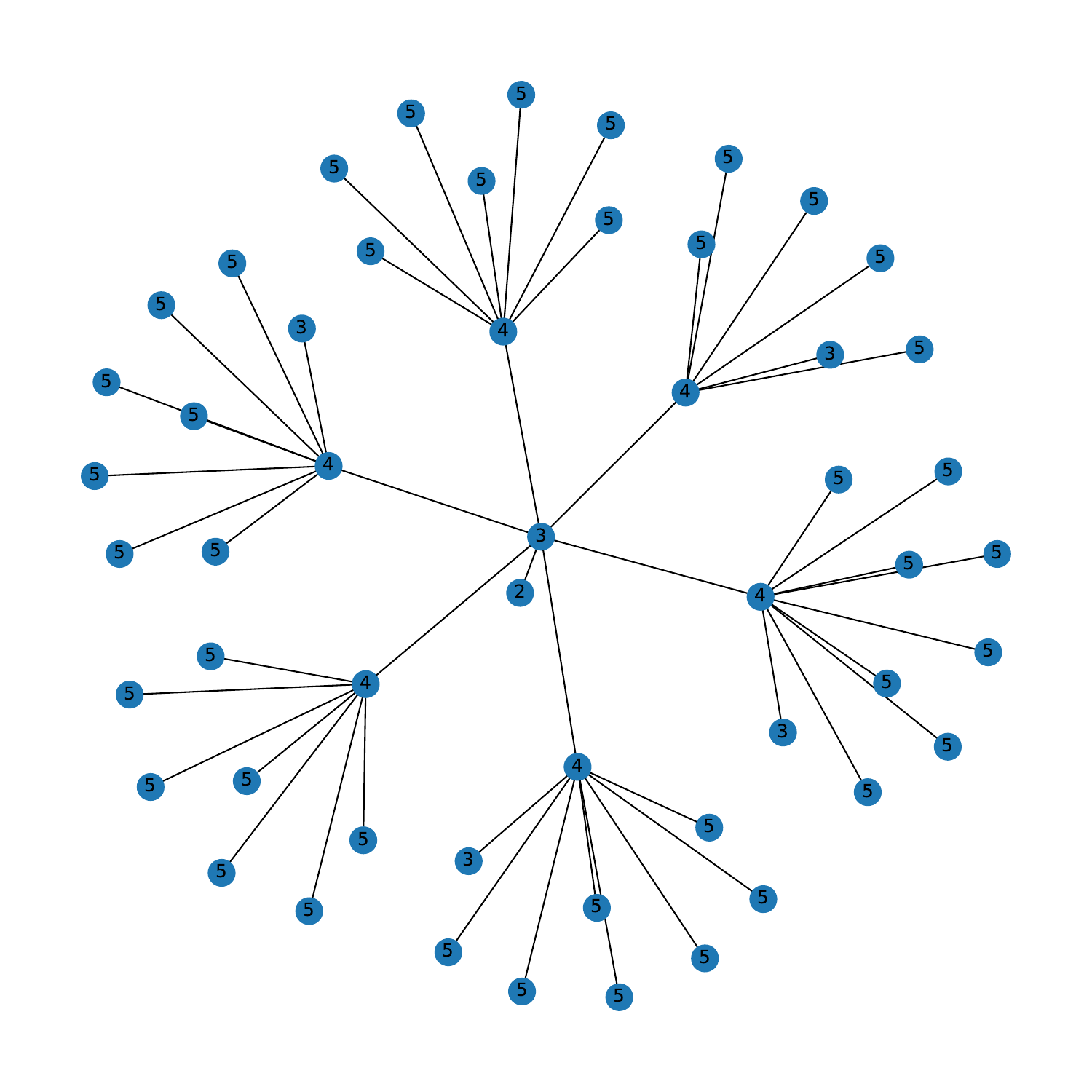}
        \caption{Depth 3}
        \label{S3_PGs3}
    \end{subfigure} 
    \begin{subfigure}{0.47\textwidth}
        \centering
        \includegraphics[width=0.98\textwidth]{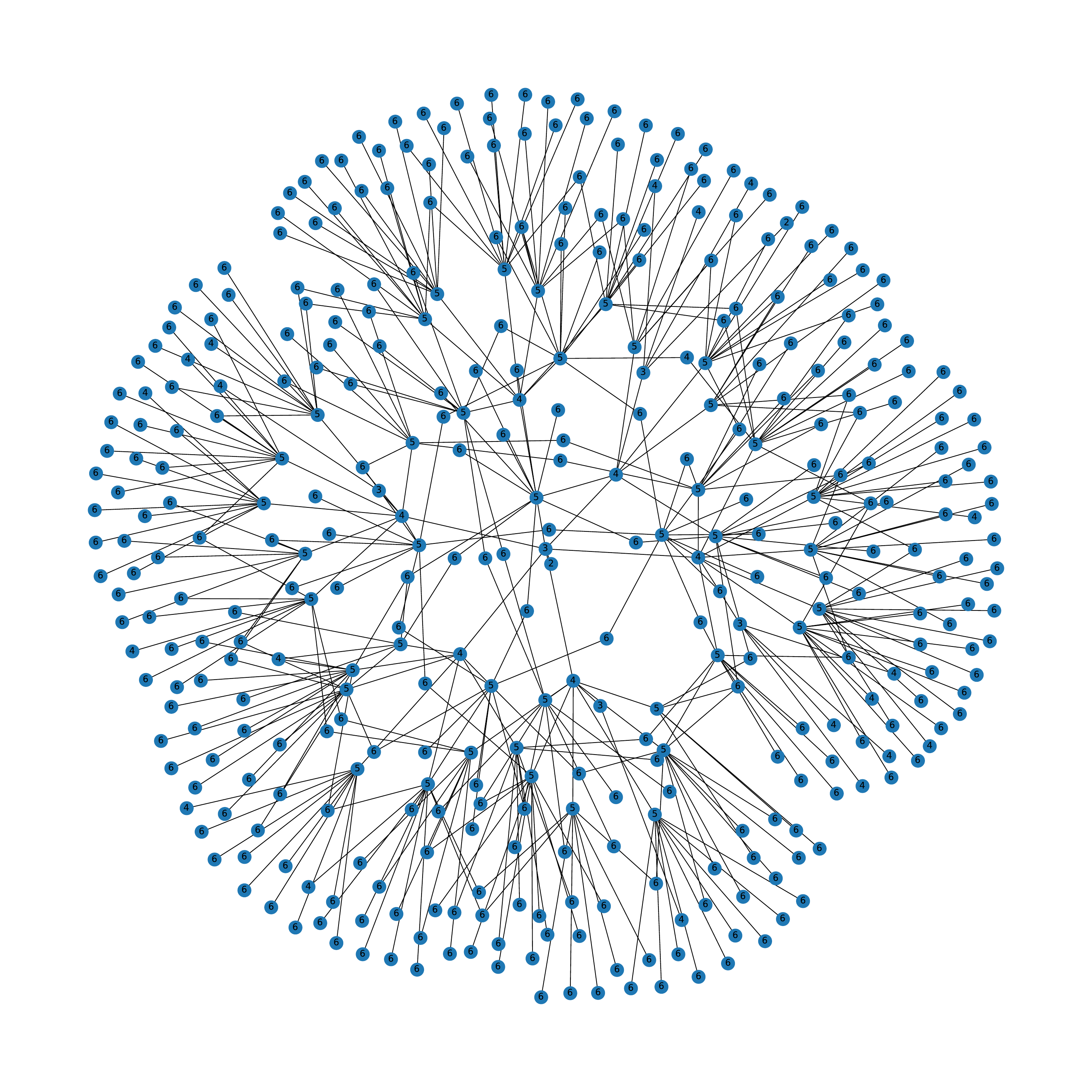}
        \caption{Depth 4}
    \end{subfigure}
    \caption{1-vertex Pachner graphs for the $S^3$ manifold, generated with $\{2-3,3-2\}$ moves up to the respective specified depths, from the initial triangulation with IsoSig: \texttt{cMcabbgqs}. The nodes are labelled with the number of tetrahedra in the respective triangulation.}\label{S3_PGs}
\end{figure}

These Pachner graph images highlight well the rapid growth rate of triangulations. 
As can be seen for $S^3$ in Figure \ref{S3_PGs}, it is only when depth 4 is considered that cycles begin to occur in this Pachner graph, at which point many cycles start to occur.
These cycles also indicate shorter paths between different triangulations, compared to returning to the initial triangulation.
The graphs also exemplify how number of tetrahedra tends to increase as one moves further from the initial triangulation in a Pachner graph.  

To probe further the graph structure we turn to tools from network science, in particular applying an analysis technique from each of its core subfields: node analysis, clustering analysis, shortest path analysis, centrality analysis, cycle analysis.
In addition, keeping track of the number of tetrahedra in each triangulation.

The Pachner graphs for the 8 selected manifolds of focus in this work were generated until memory limits deemed network analysis computationally infeasible (in each case these were RAM limits of 200GB).
Fortunately many of these Pachner graphs failed generation at the same depth of 5, making comparison between them at this depth more meaningful.
However, the final hyperbolic manifold $H_{SC}$ grew too quickly and thus failed generation earlier, this makes its comparison less useful, but it is still included for completeness.

\begin{figure}[!t]
    \centering
    \begin{subfigure}{0.47\textwidth}
        \centering
        \includegraphics[width=0.98\textwidth]{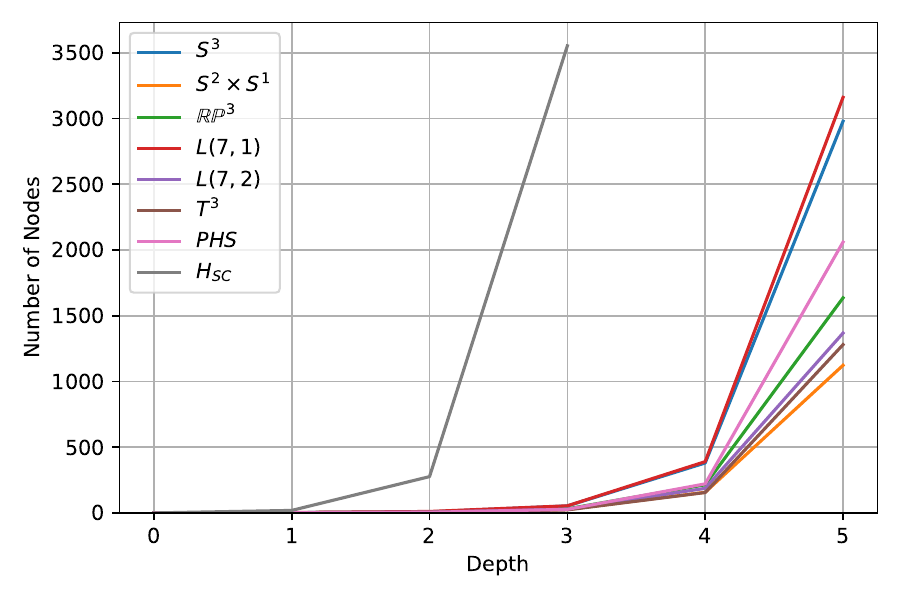}
        \caption{Number of Nodes}
    \end{subfigure} 
    \begin{subfigure}{0.47\textwidth}
        \centering
        \includegraphics[width=0.98\textwidth]{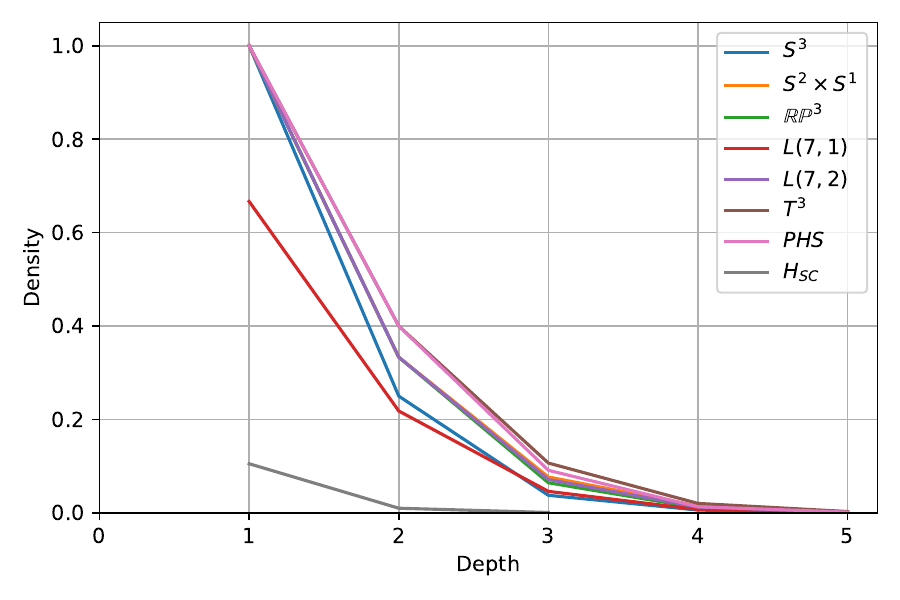}
        \caption{Density}
    \end{subfigure}
    \caption{The growth rates of the Pachner graphs with depth for the 8 selected manifolds of focus in this work. Plot (a) shows how number of nodes increases with depth, whilst (b) shows the respective decreases in graph density.}\label{Mfld_nodes}
\end{figure}

The Pachner rapid growth rate is well represented in Figure \ref{Mfld_nodes}, where the number of nodes in the Pachner graphs at each depth are shown.
This plot clearly shows the final hyperbolic manifold having a much quicker growth rate leading to memory failure at lower depth.
Generally, there will be more positions in a triangulation to apply a $2-3$ move than a $3-2$ move and so the number of tetrahedra tends to increase with depth, equivalently with more tetrahedra there are more positions to apply moves and thus this positive growth rate of triangulations in the Pachner graph appears.
Respectively the graph density drops with depth, as each time a new node is introduced to the Pachner graph with $N$ nodes, $(N-1)$ potential edges are not introduced.
As depth increases this contribution to density clearly beats any increase in density from introducing cycles, as displayed in the apparent convergence in density towards 0 with depth in Figure \ref{Mfld_nodes}.

\begin{table}[!t]
\footnotesize
\centering
\addtolength{\leftskip}{-1.75cm}
\addtolength{\rightskip}{-1.75cm}
\begin{tabular}{|c|c|c|c|c|c|c|c|}
\hline
\begin{tabular}[c]{@{}c@{}}Manifold\\ (depth)\end{tabular} & \begin{tabular}[c]{@{}c@{}}Number \\ of Nodes\end{tabular} & Density & \begin{tabular}[c]{@{}c@{}}Clustering\\ (tri, squ)\end{tabular} & \begin{tabular}[c]{@{}c@{}}Wiener Index\\ (full, norm)\end{tabular} & \begin{tabular}[c]{@{}c@{}}Centrality\\ (centre, range)\end{tabular} & \begin{tabular}[c]{@{}c@{}}Minimum Cycle Basis \\ {[[length, frequency]]} \end{tabular} & \begin{tabular}[c]{@{}c@{}} \# Tetrahedra\\ (index, min)\end{tabular} \\ \hline
$S^3$ (5)                                                  & 2979                                                       & 0.0009  & (0, 0.018)                                                      & (30734710, 6.93)                                                    & (7, 0.160)                                                           & [[4, 861], [6, 201], [8, 18]]                                                       & (0, 2)                                                                          \\ \hline
$S^2 \times S^1$ (5)                                       & 1123                                                       & 0.0023  & (0, 0.018)                                                      & (4041442, 6.41)                                                     & (11, 0.220)                                                          & [[4, 305], [6, 61], [8, 2]]                                                         & (0, 2)                                                                          \\ \hline
$\mathbb{RP}^3$ (5)                                        & 1636                                                       & 0.0017  & (0, 0.018)                                                      & (8738384, 6.53)                                                     & (13, 0.198)                                                          & [[4, 480], [6, 108], [8, 12]]                                                       & (0, 2)                                                                          \\ \hline
$L(7,1)$ (5)                                               & 3161                                                       & 0.0010  & (0, 0.022)                                                      & (33280736, 6.66)                                                    & (20, 0.201)                                                          & [[4, 1812], [6, 132], [8, 9]]                                                       & (0, 4)                                                                          \\ \hline
$L(7,2)$ (5)                                               & 1368                                                       & 0.0020  & (0, 0.019)                                                      & (5996840, 6.41)                                                     & (2, 0.202)                                                           & [[4, 417], [6, 78], [8, 9]]                                                         & (0, 2)                                                                          \\ \hline
$T^3$ (5)                                                  & 1280                                                       & 0.0030  & (0, 0.027)                                                      & (4250674, 5.19)                                                     & (13, 0.215)                                                          & [[4, 1119], [6, 22]]                                                                & (0, 6)                                                                          \\ \hline
$PHS$ (5)                                                  & 2060                                                       & 0.0017  & (0, 0.022)                                                      & (12064060, 5.69)                                                    & (6, 0.215)                                                           & [[4, 1469], [6, 49]]                                                                & (0, 5)                                                                          \\ \hline
$H_{SC}$ (3)                                               & 3553                                                       & 0.0010  & (0, 0.016)                                                      & (34673464, 5.49)                                                    & (0, 0.230)                                                           & [[4, 2489], [6, 49]]                                                                & (0, 9)                                                                          \\ \hline
\end{tabular}
\caption{Network analysis properties for the 1-vertex Pachner graphs of the 8 manifolds of focus in this work, generated to as a high depth as computationally feasible. Note that all were to depth 5, except for the last hyperbolic manifold which terminated generation at a lower depth. Analysis includes counts of the number of graph nodes, the graph density, the triangle and square clustering coefficients, the Wiener index in full and normalised form, the index of the most central node via eigenvector centrality as well as the range in centrality values across the graph, the decomposition of a minimum cycle basis, and the index of the node with the minimum number of tetrahedra and that minimum number.}
\label{tab:DeepPGNA}
\end{table}

These network analysis results are presented for the eight manifolds' 1-vertex Pachner graphs in Table \ref{tab:DeepPGNA}.
This table lists the manifolds as well as the depth they were generated to (using input initial IsoSigs from Table \ref{tab:manifoldisosigs}); each measure is discussed below with definitions. 

\noindent \textbf{Nodes:} For each Pachner graph, the \textit{number of nodes} is then listed, representing the number of unique triangulations generated. One striking fact is that the number of nodes for the depth 3 Pachner graph of the hyperbolic manifold $H_{SC}$ is bigger than the number of nodes of the depth $5$ Pachner graphs for the remaining (non-hyperbolic) manifolds. 
But for these remaining manifolds, we note the number of nodes does not seem to correlate with topological properties as the more populous graphs represent manifolds with varying curvatures, fundamental group sizes, etc; and the high number of $S^3$ nodes contradicts an argument that it is based on the number of tetrahedra in the initial IsoSig.
Respectively the graph \textit{density} is listed, where density is the ratio between the number of edges in the graph and the number of edges in the complete graph with the same number of nodes (i.e. the maximum number of possible edges), which expresses the usual approximate inverse correlation behaviour with the number of nodes.

\noindent \textbf{Clustering:} Next the graph triangle and square \textit{clustering coefficients} are listed \cite{holland1971transitivity,watts1998collective}, which represent the proportion of 3 and 4 cycles that could exist in the graph that are there; i.e. all combinations of 3 nodes could be connected into a 3-cycle, the triangle clustering coefficient is the ratio between the number of 3-cycles in the graph to the number of combinations of choosing 3 nodes.
The first observation is that all the Pachner graphs exhibit no 3-cycles, this is because any Pachner move changes the number of tetrahedra by 1, hence any 2 moves either change the number of tetrahedra by 2 or keep it the same (adding then removing, or removing then adding).
Therefore any third move following these two must change the number of tetrahedra by 1 again, such that the total change in the number of tetrahedra after 3 moves cannot be 0 and then the same as the original triangulation. 
This argument actually lifts quite nicely to any odd cycle, such that all cycles in a 1-vertex Pachner graph must be even.
The square clustering coefficients may then be non-zero, with all graphs exhibiting 4-cycles, and $T^3$ having the highest proportion.

\noindent \textbf{Shortest Paths:} The \textit{Wiener index} \cite{wiener1947structural} is the sum of the lengths of the shortest paths over all pairs of nodes\footnote{i.e. select a pair of nodes in the graph, find the shortest path between them (if there are multiple select one arbitrarily), count the number of edges in the path to get the length of this path, repeat for all pair of nodes in the graph, take the sum.}, and is hence also presented in a normalised form with respect to the number of pairs of nodes (number of nodes $N$ choose 2).
A smaller normalised Wiener index indicates the graph is easier to walk around, with a high connectivity, which is the case for the non-positive curvature manifolds with the lowest values.

\noindent \textbf{Centrality:} Centrality measures determine the most central nodes, and how skewed the connectivity is away from uniform constant degree.
The centrality measure examined in this work is \textit{eigenvector centrality}, which assigns scores to the nodes based on their respective entries in the normalised eigenvector corresponding to the unique largest eigenvalue of the graph's adjacency matrix (via the Peron-Frobenius Theorem \cite{perron1907zur,frobenius1912ueber}) \cite{bonacich_1987}.
A graph's adjacency matrix is intrinsically linked to potential paths through it, and such raising it to the $k$th power counts the number of $k$-paths between the indexed nodes.
This limiting behaviour is then dictated by the largest eigenmode, and hence the eigenvector centrality score provides a means of determining the limiting connectivity of nodes.
In Table \ref{tab:DeepPGNA}, the node index with the highest centrality score is listed, as well as the range of these scores across the graph.
Interestingly, only for $H_{SC}$ is the most central node the initial IsoSig triangulation, indicating the other Pachner graphs may be more symmetrically generated about alternative initial IsoSigs. 
In some cases, such as for $L(7,1)$ the index is as high as 20, such that it is the 20th triangulation produced when exhaustively generating the Pachner graph, and is 2 moves away from the initial triangulation; it is a $9$-tetrahedra triangulation. 
Iteratively, this measure may be used in this way to search for useful starting triangulations. 
The range of centrality scores are largest for the non-positive curvature manifolds, which then turn out to have more skewed Pachner graphs and knowledge of the most central node becomes important as a key triangulation many paths will pass through (this could be a useful information when designing paths to minimise the number of tetrahedra).

\noindent \textbf{Cycles:} Cycles in a Pachner graph represent availability of different paths to move between selected nodes.
As mentioned previously, since $3-2$ and $2-3$ moves each change the number of tetrahedra by 1, this leads to all cycles being even.
2-cycles enact the inverse of any move, so are trivial on these undirected graphs.
This causes 4-cycles to be the first cycles that can occur, whose frequency is assessed through the square clustering coefficient.
In general, a $(2k)$-cycle amounts to $k$ $2-3$ moves and $k$ $3-2$ moves and there are no a priori restrictions on their order.
Where the number of moves is small relative to the number of tetrahedra, cycles are likely caused by the more trivial case where each of the $k$ $2-3$ moves may be paired with its respective inverse $3-2$ move removing the same tetrahedra introduced, but in a different order such that the intermediary triangulations are distinct.
However as the number of moves gets larger there're less untouched tetrahedra to perform the new moves on, and then moves start to overlap in the same part of the triangulation leading to less trivial cycles.

As can be seen in the Pachner graph plots, cycles often intersect, and where two smaller cycles intersect a larger cycle can easily be created by traversing the exterior of them both.
To circumvent redundancies due to this, it is standard practice to consider a basis for the cycles \cite{liebchen2007classes}, and additionally to focus on bases where the sum of the cycle lengths is minimised: a \textit{minimum cycle basis}.
However there is still redundancy in which cycles form part of this bases, and therefore the analysis in Table \ref{tab:DeepPGNA} instead focuses on its properties, notably counting the frequencies of cycles of each length in the basis (which is unique across all minimum cycle bases).
As can be seen, all the graphs exhibit 4-cycles, some with an order of magnitude more than others.
Where the number of 4-cycles is higher the number of 6-cycles tends to be lower (except for $\mathbb{RP}^3$ which has a high number of 6-cycles, and $S^2 \times S^1$ with a surprisingly low number).
Up to depth 5 the maximum cycle size which can occur is 10, and in none of the graphs does this occur, and interestingly it is only for the manifolds with simplest fundamental group presentation that 8-cycles occur.

\noindent \textbf{Number of Tetrahedra:} The final analysis tracks the node index with the fewest number of tetrahedra, often desirable for use in manifold computation and construction, as well as listing that number of tetrahedra. 
All manifolds have index 0 as their minimum (i.e. index 0 is the first triangulation in the generation process which is the initial triangulation), expected since a minimal triangulation was selected as an initial point for the generation.
A final remark to emphasise is that the node with the minimum number of tetrahedra does not correlate with centrality, only in the case of $H_{SC}$ does the most central triangulation have the minimum number of tetrahedra.

\subsubsection{Orientable Cusped Census}\label{sec:PachnerCusped}
Beyond the eight selected 3-manifolds that we focus on in this paper, the \texttt{snappy} Orientable Closed Census provides a list of 11031 orientable closed hyperbolic manifolds. 
These however are represented as cusped manifolds which when filled have a range of tetrahedra for their initial triangulations from 9 to 37 -- too large for sensible Pachner graph analysis as the number of triangulations would grow too quickly with depth to compute to a meaningful size.

Therefore we instead consider the respective \texttt{snappy} Orientable Cusped Census, which lists all 61911 orientable cusped hyperbolic manifolds which can be represented by triangulations with at most 9 tetrahedra.
We partition off the 4815 with $< 8$ tetrahedra for an exhaustive analysis of their equivalent Pachner graphs.
Each of these manifolds are represented by an initial triangulation (and hence IsoSig) in the database.

\begin{figure}[!t]
    \centering
    \begin{subfigure}{0.47\textwidth}
        \centering
        \includegraphics[width=0.98\textwidth]{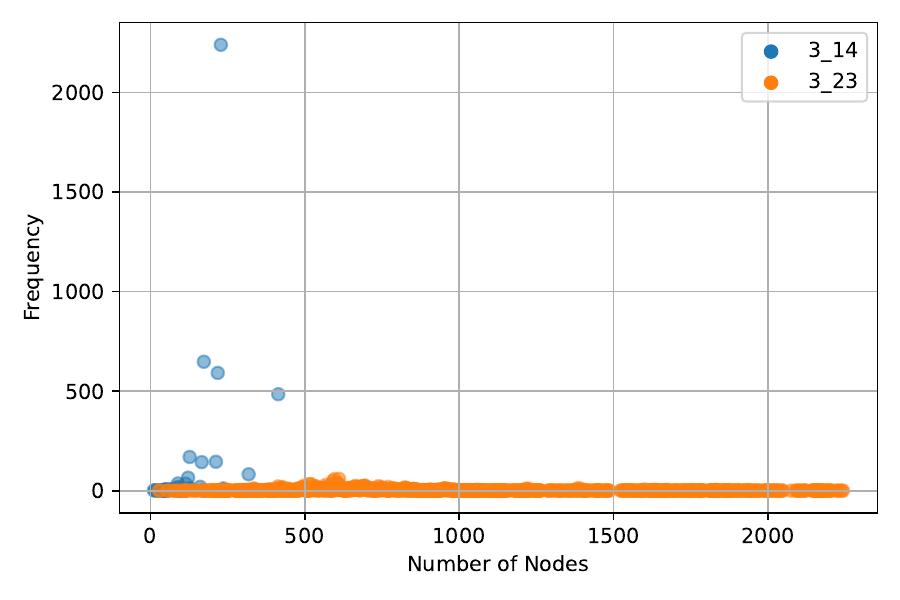}
        \caption{Nodes}\label{occ_nodes}
    \end{subfigure} 
    \begin{subfigure}{0.47\textwidth}
        \centering
        \includegraphics[width=0.98\textwidth]{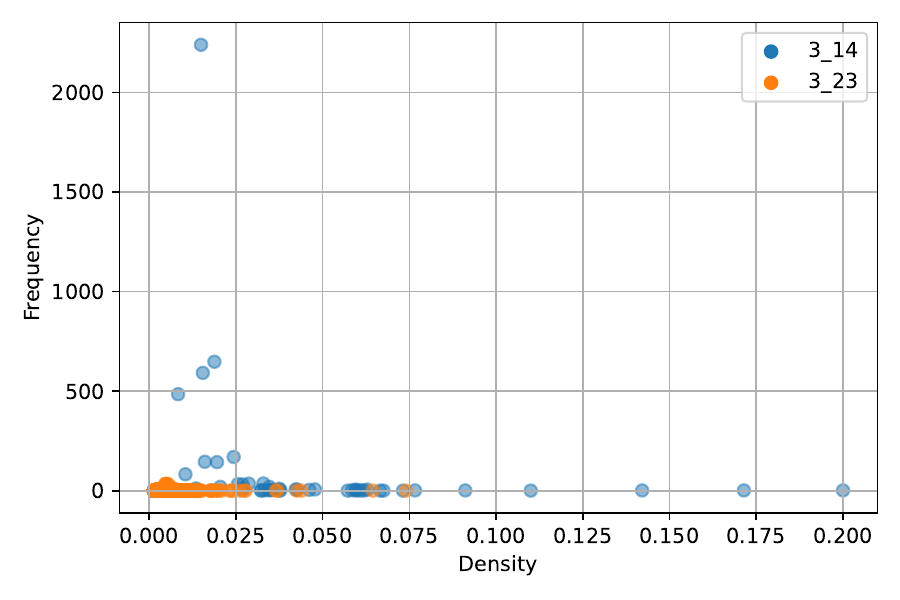}
        \caption{Density}\label{occ_density}
    \end{subfigure}\\
    \begin{subfigure}{0.47\textwidth}
        \centering
        \includegraphics[width=0.98\textwidth]{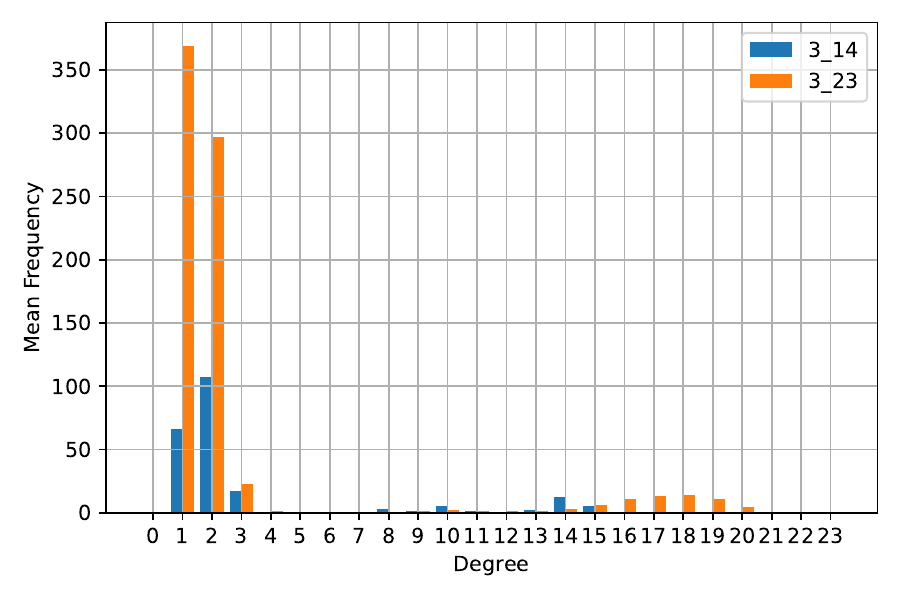}
        \caption{Mean Degree Distribution}\label{occ_degrees}
    \end{subfigure}
    \caption{Network analysis of the orientable cusped census Pachner graphs (where number of initial tetrahedra $<8$). Analysis considers building the Pachner graphs to depth 3 with either only the $\{1-4, 4-1\}$ moves (denoted 3\_14), or only the $\{2-3, 3-2\}$ (denoted 3\_23) moves, plotting histograms of: (a) number of nodes in the graph; (b) graph density; and (c) the graph degree distribution averaged across the census. Note the degrees take integer values and the vertical lines are offsetted for visibility.}\label{occ_nodesdensitydegrees}
\end{figure}

Despite the large number of manifolds in the census, with use of the high-performance compute cluster it became computationally feasible to generate Pachner graphs for \textit{all} hyperbolic manifolds in this partition of the census, and perform the equivalent network analysis.
Pachner graphs were generated for both move combinations $\{2-3, 3-2\}$ and $\{1-4, 4-1\}$, up to depth 3, which we denote `3\_23' and `3\_14' respectively. 
The subsequent bulk network analysis provides insight into the triangulation deformation structure of these Pachner moves, how they each impact the graph structure, as well as general insight into how Pachner graphs grow.
Considering the same set of network analysis techniques used to examine the depth 5 Pachner graphs for the eight manifolds of focus in Table \ref{tab:DeepPGNA}, we created plots of these properties for the census as shown in Figure \ref{occ_nodesdensitydegrees}.

\begin{figure}[tb]
\centering
    \begin{subfigure}{0.47\textwidth}
        \centering
        \includegraphics[width=0.98\textwidth] {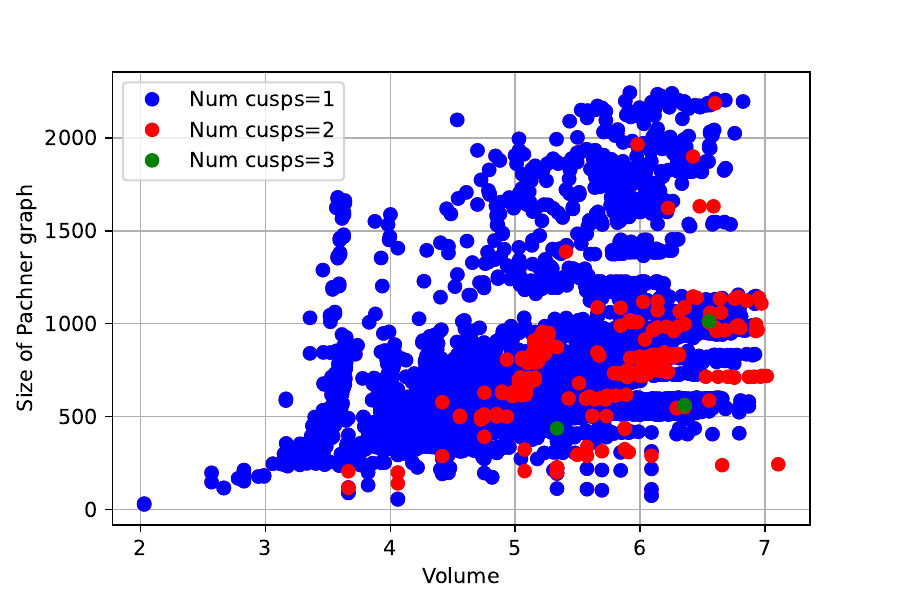}
        \caption{}
        \label{fig:PachnersizesvsVol}
    \end{subfigure}
    \begin{subfigure}{0.47\textwidth}
        \centering
        \includegraphics[width=0.98\textwidth]{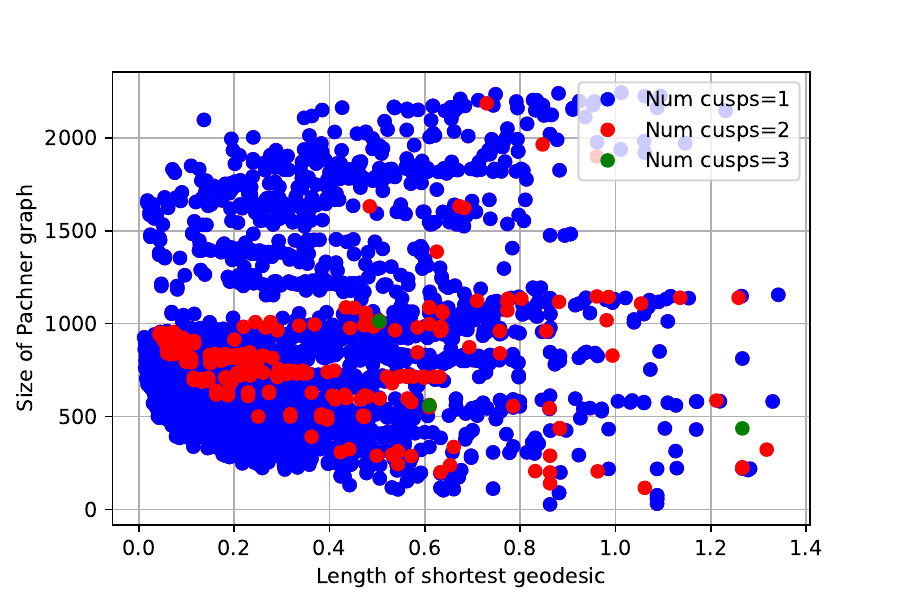}
        \caption{}
        \label{fig:PachnersizesvsShortestGeodesic}
    \end{subfigure}\\
    \begin{subfigure}{0.45\textwidth}
        \centering
        \includegraphics[width=0.98\textwidth]{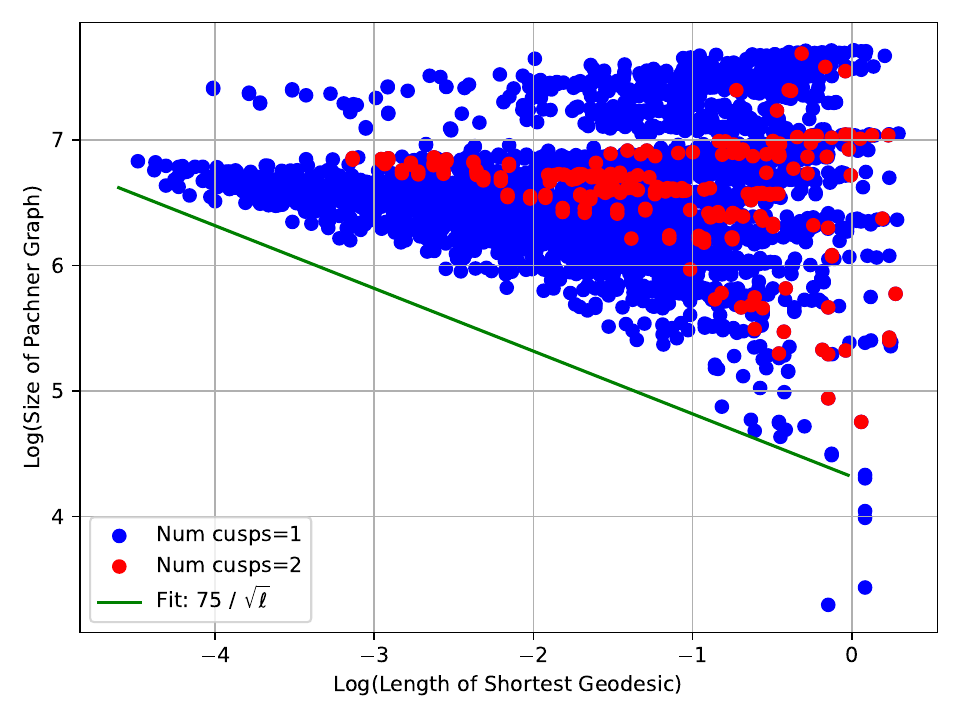}
        \caption{}
        \label{fig:logged_PachnersizesvsGeodesic}
    \end{subfigure}
    \caption{The distribution of the sizes of the Pachner graphs (i.e. number of nodes) as a function of (a) the volume and of (b) the length of the shortest geodesic (systole). In blue the points corresponding to one cusped manifolds, in red those corresponding to $2$-cusped manifolds, and in green the remaining ones. Plot (c) shows plot (b) with axes logged to inform on the power law behaviour of the lower bound, including a potential lower bound of $\frac{75}{\sqrt{\ell}}$.}
\end{figure}

\begin{figure}[tb]
\centering
    \begin{subfigure}{0.47\textwidth}
        \centering
        \includegraphics[width=0.98\textwidth]{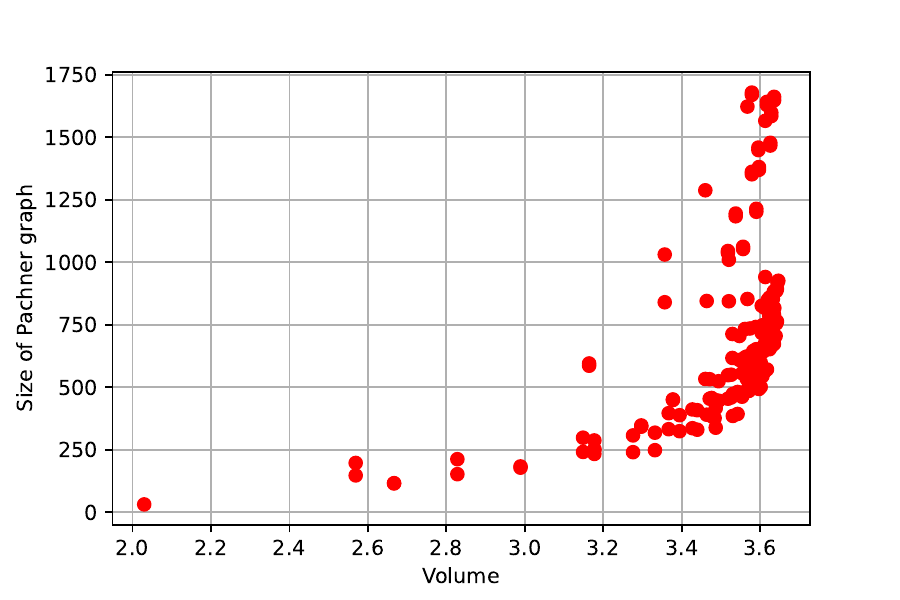}
        \caption{}\label{fig:PachnersizesvsVol_drill}
    \end{subfigure}
    \begin{subfigure}{0.47\textwidth}
        \centering
        \includegraphics[width=0.98\textwidth]{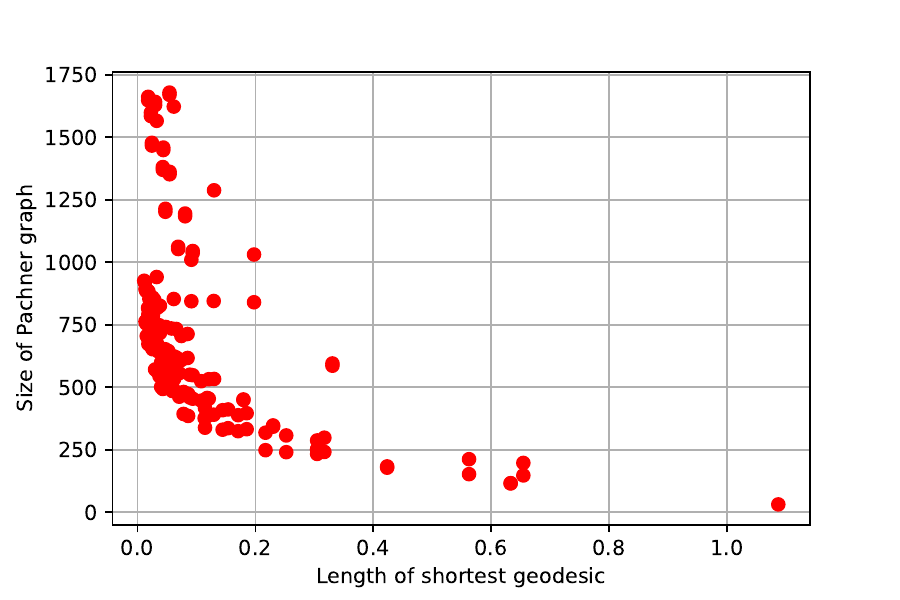}
        \caption{}\label{fig:PachnersizesvsShortestGeodesic_drill}
    \end{subfigure}
    \caption{The distribution of the sizes of the Pachner graphs (i.e. number of nodes) as a function of (a) the volume and of (b) the length of the shortest geodesic restricted to the manifolds whose drilling along the shortest geodesic is a fixed one, here ``m129" (in \texttt{snappy}'s cusped census). The symmetry is not random and is related to the known relation between the volume of a manifold and its Dehn fillings along `long' meridians.}
\end{figure}

\noindent \textbf{Nodes:} First order analysis of the graphs considers the sizes of the node and edge sets. 
A histogram for the number of nodes is shown in Figure \ref{occ_nodes}, demonstrating that the $2-3$ (and their inverse $3-2$) moves often lead to many more nodes. 
Since the number of possible $1-4$ moves for a triangulation equals the number of tetrahedra $N$, whereas the number of possible $2-3$ moves equals the number of common faces (which is half the total number of faces $4N/2 = 2N$) this matches the expected approximate scaling that $2N > N$ for triangulations of $N$ tetrahedra, causing more possible $2-3$ moves, and more nodes in the 3\_23 Pachner graphs.
The $1-4$ moves also have much higher frequencies, meaning their graph generation is more systematic, which is perhaps to be expected since the number of tetrahedra alone (which many manifolds will share) sets the number of $1-4$ moves which can be enacted, whilst the specific triangulation gluing sets the number of $2-3$ and $3-2$ moves (and the rarer $4-1$ moves).
The density in Figure \ref{occ_density} has a weak inverse relationship to the number of nodes, where adding nodes to a large graph introduces many more potential edges, and this respectively drops the density value. Since the density does drop according to this behaviour this demonstrates that the Pachner graph edge sets grow slower than the node sets.

We also analysed the size of the Pachner graphs as a function of the hyperbolic volume of the manifolds and their number of cusps (which varies from $1$ to $3$ in our dataset). The result is shown in Figure \ref{fig:PachnersizesvsVol}.
Looking at the figure one realises that there are `peaks' in the sizes of the Pachner graphs for one cusped manifolds and that these peaks correspond to the first red dots, namely the lowest volume $2$-cusped manifolds. This is no random phenomenon. The blue dot with highest size of the Pachner graph and volume $<3.63$ turns out to be ``v0117" (in \texttt{snappy}'s nomenclature, IsoSig: \texttt{hLALPkbcbefggghxwnxibk}) and drilling out its shortest geodesic one gets ``m129" (IsoSig: \texttt{eLPkbdcddhgggb}) which is the second lowest volume manifold in \texttt{snappy}'s census with number of cusps $2$ (it corresponds to one of the two leftmost red dots in the figure). This leads us to think that the size of the Pachner graph could be inversely related to the length of the shortest geodesic in the hyperbolic manifold which, by Thurston-Jorgensen's theorem on hyperbolic Dehn fillings corresponds to Dehn fillings of manifolds with one more cusp with very large filling slopes. In order to check this we plotted in Figure \ref{fig:PachnersizesvsShortestGeodesic}, the sizes of the Pachner graphs with respect to the length $\ell$ of the shortest geodesic in the cusped manifold and the plot clearly shows that such number is bounded below by an unknown decreasing function of $\ell$. 
Scaling the plot axes, as in Figure \ref{fig:logged_PachnersizesvsGeodesic}, shows the behaviour more clearly; the lower bound forms a straight line indicating a power law relationship, and the gradient is $\sim -0.5$, indicating a relationship of the form $\propto \frac{1}{\sqrt{\ell}}$ for this dataset.

Furthermore, a closer look reveals that the dots come in families, corresponding to which manifold is obtained by drilling the shortest geodesic. 
This is exemplified in Figures \ref{fig:PachnersizesvsVol_drill} and \ref{fig:PachnersizesvsShortestGeodesic_drill} where the same analysis as before is restricted to the manifolds whose drilling is a fixed one, here ``m129'' (chosen just as an example). 
All this seems to corroborate the following 

\begin{conjecture}\label{conj:systole}
The growth rate of the Pachner graph of a cusped hyperbolic $3$-manifold $M$ is bounded below by a function of the  inverse of the systole (i.e. the length of the shortest geodesic in $M$). 
\end{conjecture}

\begin{figure}[tb]
\centering
    \begin{subfigure}{0.47\textwidth}
        \centering
        \includegraphics[width=0.98\textwidth]{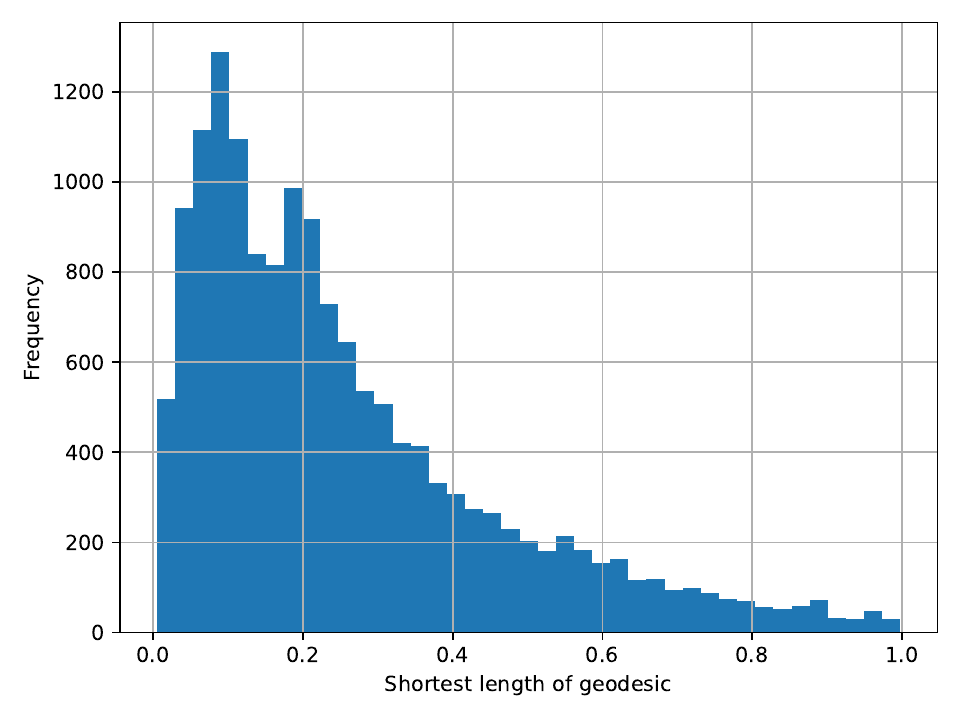}
        \caption{Histogram of Shortest Geodesic Length}\label{fig:Histogram_Shortest_geodesic}
    \end{subfigure}
    \begin{subfigure}{0.47\textwidth}
        \centering
        \includegraphics[width=0.98\textwidth]{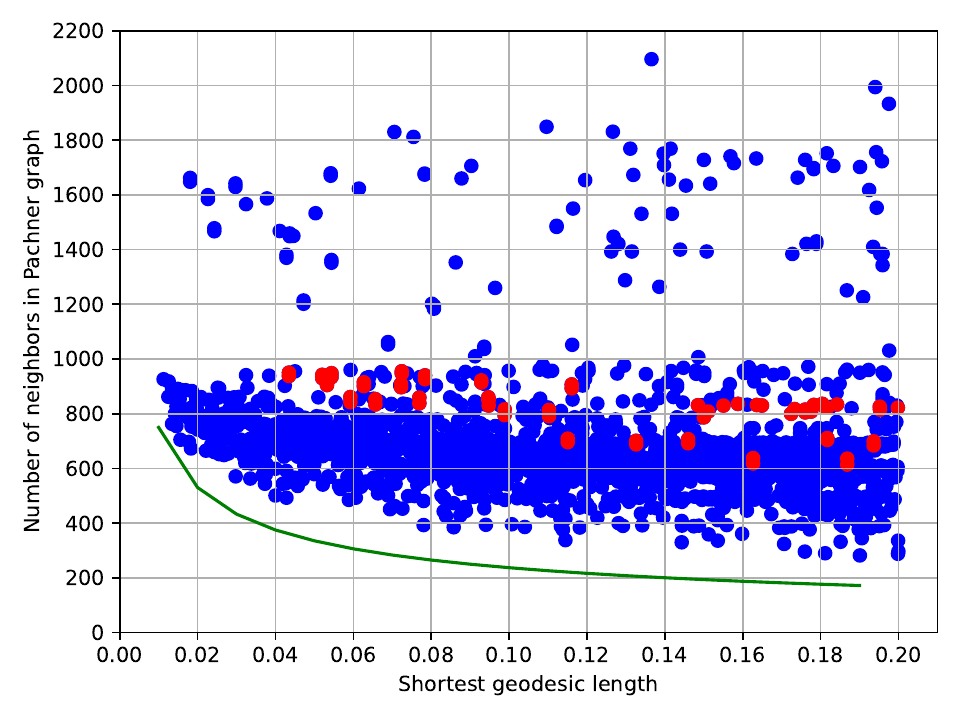}
    \caption{Zoom in on Systoles $<0.2$}\label{fig:Restricted plot of pachner size}
    \end{subfigure}
    \caption{On the left the global histogram of the shortest geodesic length on the Oriented Cusped Census of snappy with $<9$ tetrahedra. On the right the restriction of the plot of the size of the Pachner graph for manifolds with shortest geodesic less than $0.2$. In green, the graph of the function $f(\ell)=\frac{75}{\sqrt{\ell}}$ seems to approximately minor the scatter plot on the interval $[0,0.2]$, but it is plot here only as an example. 
    }
\end{figure}

One may be concerned by a potential bias in the previous analysis: as the number of tetrahedra grows the relative scarcity of triangulations with big shortest geodesic length also increases, and this may be causing the observed phenomenon. 
To check this, we first drew the histogram of the shortest geodesic lengths for the manifolds with complexity up to 8 in the Orientable Cusped Census of snappy: see Figure \ref{fig:Histogram_Shortest_geodesic}. The figure clearly shows that triangulations with shortest geodesic $<0.2$ are fairly common, thus in the right part of the same figure we restrict the image of Figure \ref{fig:PachnersizesvsShortestGeodesic_drill} to the interval $\ell\in [0,0.2]$.
The results are consistent with the observed behaviour.


\noindent \textbf{Degrees:} The node degree distribution of a graph can be described by a histogram showing the frequency of the nodes that have a certain number of edges incident to them. 
Averaging the degree distributions over the entire considered census partition produces the plot given in Figure \ref{occ_degrees} (i.e. the mean frequency of 3\_23 nodes with degree 1 across the 4815 manifolds in this census subset is 368.9).
For both move styles all graphs had at least one edge, so that a move of either form could always be performed on the initial triangulation.
This is shown by a 0 mean frequency of degree zero, so that no Pachner graphs were the trivial empty graph formed by a single node and no edges.
For both move styles the degree distributions were bimodal, with 2 peaks, indicating that some moves (in both cases) significantly open up the number of new subsequent move possibilities and thus have much higher degree in this exhaustive graph generation.
The peaks in the distributions (particularly for degree 1) are due to the leaves of the Pachner graphs created in the latter moves at the limit of the generation truncation, however they can occur otherwise as exemplified in Figure \ref{S3_PGs3}. 
The 3\_23 graphs have a larger relative size between the bimodal peaks than the 3\_14 graphs ($368.90/13.94 = 26.46 > 8.48 = 107.04/12.63$), so that the 3\_14 graphs' nodes are more skewed to higher degrees as these moves are more likely to open up many move possibilities, despite there being less possibilities demonstrated by the degree value for the second peak.
The 3\_23 graphs also have a larger difference between their two peaks ($18-1 > 14-2$), indicating that the edges are less evenly distributed between the nodes.
This is linked to the move possibility scaling, where as a $2-3$ move changes the number of tetrahedra in a triangulation from $N$ to $(N+1)$ the approximate scaling of the $2-3$ move possibilities grows at a rate of order $2N$ (i.e. the number of faces mod the number of tetrahedra per face, $4/2=2$). 
Alternatively as a $1-4$ move changes the number of tetrahedra from $N$ to $N+3$ the scaling of the $1-4$ move possibilities grows exactly at a rate of $N$.
Therefore, particularly as moves increase the number of tetrahedra, there are on average relatively more $2-3$ moves which can be made, leading to higher degree nodes\footnote{Note in this crude leading behaviour approximation we neglect contributions from the less probably inverse $3-2$ and $4-1$ moves.}.

\begin{figure}[!t]
    \centering
    \begin{subfigure}{0.47\textwidth}
        \centering
        \includegraphics[width=0.98\textwidth]{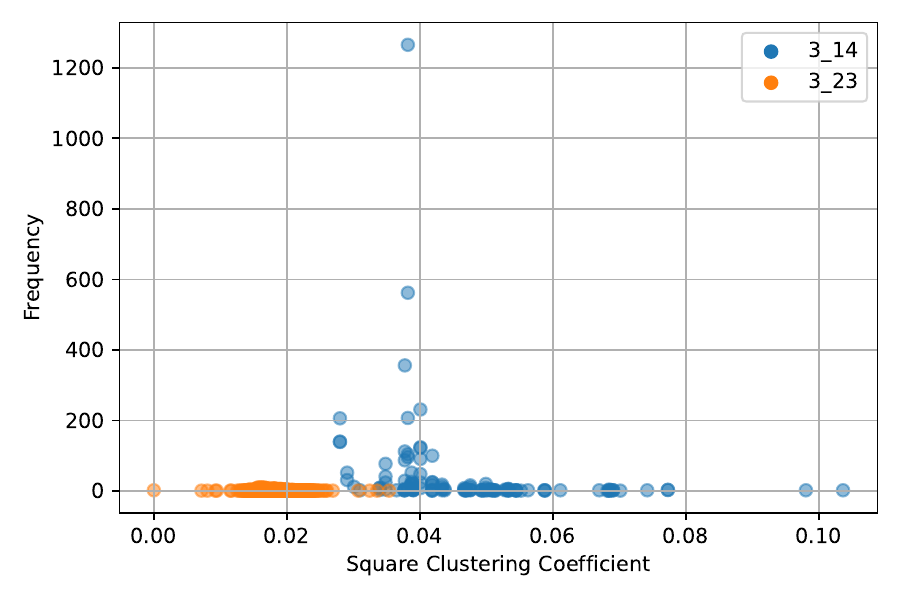}
        \caption{Clustering}\label{occ_clust}
    \end{subfigure} 
    \begin{subfigure}{0.47\textwidth}
        \centering
        \includegraphics[width=0.98\textwidth]{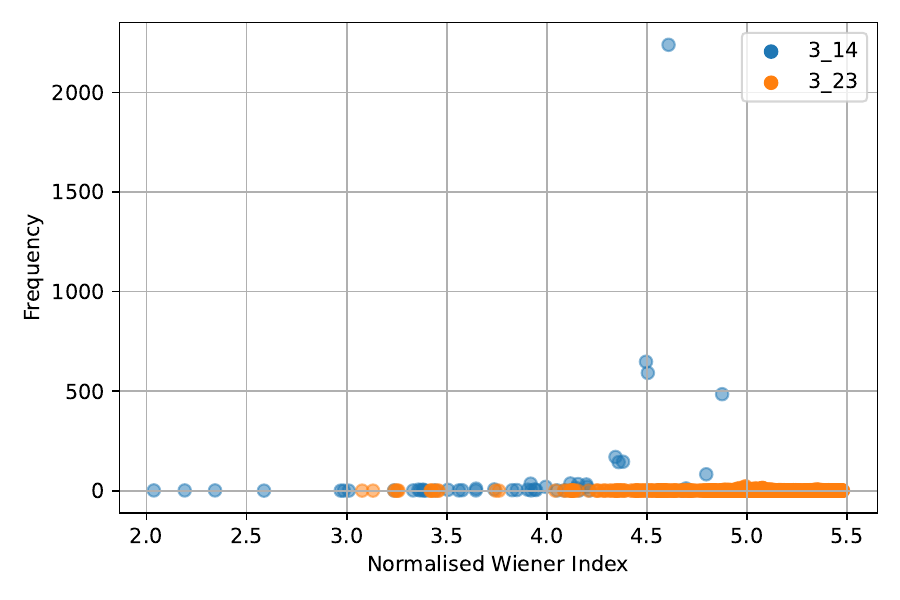}
        \caption{Shortest Path}\label{occ_wiener}
    \end{subfigure}\\
    \begin{subfigure}{0.47\textwidth}
        \centering
        \includegraphics[width=0.98\textwidth]{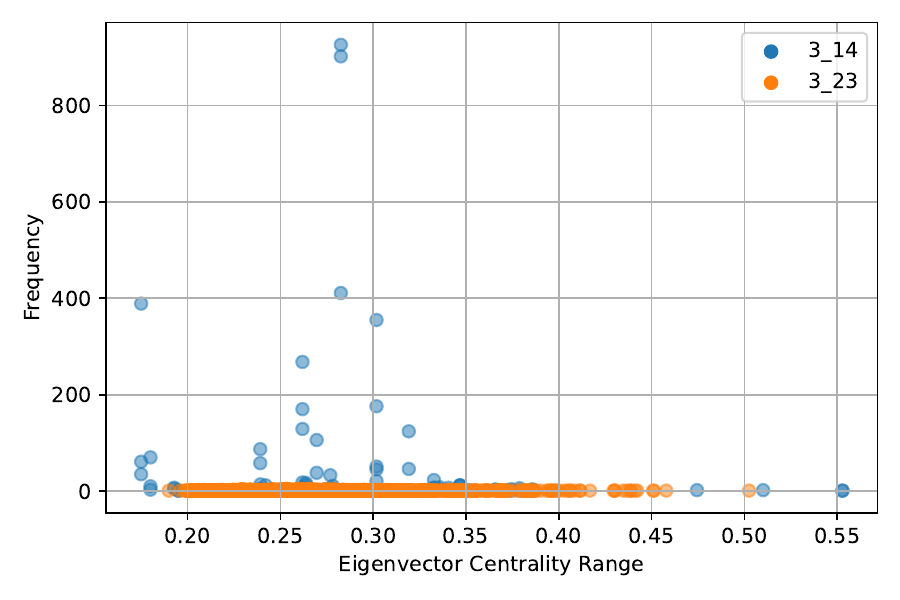}
        \caption{Centrality}\label{occ_centrality}
    \end{subfigure}
    \begin{subfigure}{0.47\textwidth}
        \centering
        \includegraphics[width=0.98\textwidth]{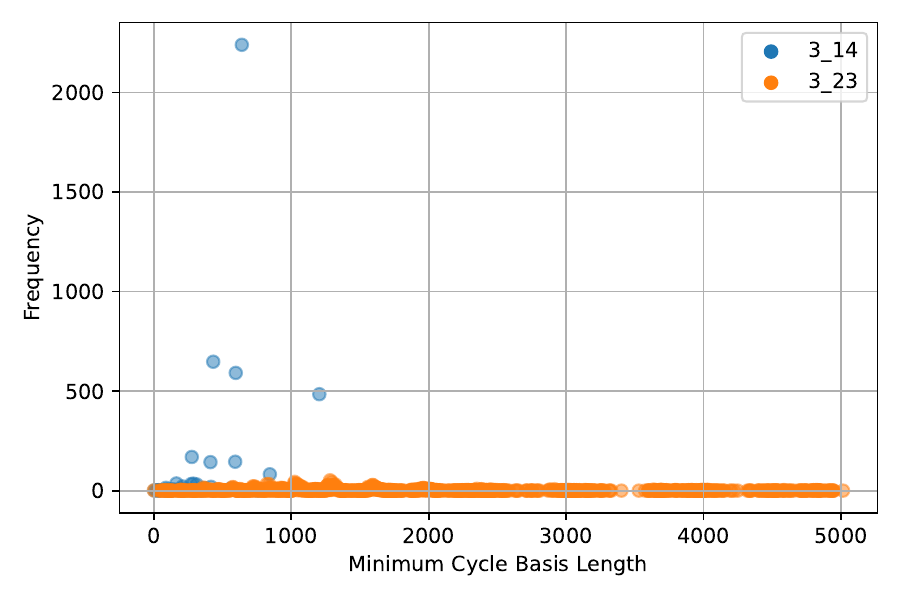}
        \caption{Cycles}\label{occ_cycle}
    \end{subfigure}
    \caption{Core network analysis properties of the Pachner graphs for the orientable cusped census partition with number of initial tetrahedra $<8$. Analysis considers building the Pachner graphs to depth 3 with either the $\{1-4, 4-1\}$ (denoted 3\_14) or the $\{2-3, 3-2\}$ (denoted 3\_23) moves, plotting histograms of: (a) the square clustering coefficients; (b) the normalised Wiener index; (c) the range of eigenvector centrality values; and (d) the length of the minimum cycle basis.}\label{occ_clustwienercent}
\end{figure}

\noindent \textbf{Clustering:} As in the previous Pachner graph analysis the triangle clustering coefficients for graphs across the census were all zero, since these graphs can only exhibit even length cycles.
However the square clustering coefficients had (min,mean,max) across the census for the 3\_14 graphs of $(0.028,0.038,0.104)$ and for the 3\_23 graphs of $(0.000,0.018,0.035)$, with only the 3\_23 graphs having examples of 0 square clustering coefficients, occurring for 2 manifolds. 
The full distributions are shown in Figure \ref{occ_clust}, they show that 4-cycles are always introduced with $1-4$ moves since all the initial triangulations have at least 2 tetrahedra so a $1-4$ move on each followed by their inverses in the other order automatically creates a 4-cycle. 
In general 4-cycles are less likely with $2-3$ moves since the gluings change more non-trivially and hence alternating inverses may be obstructed.
The 3\_14 graphs again have higher frequencies, indicating they are more likely to have similar graph structure.

\noindent \textbf{Shortest Paths:} The normalised Wiener index provides a global measure of the average shortest path length between any two nodes in the Pachner graph, that can be appropriately compared between graphs.
The (min,mean,max) values for the 3\_14 graphs were $(2.04,4.54,4.87)$ and for the 3\_23 graphs were $(3.08,5.08,5.48)$. 
In Figure \ref{occ_wiener} a histogram of these values are plotted also, showing that the 3\_14 graphs tend to have shorter shortest paths between nodes so that they can be traversed more easily, and it may be better to enact first these moves were possible in triangulation simplification methods. 
The 3\_14 graphs have higher frequencies, again supporting that graphs created by these moves are more similar and less dependent on the underlying manifold topology.

\noindent \textbf{Centrality:} 
Eigenvector centrality uses the eigenvector values corresponding to the adjacency matrix's unique largest eigenvalue to assign centrality scores.
Since these scores are relative their absolute values are less relevant, and instead focus in this work is placed on the range of values across a given Pachner graph (highest $-$ lowest centrality score).
The (min,mean,max) values of the centrality ranges for the 3\_14 graphs were $(0.175,0.271,0.553)$ and for the 3\_23 graphs were $(0.190,0.261,0.502)$; with 487 of the 4815 3\_14 graphs having the initial triangulation as the most central, and 485 times for the 3\_23 graphs. 
Therefore, of the 4815 manifolds only $\sim\frac{1}{4}$ have the initial node as the most central node in both cases, implying the generation process is highly non-symmetric. The similar centrality ranges for both move styles indicate similar adjacency eigendecompositions. 
Again higher 3\_14 frequencies indicate similar graph structures.

\noindent \textbf{Cycles:} 
As previously observed, Pachner graphs can only exhibit even cycles, and additionally 2-cycles are not possible in this Pachner graph design where undirected edges represent the ability to move between triangulations (not a unique move).
Since the graphs generated in the census were truncated at depth 3, the largest cycle size possible is 6, and hence the cycle bases for all these graphs were spanned by 4-cycles and 6-cycles.
The average number of 4 and 6-cycles in each minimum cycle basis for the census partition's generated Pachner graphs were $(155.16,0.00)$ for 3\_14 graphs and $(354.58,12.87)$ for 3\_23 graphs. 
These results show that although the Pachner graphs are only generated to depth 3 they are higher interconnected with a huge number of cycles.
The 3\_23 graphs have more cycles on average as there are many more move choices, and hence much more path degeneracy.
Surprisingly the 3\_14 graphs never exhibit 6-cycles, presumably because in all scenarios they are visible they were decomposable into 4-cycles in the basis length minimisation process (as the reordered inversion of any 3 $1-4$ moves to make a 6-cycle would have all reorderings of 2 moves making sub 4-cycles).
Higher frequencies for these 3\_14 graphs again indicate similar graph structures with similar cycle decompositions, less dependent on underlying topology.

\begin{figure}[!t]
    \centering
    \begin{subfigure}{0.47\textwidth}
        \centering
        \includegraphics[width=0.98\textwidth]{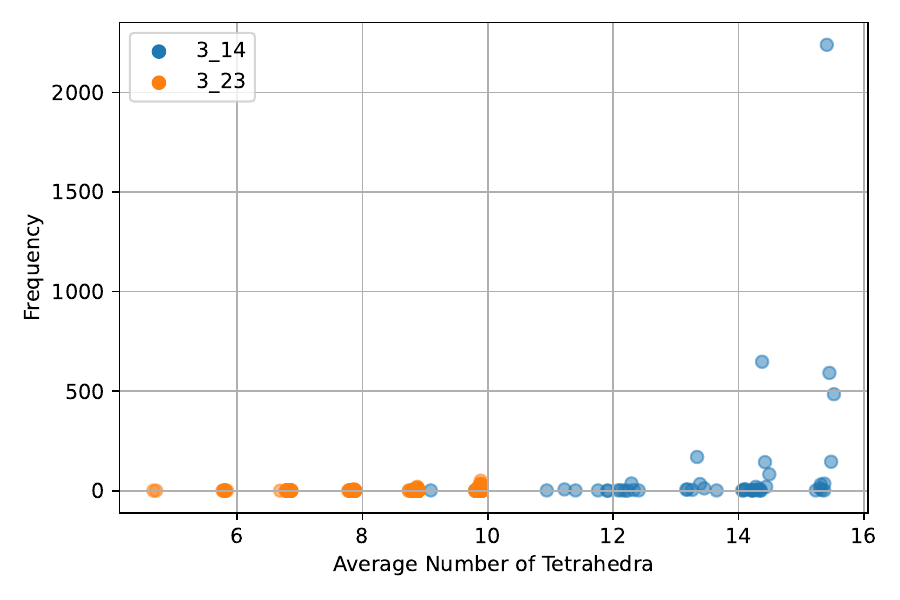}
        \caption{Number of Tetrahedra}\label{occ_avgnumtet}
    \end{subfigure} 
    \begin{subfigure}{0.47\textwidth}
        \centering
        \includegraphics[width=0.98\textwidth]{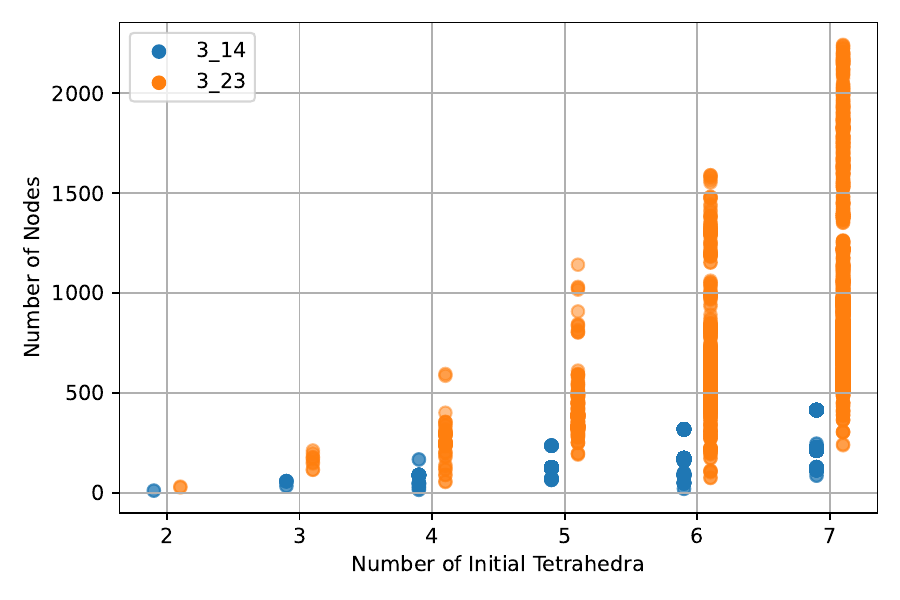}
        \caption{Pachner Graph Size}\label{occ_nodenumtet}
    \end{subfigure}
    \caption{Distributions of number of tetrahedra in triangulations across the 4815 Pachner graphs of the orientable cusped census partition with number of initial tetrahedra $<8$. Analysis considers building the Pachner graphs to depth 3 with either the $\{1-4, 4-1\}$ (denoted 3\_14) or the $\{2-3, 3-2\}$ (denoted 3\_23) moves, plotting (a) a histogram of the average number of tetrahedra across each Pachner graph; and (b) the number of tetrahedra in the initial triangulation against the number of nodes in the Pachner graph. Note the number of initial tetrahedra take integer values and the vertical lines are offsetted for visibility.
    }\label{occ_numtet}
\end{figure}

\noindent \textbf{Number of Tetrahedra:} 
The number of tetrahedra in each triangulation across a Pachner graph can only take 1 of 4 values (as only generated up to depth 3), based on an initial triangulation with $N$ tetrahedra and a move set of $\{2-3,3-2\}$ then the values are $\{N,N+1,N+2,N+3\}$, alternatively for move set $\{1-4,4-1\}$ the values are $\{N,N+3,N+6,N+9\}$.
This is because the moves always change the number of tetrahedra by a fixed amount, and therefore the number of tetrahedra in the initial triangulation is paramount to determining the spectrum, which is selected for the dataset to be a local minimum.
We note that this is not always a unique minimum as the frequency of triangulations with as many tetrahedra as the initial triangulation is occasionally $>1$.
The number of tetrahedra in the initial triangulation takes values of $(2,3,4,5,6,7)$ across the census partition, with respective frequencies $(2,9,56,234,962,3552)$.
Therefore most manifolds are initialised with triangulations with 7 tetrahedra.

Figure \ref{occ_avgnumtet} shows a histogram of the average number of tetrahedra across each Pachner graph, plotted for all members of the census partition.
It shows that the 3\_14 graphs have a higher average number of tetrahedra, expected since the $1-4$ moves introduce more tetrahedra than the $2-3$ moves (and both these are more common than their inverses).
There are also higher frequencies for these graphs as the graph structures are similar with each tetrahedron available to perform a $1-4$ move irrespective of its gluing.
Moreover, Figure \ref{occ_nodenumtet} shows how the number of tetrahedra in the initial triangulation correlates with the number of nodes in the generated Pachner graph.
Generally, more initial tetrahedra leads to more nodes and larger graphs, with larger ranges for the graph sizes as the initial number increases; this is related to more tetrahedra providing more possibilities for moves, facilitating the rapid graph growth with depth.
Since the number of places to apply $2-3$ moves grows quicker than for $1-4$ moves (each new tetrahedra introduce $>1$ new common faces) these 3\_23 graphs are more responsive to the initial number of tetrahedra but also exhibit a greater range of values as the combinatorics is more susceptible to limitations caused by complicated gluings.

Considering the limiting behaviour in the $3\_23$ graphs of a tree-like Pachner graph (no IsoSigs are repeated during generation) for the case of only considering $2-3$ moves (ignoring the $3-2$ inverse move), 
if one were to start with a triangulation of $N$ tetrahedra, in this limiting case this means there are $2N$ common faces to perform $2-3$ moves and hence $2N$ nodes at depth 1 of the Pachner graph. 
Each of these triangulations will have $N+1$ tetrahedra, and hence $2(N+1)$ faces for $2-3$ moves, leading to $4N(N+1)$ nodes at depth 2.
Similarly, these triangulations will now have $N+2$ tetrahedra, $2(N+2)$ faces for $2-3$ moves, leading to $8N(N+1)(N+2)$ nodes at depth 3.
In this limiting case one would then expect $8N(N+1)(N+2) + 4N(N+1) + 2N + 1 = 8N^3 + 28N^2 + 22N + 1 = 8N^3 + \mathcal{O}(N^2)$ nodes of the Pachner graph.
Mapping this leading order behaviour to the initial $N$ values considered in Figure \ref{occ_nodenumtet}, the predicted number of nodes would respectively be $(64, 216, 512, 1000, 1728, 2744)$ for the $N$ range $[2,7]$, which does approximately match the observed maximum counts in the plots.
The deviations from the maxima are when this limiting tree-like Pachner graph approximation breaks down, as triangulations are repeated in generation and the graph exhibits cycles (as well as minor effects from the less frequently occurring $3-2$ moves).
Equivalently in the $3\_14$ graphs for the limit case of only $1-4$ moves in a tree-like Pachner graph, an initial triangulation with $N$ tetrahedra has exactly this many tetrahedra to perform $1-4$ moves, leading to $N$ Pachner graph nodes at depth 1. 
Respectively there are then $N(N+1)$ nodes at depth 2, and $N(N+1)(N+2)$ at depth 3, for an overall estimate of $N^3 + 4N^2 + 4N + 1 = N^3 + \mathcal{O}(N^2)$.
The respective estimates would then be $(8, 27, 64, 125, 216, 343)$, again matching the observed order of the maxima in Figure \ref{occ_nodenumtet}.

\section{Machine Learning IsoSigs}
In recent decades, machine learning methods have seen a wealth of successes across many fields academia.
This can largely be attributed to the exponential growth in computational power, allowing the generation of ever-larger datasets for broader and more thorough analysis.

Mathematical data has a unique advantage in that it is exact, there is no noise from experimental apparatus, and as such any statistical patterns in data can be interpreted to have more direct meaning about the data's underlying structure.
However, new mathematical data can often be generated infinitely, and thus exhaustive computation becomes impossible, relying on appropriate sampling. 
Hence there is a true need for statistical methods when performing analysis on these large datasets. The computational statistical methods we explore here come under \textit{machine learning}.

In this work we examine the IsoSig representation of 3-manifold triangulations, testing how well can machine learning methods learn to identify manifolds from these IsoSigs.
Already, as demonstrated in §\ref{algorithms}, identifying the equivalence of the triangulations is highly computationally intensive. 
Where sequences of Pachner moves between two triangulations are not known to exist, the results from machine learning can provide a probability of existence of such a sequence moves.
Successes in learning IsoSig equivalence for 3-manifolds then opens up the scope for development into learning manifold properties and finding explicitly the connecting paths between triangulations. 

For these investigations the IsoSig input data is paired with an output label which identifies to which manifold the triangulation belongs to. 
The machine learning problem is in the style of supervised classification, where the network seeks to sort IsoSigs corresponding to different manifolds into separate classes.
In this setup, the networks are learning the combinatoric structure of triangulations which differentiates the manifolds, directly linked to their respective topology.

Specifically, for each of the eight 3-manifolds considered, a large Pachner graph was generated to depth 8 as described in §\ref{data}.
From there the frequency distribution of IsoSig lengths was computed across all IsoSigs in each manifold's Pachner graph.
Then an appropriate length was selected such that sufficiently large samples could be taken from each Pachner graph to produce datasets of fixed length IsoSigs\footnote{As noted previously, initial testing showed that padding mixed length IsoSigs caused the architectures to focus exclusively on the distribution of lengths, which was sufficiently different between manifolds to lead to good classification results. This simplification was avoided by considering a fixed length.}; the fixed length considered was 30 characters.

For computational processing these IsoSigs were converted to arrays of integers, using one-hot encoding of each character in the IsoSig.
Since the vocabulary consists of $64$ characters \{\texttt{a,b,c,...z,A,B,C,...,Z,0,1,2,...9,+,-}\}, the encoding of each character is a $64$-dimensional sparse vector of 0s with a single 1 in the dimension indexing that character in the vocabulary.
Therefore, because all IsoSigs considered for ML were $30$ characters long, each IsoSig was encoded by a $30 \times 64 = 1920$-dimensional vector.
Despite this representation being rather memory intensive, requiring more computational resources for training any architecture, the naturalness of the one-hot encoding style (where each character is an independent dimension) makes learning clearer and more interpretable\footnote{The other natural encoding of just mapping each character to a number $0 \mapsto 63$ was tested however lead to lower learning performance scores.}.

\subsection{Differentiating Manifolds}\label{sec:diffman}
In this section, the machine learning architectures of \textit{neural networks} (NNs) and \textit{transformers} are both applied to the problem of binary classifying IsoSigs belonging to different manifolds.
For the eight 3-manifolds focused on in this work, this creates ${8\choose2}=28$ pairs of manifolds to design paired data for, and perform binary classification on.
Creating 28 binary datasets to learn: \{\texttt{IsoSig}\} $\longmapsto \{0, 1\}$, the performance can then be compared between different manifold pairings, providing insight into what topological properties about the manifolds, encoded in the IsoSigs, the networks use to differentiate between manifolds.

The NN architecture is the prototypical machine learning model, it is very adaptable, and used ubiquitously with great success throughout many fields.
Any NN is made up of sets of neurons organised into layers, where each neuron takes in a vector of inputs, and outputs a single number.
The neuron acts on the input vector $\textbf{x}$ with a linear function, then a non-linear `activation' function, such that: $\textbf{x} \mapsto act(\textbf{w}\cdot\textbf{x}+b)$ for the trainable parameters of weights vector \textbf{w} and bias scalar $b$.
The full NN function of $L$-layers then takes the form
\begin{equation}
    f_{NN} := \text{softmax}(W_{i_L,i_{L-1}} \cdot \alpha(\ ...\ \alpha(W_{i_1,i_{input}} \cdot x_{input} + b_{i_1})\ ...\ ) + b_{i_L})\;.
\end{equation}
Demonstrating the vectorised linear action of matrix multiplication with the weight matrix $W_{i_m i_n}$ between layers $i_m$ and $i_n$, then addition of a bias vector $b_{i_m}$.
Followed by the $\alpha$ activation function used to provide the alternating non-linear action; here we use the Leaky-ReLU activation function, which has significant speed advantages due to its simplicity, and effectively makes the training process the fitting of a piecewise linear function.
Finally the softmax activation function is used to bound outputs into the range $[0,1]$ such that they produce a probability distribution for the classification task.
These activation functions are defined
\begin{align}
    \text{Leaky-ReLU}(x)_i & = \begin{cases} x_i & x>0\\ \varepsilon x & \text{otherwise} \end{cases}\,,\\
    \text{softmax}(x)_i & = \frac{e^{x_i}}{\sum_j e^{x_j}}\,,
\end{align}
for some small factor $\varepsilon$ ($\sim 0.01$).
These architectures are very adaptable, possess some universal approximation theorems supporting their use, and have practically been particularly successful across other fields of mathematics.
More information on this architecture can be found at \cite{gurney2018introduction}.

Conversely, the transformer architecture is much more sophisticated, and involves many stages of sub-NNs for encoding and decoding, primarily relying on the concept of attention \cite{vaswani2023attention} to learn context amongst the input strings -- note transformers do not require the tensorial encoding of the IsoSigs, but do require substantial computational resources, so are tested less extensively in this work.

Generally, transformer models are another class of supervised architectures designed to handle sequential language-like data using attention mechanisms instead of traditional recurrence or convolution. 
At their core, they employ an encoder-decoder structure, though many applications focus on encoder-only (e.g., BERT) or decoder-only (e.g., GPT) variants. 
The multi-head self-attention mechanism is the cornerstone, which allows the model to dynamically weigh the importance of all input tokens relative to each other. 
For an input sequence \texttt{abcd...} attention scores are computed using query $Q$, key $K$, ad value $V$ matrices, from which the attention importance weighting between letters of the input is computed as
\begin{equation}
    \text{Attention}(Q,K,V) := \text{softmax}(\frac{QK^T}{\sqrt{d_k}})V\;,
\end{equation}
for $d_k$ the dimensionality of the key vectors.
This mechanism enables Transformers to capture long-range dependencies efficiently.

Transformers work via alternating small NNs providing trainable parameters for the encoding and decoding, and attention layers which preserve the interrelations within the sequence. 
Transformers are trained using a supervised loss function in a similar manner to NNs, typically using cross-entropy loss, over a labelled dataset. 
They are also particularly popular as the attention mechanism allows parallelisation in computation, making larger models more feasible than traditional recurrent methods. 
More information on this architecture can be found at \cite{lin2022survey}.

As motivated in §\ref{algorithms}, the task of identifying equivalent triangulations is infeasibly expensive, and rigorous algorithm design is impractical.
Therefore efficient statistical methods may prove useful in informing on the likelihood of equivalence, thereby practically informing on optimal resource allocation for more rigorous checks; i.e. if a statistical method says a pair IsoSigs are very likely to be equivalent, it may be worth committing resources to compute the connecting path where it would be les likely found for other pairs.
Recent years have seen ML dominate the development of computational statistical methods across fields, and in this section we provide the first investigation into its efficacy on this equivalence problem and for general 3-manifolds as a whole.
In doing so we consider\footnote{More specialised architectures may well perform better, but we leave this experimentation to future work, since here we are performing the first tests of these methods for this problem.} the ubiquitous NN architecture, and the transformer architecture which has shown particular success on language-like data (alike our IsoSig encodings).

\subsubsection{Neural Networks}
\label{nns}
The NNs used in these investigations were built using \texttt{tensorflow} \cite{tensorflow2015whitepaper}, with architecture of three dense layers (sizes 256, 128, 64) then a final dense layer of size 2 for the one-hot encoded binary output; all layers used leaky-ReLU activation (factor 0.01), L2 regularisation, and dropout (factor 0.01) -- except for the output layer using softmax activation for the binary classification. 
The NNs were trained with an Adam optimiser to minimise the binary cross-entropy loss for 30 epochs with a batch size of 64 and a learning rate of 0.001.

Binary classification performances were measured with accuracy (the proportion of correctly classified test IsoSigs), and Matthew's correlation coefficient MCC (an unbiased alternative to accuracy). Both measures evaluate to 1 for perfect learning, whereas for no learning in a binary classification problem accuracy is 0.5 and MCC is 0.
Each binary classification investigation was performed 5 times on 5 different partitions of the IsoSig data into train and test sets, providing statistical confidence in the mean performance measures reported (this process is known as cross-validation).

\begin{figure}[!t]
    \centering
    \begin{subfigure}{0.24\textwidth}
        \centering
        \includegraphics[width=1.25\textwidth]{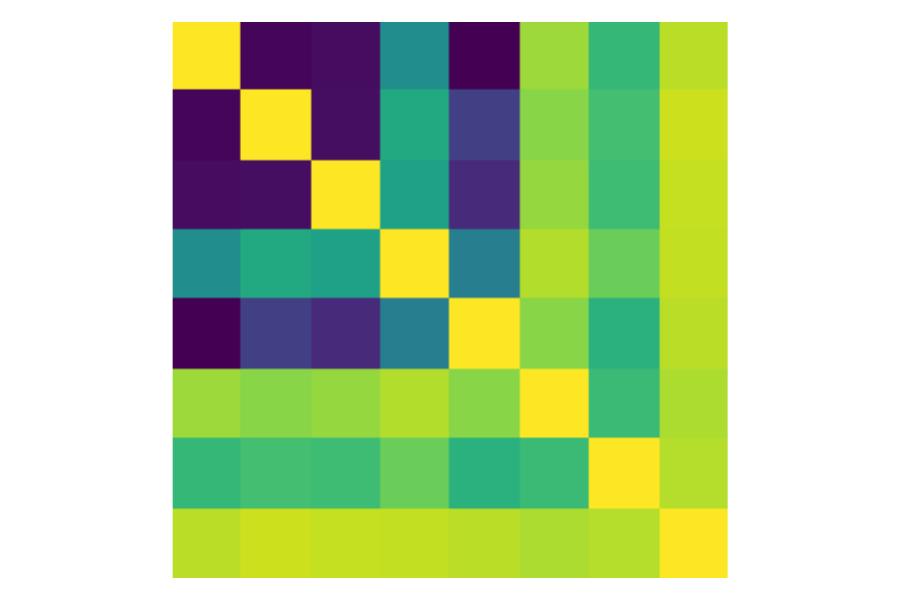}
        \caption{Accuracy}
        \label{BC_acc_23}
    \end{subfigure} 
    \begin{subfigure}{0.74\textwidth}
        \centering	
        \tiny
        $\begin{pmatrix} 1.0 & 0.49 & 0.5 & 0.737 & 0.484 & 0.923 & 0.828 & 0.947 \\ 0.49 & 1.0 & 0.502 & 0.796 & 0.582 & 0.909 & 0.845 & 0.96 \\ 0.5 & 0.502 & 1.0 & 0.778 & 0.545 & 0.919 & 0.839 & 0.955 \\ 0.737 & 0.796 & 0.778 & 1.0 & 0.705 & 0.941 & 0.883 & 0.953 \\ 0.484 & 0.582 & 0.545 & 0.705 & 1.0 & 0.908 & 0.814 & 0.948 \\ 0.923 & 0.909 & 0.919 & 0.941 & 0.908 & 1.0 & 0.836 & 0.937 \\ 0.828 & 0.845 & 0.839 & 0.883 & 0.814 & 0.836 & 1.0 & 0.943 \\ 0.947 & 0.96 & 0.955 & 0.953 & 0.948 & 0.937 & 0.943 & 1.0 \end{pmatrix}$
    \end{subfigure}\\
        \begin{subfigure}{0.24\textwidth}
        \centering
        \includegraphics[width=1.25\textwidth]{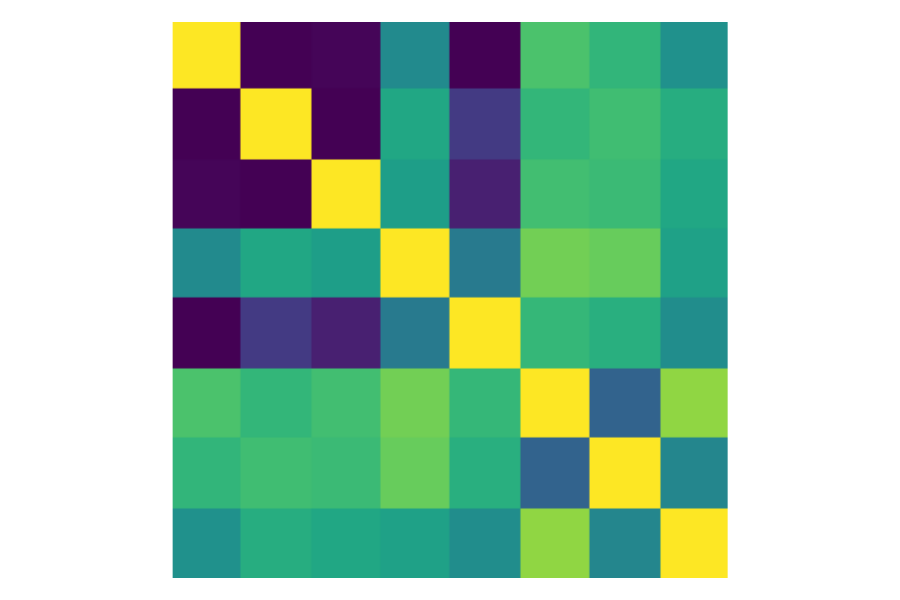}
        \caption{MCC}
    \end{subfigure} 
    \begin{subfigure}{0.74\textwidth}
        \centering	
        \tiny
        $\begin{pmatrix} 1.0 & 0.001 & 0.011 & 0.477 & 0.0 & 0.715 & 0.658 & 0.504 \\ 0.001 & 1.0 & -0.003 & 0.594 & 0.166 & 0.663 & 0.692 & 0.623 \\ 0.011 & -0.003 & 1.0 & 0.56 & 0.083 & 0.695 & 0.68 & 0.594 \\ 0.477 & 0.594 & 0.56 & 1.0 & 0.409 & 0.788 & 0.767 & 0.57 \\ 0.0 & 0.166 & 0.083 & 0.409 & 1.0 & 0.664 & 0.629 & 0.489 \\ 0.715 & 0.663 & 0.695 & 0.788 & 0.664 & 1.0 & 0.317 & 0.834 \\ 0.658 & 0.692 & 0.68 & 0.767 & 0.629 & 0.317 & 1.0 & 0.458 \\ 0.504 & 0.623 & 0.594 & 0.57 & 0.489 & 0.834 & 0.458 & 1.0 \end{pmatrix}$
    \end{subfigure}
    \caption{Matrices and their corresponding heat maps representing the (a) mean accuracy and (b) mean MCC scores for the binary classification between the datasets of length 30 IsoSigs generated using  $\{2-3, 3-2\}$ moves; averaged over the 5-fold cross-validation runs. Lighter colours represent higher performance scores, where accuracy evaluates in the range $[0,1]$, and MCC in the range $[-1,1]$. The matrix is manifestly symmetric as the pair order is irrelevant, and the diagonals are trivially 1 where no classification was required. The manifolds are ordered: $(S^3, S^2 \times S^1, \mathbb{RP}^3, L(7,1), L(7,2), T^3, PHS, H_{SC})$.}\label{BC_MCCs}
\end{figure}

The performance measure results are arranged in a matrix format in Figure \ref{BC_MCCs}. Each entry in the matrix $M_{ij}$ represents the mean measure score across the cross-validation runs for the binary classification of IsoSigs pairs between manifolds with index $i$ and index $j$ in the list: $(S^3, S^2 \times S^1, \mathbb{RP}^3, L(7,1), L(7,2), T^3, PHS, H_{SC})$\footnote{The diagonal cases $M_{ii}$ do not represent sensible experiments so had values set to 1, whilst the symmetric nature of the investigations was enforced by setting the lower diagonal equal to the upper diagonal.}.
Heat maps of these matrices are also shown for ease of comparison.

\begin{figure}[!t]
    \centering    
    \includegraphics[width=0.7\textwidth]{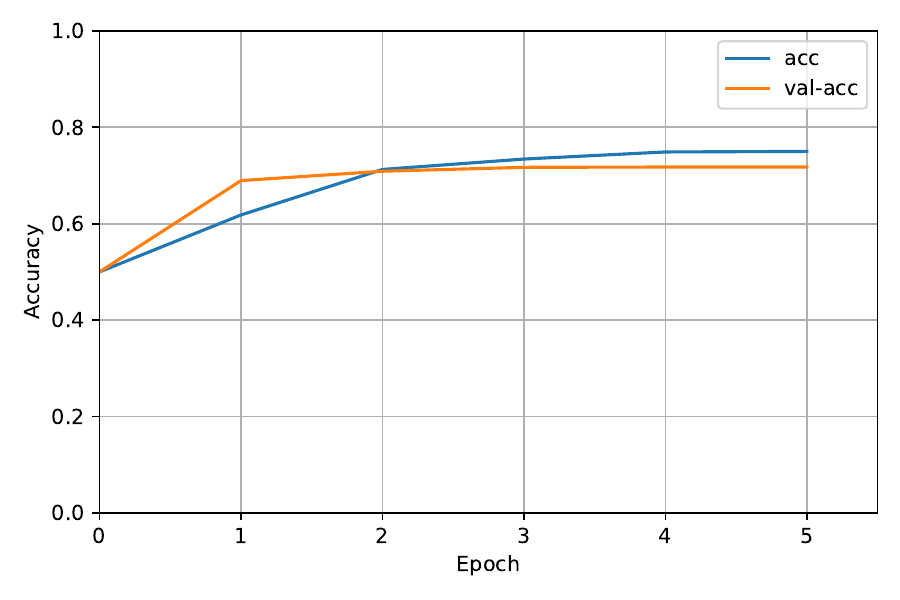}
    \caption{Example NN binary classification training curve, between length 30 IsoSigs for manifolds $L(7,1)$ and $L(7,2)$ generated using $\{2-3, 3-2\}$ moves up to depth 8. The plot shows that the accuracy performance measure for both training and validation sets increased quickly, converging to the final performance within a few epochs.}
    \label{training_curve}
\end{figure}

For each matrix the off-diagonal scores respectively have (min, mean, max): (a) Accuracy (0.4835, 0.800, 0.960), (b) MCC (-0.003, 0.487, 0.834); showing a range from no learning to near-perfect learning. 
Comparing between manifolds, the positive curvature manifolds are harder to differentiate amongst each other. 
The other binary classifications have surprisingly high performances in (a) and (b), noting particularly the success of the $L(7,1)$ vs $L(7,2)$ classification which are notoriously difficult manifolds to differentiate (they are homotopy equivalent but not homeomorphic).
To exemplify the learning behaviour, for this classification the accuracy performance measure was tracked throughout the training process on the training data used for updating the NN parameters and the independent validation data. 
The accuracy scores throughout training, evaluating at each epoch, are shown in Figure \ref{training_curve}.
This plot highlights that the performance improves quickly as the general behaviour can be learnt, but does plateau at the limiting score -- indicating there are triangulations which are truly difficult to differentiate the manifold structure between.

Interestingly, the manifold which is best distinguished from all the other $7$ is the Weeks manifold, the only hyperbolic one in the lot. 
At the opposite, $S^3,S^2\times S^1$ and $\mathbb{RP}^3=L(2,1)=SO(3)$ are not distinguished at all.

\paragraph{NN Gradient Saliency}\mbox{}\\
We performed an in-depth analysis of the gradient saliency values used by our NNs to distinguish $L(7,1)$ vs $L(7,2)$. 
A trained NN represents a function from inputs to outputs, in the case of these problems this is an input of a single IsoSig, one-hot encoded from length 30 in a length 64 alphabet to a length 1920 vector and an output of a length 2 vector with each entry evaluating in the range $[0,1]$ (as a probability distribution over the 2 classes) whose larger value indicates the predicted manifold the input IsoSig triangulation represents.

Computing the derivatives of the NN function outputs with respect to each of the inputs (i.e. for each index of each encoded letter in the IsoSig) provides a gradient function for each input. 
Evaluating these functions on each IsoSig in the test data gives numerical values for how sensitive the output classification is to each input value.
Averaging over all inputs in the test data, and across multiple cross-validation runs, then normalising the output gradients, provides relative gradient measures with greater statistical confidence. 

\begin{figure}[!t]
    \centering
    \includegraphics[width=0.7\textwidth]{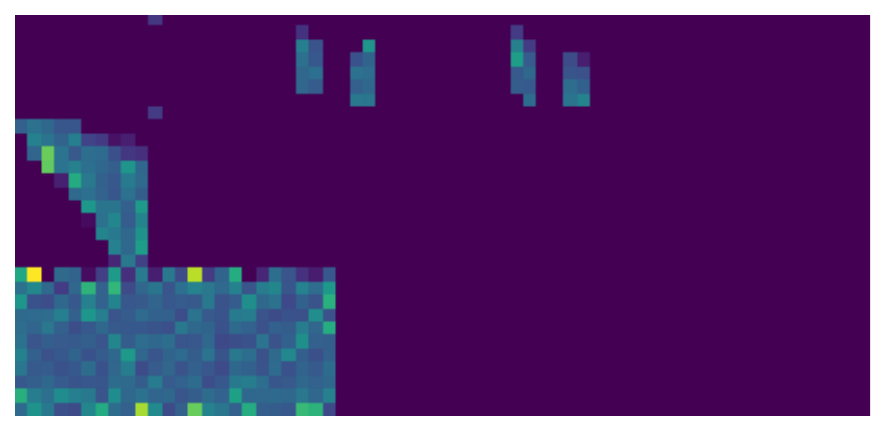}
    \caption{Image representation of the gradient saliency results for NN binary classification between IsoSigs from manifolds $L(7,1)$ and $L(7,2)$ generated using $\{2-3, 3-2\}$ moves. Each row is the encoding of one of the 30 IsoSig characters in the length $64$ alphabet, reshaping the 1920 length input vector into a $(30,64)$ size matrix. Lighter colours indicate a more significant sensitivity to that input for the output classification.}
    \label{saliency_image}
\end{figure}

For the $L(7,1)$ vs $L(7,2)$ classification, these gradient saliency values were averaged over 100 cross-validation runs, and the final gradient values are represented in Figure \ref{saliency_image}, where each of the 30 rows is an IsoSig character one-hot encoded into the length 64 alphabet.
Lighter colours indicate more sensitivity of the classification output to that input. The results in the figure show that the most important parts of the IsoSig input are later characters in the IsoSig and that these are more likely to be lower case letters from the start of the alphabet.

The rough behaviour of the Figure \ref{saliency_image}, respects the IsoSig encoding structure as described in Appendix §\ref{sec:isosigencoding}. The first row is $\pi(n)$, and as all these triangulations have 10 tetrahedra only the 10th character in the alphabet (\texttt{j}) has a non-zero saliency and a slightly brighter colour.
The next rows show importance on characters deeper in the alphabet, matching the triple encoding structure of the type sequence where higher index characters are regularly expressed.
Following this is a wide-diagonal-like structure up to the 10th character, which coordinates with the destination sequence being bound by $n$, which here is 10.
Since the tetrahedra are numbered as the destination sequence is built this explains why higher index characters are more important later on as the earlier new tetrahedra seen would be allocated lower numbers.
The final section of Figure \ref{saliency_image} shows a box of important values, up to the 23rd character, matching the permutation sequence behaviour where there're only 23 permutations in $S_4$ for relabelling vertices.
This part shows that all potentially observed permutations are important for the architecture's classification, and since the brightest colour are exhibited in this part this information tends to be the most useful for distinguishing triangulations of these geometries.

\begin{figure}[!t]
    \centering
    \begin{subfigure}{0.47\textwidth}
        \centering
        \includegraphics[width=0.98\textwidth]{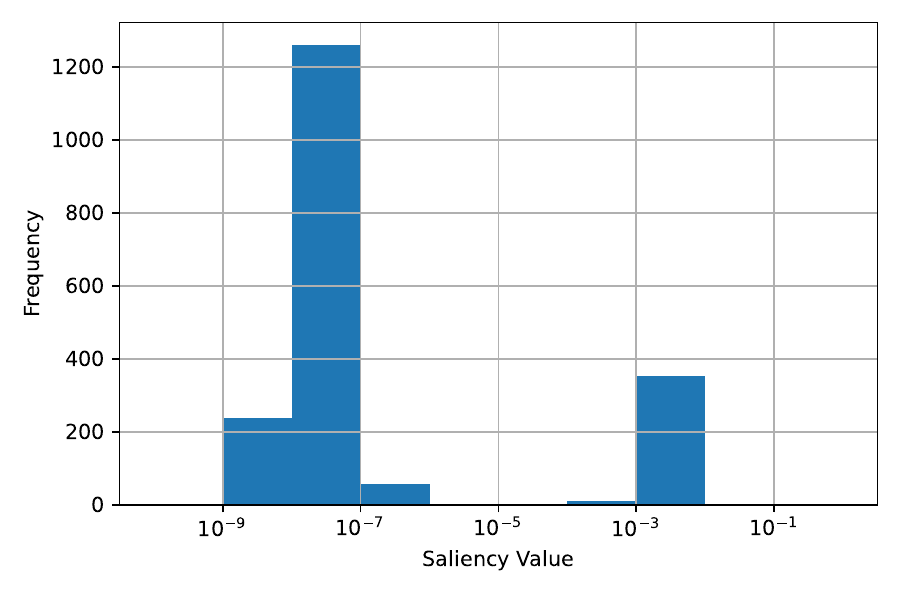}
        \caption{}
        \label{saliency_gradients}
    \end{subfigure} \\
    \begin{subfigure}{0.47\textwidth}
        \centering
        \includegraphics[width=0.98\textwidth]{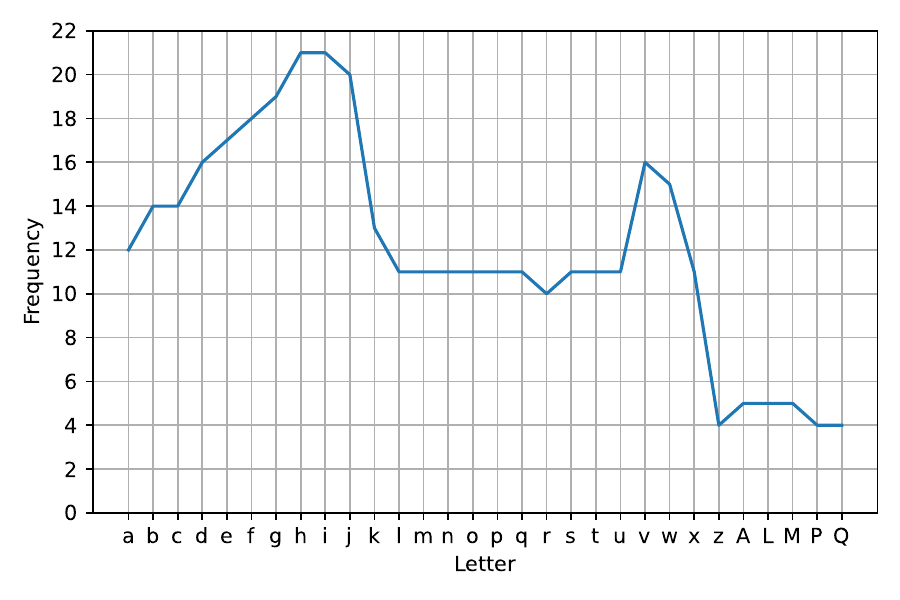}
        \caption{}
        \label{saliency_letters}
    \end{subfigure}
    \begin{subfigure}{0.47\textwidth}
        \centering
        \includegraphics[width=0.98\textwidth]{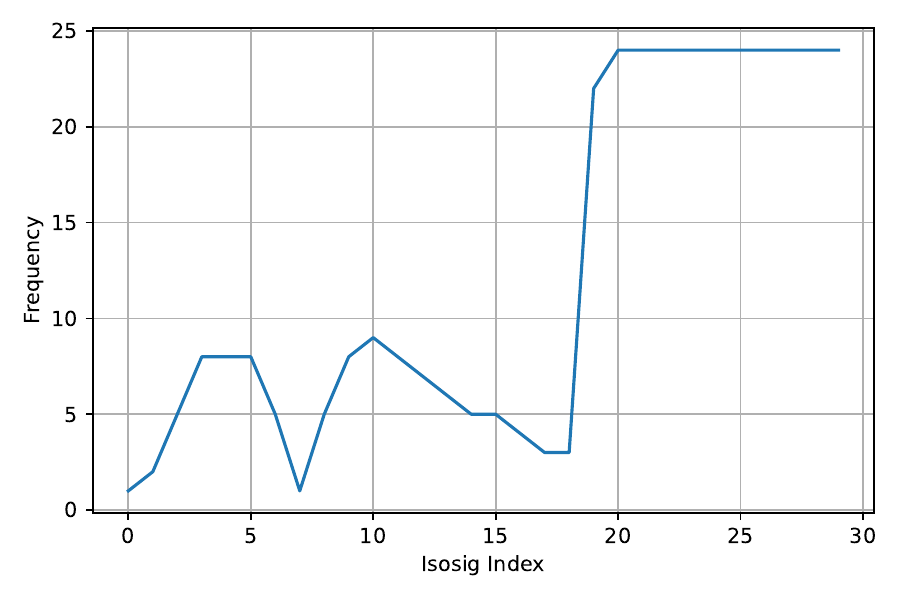}
        \caption{}
        \label{saliency_characters}
    \end{subfigure}
    \caption{Saliency analysis for the $L(7,1)$ vs $L(7,2)$ binary classification, averaged over 100 NNs. A histogram of the average saliency absolute values across the 1920 input vectors is shown in (a), demonstrating a partition into more/less important features. Subsequently (b) shows a histogram of the more important features with saliency value $>10^{-4}$ sorted according to the letter in the alphabet they correspond to; whilst (c) shows a histogram of the same set sorted according to which of the 30 characters in the input IsoSig they correspond to.}
    \label{saliencyanalysis}
\end{figure}

Further examining the actual saliency values (represented by colours in Figure \ref{saliency_image}), averaged over all of the cross-validation runs as well as each runs test dataset, we can plot a histogram of the average saliency size, as shown in Figure \ref{saliency_gradients}.
This image shows a partition of values with a selection having normalised absolute gradients $>10^{-4}$ and the majority having values $<10^{-6}$.
There is therefore a clear selection of `important' inputs and `non-important' ones.
Taking the input indices corresponding to the larger `important' saliency values $>10^{-4}$, these can be sorted into what character in the alphabet they correspond to (as plotted in a histogram in Figure \ref{saliency_letters}), or into what character in the length 30 IsoSig they correspond to (as plotted in a histogram in Figure \ref{saliency_characters}).

Figure \ref{saliency_letters} shows that the  363 `important' inputs correspond to only a subset of the full alphabet -- predominantly the beginning of this alphabet with a few letters deeper in.
Equivalently, these same important input features have much higher frequencies towards the end of the IsoSig as shown in Figure \ref{saliency_characters}.
From this saliency analysis the NNs are hence relying most on the final entries of the IsoSigs (the permutation sequence) to determine the manifold, and the value of these entries becomes discriminant.
Additionally the occurrence of specific letters (like \texttt{h} or \texttt{i} with the highest frequency in Figure \ref{saliency_letters}) can significantly help the architectures' confidence in classification.


\subsubsection{Transformers}
Since the IsoSig data closely resembles that of text, one may sensibly question how transformer architectures, which have seen great performance in many natural language processing problems \cite{openai2023gpt4,Douglas:2023olt}, perform when classifying this data.

A popular baseline for performing binary classification tasks with transformers is to use the BERT (Bidirectional Encoder Representations from Transformers) architecture \cite{devlin2019bert}.
In a similar vein to that employed in §\ref{nns}, the BERT model is now used to classify whether a given pair of IsoSigs represent the triangulations of the same manifold. 
This model is implemented using the \texttt{transformers} library \cite{wolf-etal-2020-transformers}.

Transformers are primarily successful as large language models, and due to this tend to have a very large number of trainable parameters. 
To make training feasible, their use often relies upon the principle of \textit{transfer learning}, where the models are pre-trained on a large standard set of data, and then retrained on a smaller specialised set.
This use of pre-trained models can avoid tens of thousands of core hours in retraining, and is what makes these large language models feasible to use at all.
However, the pre-trained models publicly available are trained on language data (often English language text).
This hence, raises appropriate concerns about the success of transfer learning onto this IsoSig 3-manifold `language' data, which does not resemble the syntactical and grammatical structure of the English language.

Nevertheless, due to their paramount successes in recent times, it is interesting to examine their performance in this context.
To do so we take the binary classification of the $L(7,1)$ and $L(7,2)$ datasets, as performed by the NNs to Accuracy and MCC scores of (0.705, 0.409) respectively, and repeat the training and testing using the BERT model.
Using a single layer output classifier to change the standard 768 output dimension of the BERT model to the 2 binary classes, along with equivalent hyperparameters such as training with the Adam optimiser and cross-entropy loss for 30 epochs with a batch size of 64 and learning rate of 0.001.
The transformer performance scores were:
\begin{equation}
\begin{split}
    \text{Accuracy} & = 0.689\;,\\
    \text{MCC} & = 0.381\;.
\end{split}
\end{equation}
Unfortunately, whilst some learning was achieved, the performance was not as good as the NN models, likely due to the inappropriate nature of the English language pre-training misleading the classification.
In addition to worse performance, the transformers also required significantly more computational resources (time, power, memory).
To exemplify this, whilst the NNs took the order of seconds to train and test, this simplest transformer model took the order of hours, and required orders of magnitude more memory to run. 

Therefore, although it was relevant experimentation, these results support the motivation to focus on the NN architecture for learning of triangulation IsoSig properties.

\subsection{Differentiating Knots}\label{sec:knots}
By a famous result of Gordon and Luecke, the complement of a knot entirely determines the knot \cite{GordonLuecke}. 
This allows one to identify knots with $3$-manifolds, as the knot complement embedded within $S^3$; and with this motivation the workflow we have set-up may also be directly extended to study and distinguish knots. 

\begin{figure}[!t]
    \centering
    \begin{subfigure}{0.24\textwidth}
        \centering
        \includegraphics[width=1.25\textwidth]{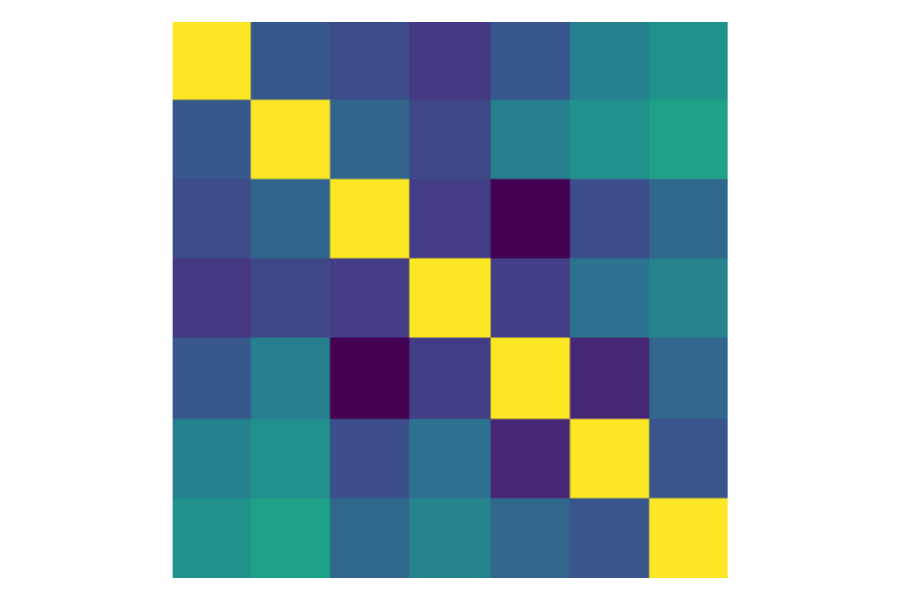}
    \end{subfigure} 
    \begin{subfigure}{0.74\textwidth}
        \centering	
        \tiny
        $\begin{pmatrix} 1.0 & 0.632 & 0.612 & 0.576 & 0.629 & 0.714 & 0.752 \\ 0.632 & 1.0 & 0.657 & 0.6 & 0.71 & 0.748 & 0.782 \\ 0.612 & 0.657 & 1.0 & 0.583 & 0.491 & 0.612 & 0.662 \\ 0.576 & 0.6 & 0.583 & 1.0 & 0.584 & 0.681 & 0.718 \\ 0.629 & 0.71 & 0.491 & 0.584 & 1.0 & 0.544 & 0.658 \\ 0.714 & 0.748 & 0.612 & 0.681 & 0.544 & 1.0 & 0.626 \\ 0.752 & 0.782 & 0.662 & 0.718 & 0.658 & 0.626 & 1.0 \end{pmatrix}$
    \end{subfigure}\\
    \begin{subfigure}{0.24\textwidth}
        \centering
        \includegraphics[width=1.25\textwidth]{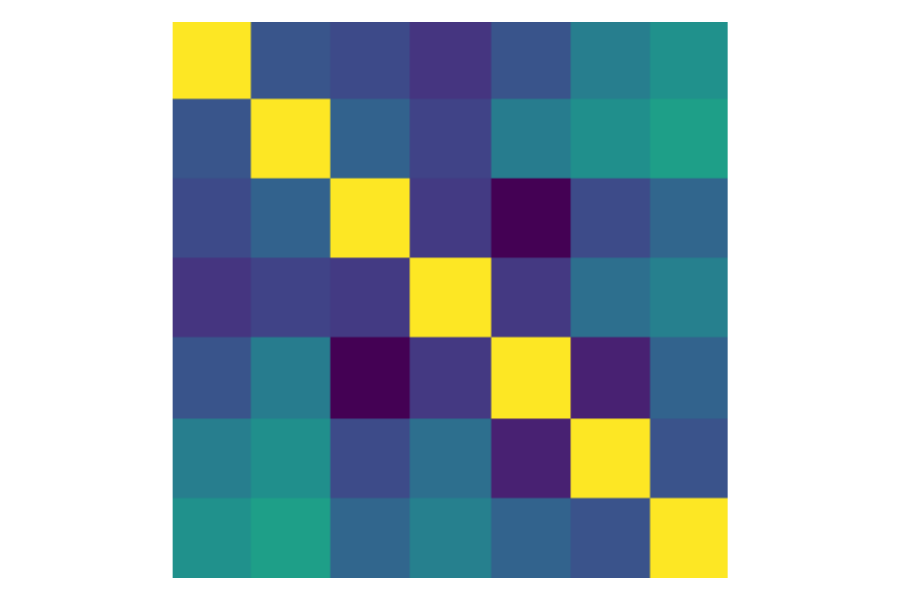}
    \end{subfigure} 
    \begin{subfigure}{0.74\textwidth}
        \centering	
        \tiny
        $\begin{pmatrix} 1.0 & 0.265 & 0.224 & 0.153 & 0.259 & 0.429 & 0.505 \\ 0.265 & 1.0 & 0.313 & 0.201 & 0.42 & 0.499 & 0.565 \\ 0.224 & 0.313 & 1.0 & 0.168 & 0.0 & 0.227 & 0.33 \\ 0.153 & 0.201 & 0.168 & 1.0 & 0.168 & 0.366 & 0.437 \\ 0.259 & 0.42 & 0.0 & 0.168 & 1.0 & 0.092 & 0.316 \\ 0.429 & 0.499 & 0.227 & 0.366 & 0.092 & 1.0 & 0.254 \\ 0.505 & 0.565 & 0.33 & 0.437 & 0.316 & 0.254 & 1.0 \end{pmatrix}$
    \end{subfigure}
    \caption{Matrices and corresponding heat maps representing the (a) mean accuracy and (b) mean MCC scores for the binary classification between the datasets of length 25 knot complement IsoSigs generated using $\{2-3, 3-2\}$ moves; averaged over the 5-fold cross-validation runs. Lighter colours represent higher performance scores, where accuracy evaluates in the range $[0,1]$, and MCC in the range $[-1,1]$. The matrix is manifestly symmetric as the pair order is irrelevant, and the diagonals are trivially 1 where no classification was required. The manifolds are ordered: (Unknot, Trefoil, FigureEight, $5_1$, $5_2$, DT $6a_3$, DT $8n_1$).}\label{Knot_MCCs}
\end{figure}

In this work, we again put focus on a selection of knots, in particular focusing on the most fundamental knots with a low number of crossings.
These knots are labelled: \{Unknot, Trefoil, FigureEight, $5_1$, $5_2$, DT $6a_3$ (a.k.a. $6_1$), DT $8n_1$ (a.k.a. $8_2$)\}, where the final knot is hyperbolic non-alternating.
The figure eight and the $5_2$ knot are the smallest hyperbolic knots in the $3$-sphere and indeed it turns out that they are more difficult to distinguish. In a similar way one could expect some difficulties in distinguishing the torus knots $3_1$ (the trefoil) and $5_1$ (the cinquefoil) but it is not so much the case (see Figure \ref{Knot_MCCs}). 

Again all the initial IsoSigs are 1-vertex, and are generated using only $2-3$ and $3-2$ moves.
The initial IsoSigs for these knot complement 3-manifolds are given in Table \ref{KnotInit_properties}, along with the respective number of tetrahedra in each initial triangulation.

\begin{table}[!t]
\centering
\begin{tabular}{|c|c|c|}
\hline
Knot         & Initial IsoSig             & \# Tetrahedra \\ \hline
Unknot       & \texttt{cMcabbgds}         & 2             \\ \hline
Trefoil      & \texttt{cPcbbbadu}         & 2             \\ \hline
Figure Eight & \texttt{cPcbbbiht}         & 2             \\ \hline
$5_1$        & \texttt{dLQbcccaekv}       & 3             \\ \hline
$5_2$        & \texttt{dLQbcccdero}       & 3             \\ \hline
DT $6a_3$=$6_1$ knot    & \texttt{eLPkbcddddcwjb}    & 4             \\ \hline
DT $8n_1$=$8_2$ knot    & \texttt{fLLQcbcdeeedowxxd} & 5             \\ \hline
\end{tabular}
\caption{The knot IsoSigs used to start the generation of Pachner graphs using only $2-3$ and $3-2$ moves. Each knot was used to construct a knot-complement 3-manifold, which was then triangulated to give a 1-vertex triangulation. The respective IsoSigs are listed in this table, along with the number of tetrahedra in this initial triangulation.}
\label{KnotInit_properties}
\end{table}

The Pachner graphs were generated up to a depth that provided a suitable amount of data for a fixed IsoSig length. 
Again searching for 2000 IsoSigs per Pachner graph, the selected IsoSig length was 25 characters.
From there the same NN architecture was trained to perform binary classification between the knot complement manifolds.
The binary classification learning was performed for each pairing of knot-complement 3-manifolds, where the accuracy and MCC scores are displayed in Figure \ref{Knot_MCCs}.

The results show comparable performances between knots, however with worse maximum performances compared to the general 3-manifold learning in §\ref{sec:diffman}.
Some pairings had close to no learning, better exemplified by the near-zero MCC scores for the (FigureEight, $5_2$) and the ($5_2$, DT $6a_3$) pairings.
Generally, as shown by the images, the colours are darker closer to the diagonal, indicating worse performance, and suggesting that in general the more disparate the crossing structure of the knot the easier the NNs can distinguish them. 
Accuracies for identifying the Unknot are not as high as in the focused Unknot study in \cite{Gukov:2020qaj}. 

\subsubsection{Supplementing with Dehn Surgery}\label{sub:dehnsurg}
In addition to learning with their knot-complement manifolds, we also perform two Dehn surgery operations (0-surgery and 1-surgery) whereby a tubular neighbourhood of a knot is drilled out and then filled with a solid torus \cite{Dehn}.

Dehn surgery is a standard operation in topology allowing modification of the topology of a $3$-manifold by a `cut and paste' operation on the $D^2\times S^1$-neighborhood of a knot consisting in re-gluing the solid torus as $S^1\times D^2$. This operation is also called $0$-Dehn surgery and more in general $n$-Dehn surgery consists in regluing the solid torus after twisting the image of the meridian $n$-times along the longitude. 
In particular operating $0$ or $1$ surgeries along knots in $S^3$ one produces closed oriented $3$-manifolds with the homology respectively of $S^2\times S^1$ and of $S^3$. These manifolds are themselves canonically associated to the knots so their homeomorphism type is by itself a knot invariant, albeit a complicated one. 
Here we try and distinguish knots by distinguishing these surgery manifolds. 

Focus is put on the differentiation between the Unknot and Trefoil binary classification problem, which for the knot complement $3$-manifolds classified with accuracy 0.632 and MCC 0.265, neither particularly impressive scores.
In this extension, the initial knot complement triangulation IsoSigs, as given in Table \ref{KnotInit_properties}, had either 0-surgery or 1-surgery performed on them, to produce further $3$-manifolds whose IsoSigs are given in Table \ref{SurgeryInit_properties}. We remind the reader that the $0$ and $1$ surgeries on the Unknot are respectively $S^2\times S^1$ and $S^3$ and on the left-handed Trefoil they are two Seifert fibered manifolds, respectively $(S^2,(2,1),(3,1),(6,-5))$ and $(S^2,(2,1),(3,1),(7,-6))$ (the latter \emph{not} being the Poincar\'e homology sphere, which is the $1$ surgery over the \emph{right}-handed Trefoil). Each of these was used to seed Pachner graph generation to generate databases of 2000 length 30 IsoSigs.
Noting that in this subsection length 30 was used again instead of the length 25 as used in §\ref{sec:knots}, this is because the initial IsoSigs after surgery lead to length distributions with very few length 25 IsoSigs\footnote{This is also why just these two simplest knots are considered since the length distributions for some knot complements after surgery did not even include any length 25 or 30 IsoSigs.}.
Another comment worth emphasising is that whilst the IsoSigs in Table \ref{SurgeryInit_properties} have a range of initial number of tetrahedra, since only a fixed length IsoSig is taken from the graph for ML these IsoSigs will also correspond to triangulations with the same number of tetrahedra.
Additionally, it is worth noting that the initial IsoSig from 0-surgery of the Unknot is a triangulation with 3-vertices (whilst the others are still all 1-vertex), this means this Pachner graph is a no longer a 1-vertex Pachner graph.

\begin{table}[tb]
\centering
\begin{tabular}{|c|c|c|c|}
\hline
Manifold                 & Surgery & IsoSig                          & \# Tetrahedra \\ \hline
\multirow{2}{*}{Unknot}  & 0       & \texttt{cMcabbjaj}              & 2             \\ \cline{2-4} 
 & 1       & \texttt{cMcabbgag}              & 2             \\ \hline
\multirow{2}{*}{Trefoil} & 0       & \texttt{gvLQQcdefeffnjndspx}    & 6             \\ \cline{2-4} 
& 1       & \texttt{hLLLQkcdefgfgghsdaenjw} & 7             \\ \hline
\end{tabular}
\caption{The surgery Isosigs used to start generation of Pachner graphs using only $2-3$ and $3-2$ moves. 
}
\label{SurgeryInit_properties}
\end{table}

With these datasets of 2000 length 30 IsoSigs corresponding to surgeries of these knot complements, the respective ML binary classification between the knots was repeated for each surgery operation.
The NN architecture used was the same as throughout this work, as quoted in §\ref{nns}.
The 5-fold cross-validation performance scores are given in Table \ref{SurgeryResults}, including the values for the knot complements without surgery repeated from §\ref{sec:knots}.

\begin{table}[tb]
\centering
\begin{tabular}{|c|c|c|}
\hline
Investigation   & Accuracy & MCC   \\ \hline
Knot Complement & $0.632$    & $0.265$ \\ \hline
0-surgery       & $0.935$    & $0.871$ \\ \hline
1-surgery       & $0.974$    & $0.949$ \\ \hline
\end{tabular}
\caption{NN binary classification performance scores for the knot complements, the 0-surgeries, and the 1-surgeries, between the Unknot and Trefoil knots. Scores were averaged over 5 cross-validation runs.}
\label{SurgeryResults}
\end{table}

Performances are substantially higher after surgery, for either surgery actions. 
It is interesting to observe that since $0$ and $1$-surgeries over the unknot are $S^2\times S^1$ and $S^3$ respectively, the problem discussed here is yet another instance of the problem analysed in §\ref{sec:diffman} and the accuracies for the $0$-surgery and $1$-surgery cases are comparable with those seen in Figure \ref{BC_MCCs}.
These performances are quite high and show that one can easily complement the initial performances for the simple knot complements with further parallel analyses. 
The potential for distinguishing knots seems therefore quite high.

\section{Conclusions and Outlook}
In this work triangulations of 3-manifolds were studied via their IsoSig representation. 
A selection of 3-manifolds were chosen to be of focus in this study, along with a selection of knot complements from knots exhibiting a low number of crossings.
Their databases of triangulations were generated through their Pachner graphs, exhaustively performing $\{2-3, 3-2\}$ moves on an initial triangulation up to some pre-specified depth.
The graph structure of these graph was also analysed through a variety of network analysis techniques.

Pachner graph network analysis made some preliminary connections between the 3-manifold topological properties and how their respective Pachner graph grows, particularly using eigenvector centrality and decomposition of their minimum cycle bases.
This included proposal of a conjectural relation (Conjecture \ref{conj:systole}) between the growth rate of the Pachner graph and the systole of a hyperbolic $3$-manifold. 
Machine learning results showed that through a one-hot encoding of the IsoSigs that simple neural network architectures could learn to well distinguish which manifold a triangulations' IsoSig corresponded to -- a feat very difficult by eye due to the high information density of this encoding scheme.
Gradient saliency analysis then demonstrated that the neural networks were using the IsoSig characters in a manner which respected the encoding schemes decomposition into type, destination, and permutation sequences.

Transformers could not reach the same performance scores as the neural networks, likely due to the highly unique language-style of this representation.
However, the neural network architectures could generalise well to differentiating these IsoSig representations of knot complements, and particularly well manifolds after the operation of Dehn surgery.

At this work's respective \href{https://github.com/edhirst/IsoSigPGML.git}{GitHub}, code functionality is provided for generating and visually representing Pachner graphs for generic 3-manifolds; as well as performing a range of network analysis techniques. Scripts for repeating the machine learning investigations performed here are also provided, as well as databases of IsoSigs for the use of interested readers.

Future work would aim to make the connections between Pachner graph network properties and manifold topological properties more concrete.
Further to this, machine learning methods may be applied directly to the Pachner graphs, approximating optimal search techniques for problems such as minimising number of tetrahedra in a representation, or conversely using graph-neural-network techniques to study them.
Additionally the study of more manifolds, and of triangulations deeper into the Pachner graphs would be interesting to corroborate results observed here.

\section*{Acknowledgments}
F.C. was supported by CIMI Labex ANR 11-LABX-0040 at IMT Toulouse within the program ANR-11-IDEX-0002-02.
YHH would like to thank STFC for grant ST/J00037X/2 and a Leverhulme project grant.
E.~Heyes would like to thank SMCSE at City, University of London for the PhD studentship, as well as the Jersey Government for a postgraduate grant.
E.~Hirst acknowledges support from Pierre Andurand over the course of this research. 
This research utilised Queen Mary's Apocrita HPC facility \cite{apocrita}, supported by QMUL Research-IT.

\addcontentsline{toc}{section}{References}
\bibliographystyle{utphys}
\bibliography{references}

\newpage
\appendix

\section{IsoSig Encoding}\label{sec:isosigencoding}
Isomorphism Signatures, or IsoSigs, as introduced in \cite{burton2011simplification}, are a way of representing 3-manifold triangulations. 
In addition to uniquely representing each isomorphism class, they are both low memory and fast to compute/check/search.

The algorithm to produce an IsoSig for a given 3-manifold triangulation with $N$ tetrahedra starts with the generation of \textit{canonical labellings}.
Before defining a canonical labelling we will first define a labelling, which is generated as follows:
\begin{enumerate}
    \item Arbitrarily label the tetrahedra from 0 to $N-1$.
    \item Arbitrarily label each tetrahedra's four vertices from 0 to 3. Respectively the face opposite each vertex $i$ is given the same label $i$.
    \item Using the notation $A_{t,f}$ to equal the label of the tetrahedra glued to face $f$ on tetrahedra $t$ (such that $A_{t,f} \in \{0,...,N-1\}$); list the gluings as: \\ $(A_{0,0},A_{0,1},A_{0,2}, A_{0,3},A_{1,0},...A_{N-1,3})$ to define the labelling.    
\end{enumerate}
The are many possible labellings that uniquely define a triangulation, we first reduce this set to only consider canonical labellings, which are defined as labellings with the properties:
\begin{itemize}
    \item Tetrahedron $t$ first appears in the labelling before tetrahedron $t'$ iff $t<t'$, for $t,t' \geq 1$ since 0 is implicitly first.
    \item Where tetrahedron $t$ first appears in the list at position $A_{t',f'}$, the labelling of tetrahedra $t$'s vertices matches the labelling of tetrahedra $(t')$'s  -- i.e. each of the three vertices in face $f'$ of tetrahedra $t$ has the same vertex labelling as the vertex it is glued too in tetrahedra $t'$, the remaining fourth vertex labelling is then inferred from these three.
\end{itemize}

For any triangulation, a canonical labelling can be directly generated by first selecting any of the $N$ tetrahedra as the $0$th tetrahedron; then arbitrarily labelling its vertices (of which there're $4!=24$ ways).
From this choice the remainder of the canonical labelling can be inferred by following the gluings list -- each time a new tetrahedron is introduced label it as the next one in the list and label its vertices using the identity map such that its vertex labels match the current tetrahedron in accordance with the canonical labelling definition above.
Through this labelling process, it can be seen that there're $24N$ canonical labellings for any triangulation \cite{burton2011simplification}.

This process of allocating a canonical labelling is clearly still redundant.
To resolve this many-to-one map from labellings to triangulations, all canonical labellings are computed, encoded as IsoSigs, then the lexicographically smallest encoding is selected to represent the triangulation.
How the labellings are encoded is what gives the IsoSigs their character-based / language-like representation, as well as being significantly more efficient.

A first step in improving the encoding efficiency is using an alphabet, and when the number of tetrahedra $N$ is large this becomes particularly effective.
The alphabet\footnote{To encode larger numbers $d = \lfloor \text{log}_{64}(N) \rfloor + 1$ characters can be used by just converting all numbers to base 64, fixing $d$ characters for all encoded integers, and using padding when necessary. To exemplify, for $N=80$, then $d=2$ characters are needed, and the number $12$ would be written as $(12, 0)$ in the $(64^0, 64^1)$ basis, such that the encoding would be \texttt{ma}.} implemented is in base 64, such that:
\begin{table}[!h]
\centering
\begin{tabular}{|c|c|c|c|c|c|c|c|c|c|c|c|}
\hline
Integer   & 0          & $\cdots$ & 25         & 26         & $\cdots$ & 51         & 52         & $\cdots$ & 61         & 62         & 63         \\ \hline
Character & \texttt{a} & $\cdots$ & \texttt{z} & \texttt{A} & $\cdots$ & \texttt{Z} & \texttt{0} & $\cdots$ & \texttt{9} & \texttt{+} & \texttt{-} \\ \hline
\end{tabular}
\end{table}

Redundant information in the gluings list can also be omitted.
In this process, the gluings list through the triangulation faces $A_{t,f}$ is re-encoded into three parts: (i) destination sequence, (ii) type sequence, (iii) permutation sequence.
Next we describe each of these part individually:
\begin{enumerate}
    \item[(i)] Destination Sequence: Since each face gluing involves two faces, the gluings list for each face naturally includes the information about each gluing twice. Therefore, from the gluings list $(A_{0,0},A_{1,0},...,A_{N-1,3})$, each time a face $A_{t,f}$ is glued to another face $A_{t',f'}$ which has already been listed (i.e. $t' < t$, or $t' = t$ and $f' < f$), the entry $A_{t,f}$ is omitted from the gluings list. In addition for more general triangulations with boundary, a character $\partial$ may be used instead of the tetrahedron number for $A_{t,f}$ to evaluate to, when the face $f$ is a boundary of the triangulation. The destination sequence is then a list of all the unique gluings (and boundaries).
    \item[(ii)] Type Sequence: Each term in the destination sequence is then classified as one of three types. Either type 0 $\rightarrow$ the face is a boundary ($A_{t,f} = \partial$); type 1 $\rightarrow$ the face is glued to a new tetrahedron; or type 2 $\rightarrow$ the face is glued to a previously seen tetrahedron. The type sequence is then a list of the numbers 0, 1, or 2 for each term in the destination sequence.
    \item[(iii)] Permutation Sequence: For each non-boundary face, the gluing information is described by a permutation of the vertices from the face in the first tetrahedron to the face in the second (each of the three vertices in the face are appropriately paired, and the remaining fourth vertex in each case can then be inferred). Indexing all the 24 possible permutations $p_{t,f} \in S_4$ from $0$ up to $23$ (i.e. $(0,1,2,3) \vcentcolon = 0, (0,1,3,2) \vcentcolon = 1$, etc.), the permutation sequence is then a number for each non-boundary term in the destination sequence.
\end{enumerate}

The full IsoSig is then a further reduction of this information. The map $\pi(\cdot)$ takes any integer into its alphabet encoding; if $N<63$ the IsoSig starts with $\pi(N)$; when $N\geq 63$, the sequence starts\footnote{Since 63 is the highest number which can be encoded if the sequence starts $\pi(63)$ it is considered as `63 or above', then the next character is taken to be $\pi(d)$ followed by the full $\pi(N)$ (which may be 63 again).} with 
$\pi(63) \pi(d) \varepsilon(N)$, where $\varepsilon(\cdot)$ encodes each of the digits of $N$ written modulo $64$ (so when $N<64$ the $d=1$ encoding is implicit).


Next the type sequence is encoded and listed, collecting three consecutive terms into each printable character via $\pi(\tau_i + 4\tau_{i+1} + 16\tau_{i+2})$  (since each integer $\tau \in \{0,1,2\}$), padding the end with zeros if necessary.
Then the destination sequence is encoded and listed via $\pi(A_{t,f})$, but \textit{only} for the type 2 faces (since type 0 are boundary, and type 1 can be deduced as the next tetrahedra label since the labelling is canonical).
Finally, the permutation sequence is encoded and listed, but again \textit{only} for the type 2 faces (since the type 0 faces have no gluing, and the type 1 faces are implicitly the identity permutations since the labelling is canonical).

This produces the encoding of the canonical labelling, however as previously mentioned there are $24N$ of these.
To uniquely specify the IsoSig, each of these encoded canonical labellings is computed, and the lexicographically smallest is selected as the IsoSig.
The lexicographical ordering used is the ASCII computational standard \cite{ansi1975}. 

\begin{wrapfigure}[17]{r}{0.30\textwidth}
    \begin{center}
    \includegraphics[width=0.25\textwidth]{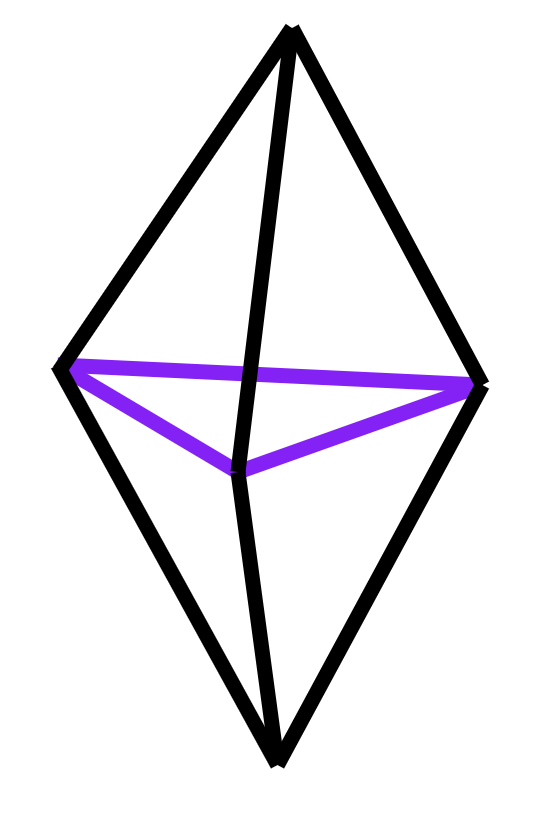} 
    \caption{The triangular bipyramid, an example 3-manifold triangulation, where the coloured edges indicate the glued face.}
    \label{fig:bipyramid}
    \end{center}
\end{wrapfigure}

\paragraph{Example: Triangular Bipyramid}\mbox{}\\
As a simple example let us consider the triangular bipyramid, as depicted in Figure \ref{fig:bipyramid}. 
This has two tetrahedra, which is far less than 63, so the encoding begins with \texttt{c}$=\pi(2)$ (the $d=1$ is implicit). 

Arbitrarily starting with the top tetrahedron as our tetrahedron 0, let us choose the vertex labelling such that the single glued face has vertices 0, 1, and 2; with the remaining vertex labelled 3 (which matches the labelling of the glued face as it is opposite).
The gluings list then reads $A_{0,0} = \partial$, $A_{0,1} = \partial$, $A_{0,2} = \partial$, $A_{0,3} = 1$,  now this is the first occurrence of this new tetrahedron (bottom) and hence sets its tetrahedron label to the next integer which is 1; furthermore the canonical nature of the labelling sets the vertex labellings for this second tetrahedron, each vertex on the glued face in tetrahedron 0 is glued to a vertex of the same labelling in tetrahedron 1.
The gluings list hence continues $A_{1,0} = \partial$, $A_{1,1} = \partial$, $A_{1,2} = \partial$, $A_{1,3} = 0$.
The overall gluings list is then: $(\partial,\partial,\partial,1,\partial,\partial,\partial,0)$.

This reduces to the destination sequence: $(\partial,\partial,\partial,1,\partial,\partial,\partial)$, removing the repetition of the initial tetrahedron at the end of the gluings list.
The type sequence is then defined as: $(0,0,0,1,0,0,0)$ over the destination sequence, i.e. there are no non-trivial gluings (type 2's). 
Finally the permutation sequence is computed, which is only defined for non-boundary terms in the destination sequence: $(0)$.

These sequences are then encoded and combined, after the initial \texttt{c}, the type sequence is then split into threes and padded as $((0,0,0),(1,0,0),(0,0,0))$, becoming $\pi(0 + 4 \times 0 + 16 \times 0) = \pi(0) = $ \texttt{a}, $\pi(1 + 4 \times 0 + 16 \times 0) = \pi(1) = $ \texttt{b}, $\pi(0 + 4 \times 0 + 16 \times 0) = \pi(0) = $ \texttt{a}, which together is \texttt{aba}.
Following this, since there are no type 2 faces, the destination and permutation sequences are trivially empty, and this defines the full encoded canonical labelling as \texttt{caba}, which to check is canonical we would have to consider all $24 \times 2 = 48$ canonical labellings and choose the lexicographically smallest (which in this case it is!).

\section{Pachner Graphs}\label{sec:extra_PGs}
The Pachner graphs generated with moves $\{2-3,3-2\}$ up to depths 3 and 4, for the remaining 8 considered manifolds, starting from triangulations with IsoSigs specified in Table \ref{tab:manifoldisosigs}.
The nodes represent triangulations and are labelled with the triangulation's respective number of tetrahedra, edges represent Pachner moves which transform the triangulations they connect into each other.

The smallest closed hyperbolic manifold $H_{SC}$ had too many IsoSigs in its Pachner graph to plot at even depth 3, therefore the depths 1 and 2 are shown instead.
To exemplify the speed of Pachner graph growth, the depth 5 $S^3$ Pachner graph is also given at the end.

\begin{figure}[b]
    \centering
    \begin{subfigure}{0.47\textwidth}
        \centering
        \includegraphics[width=0.98\textwidth]{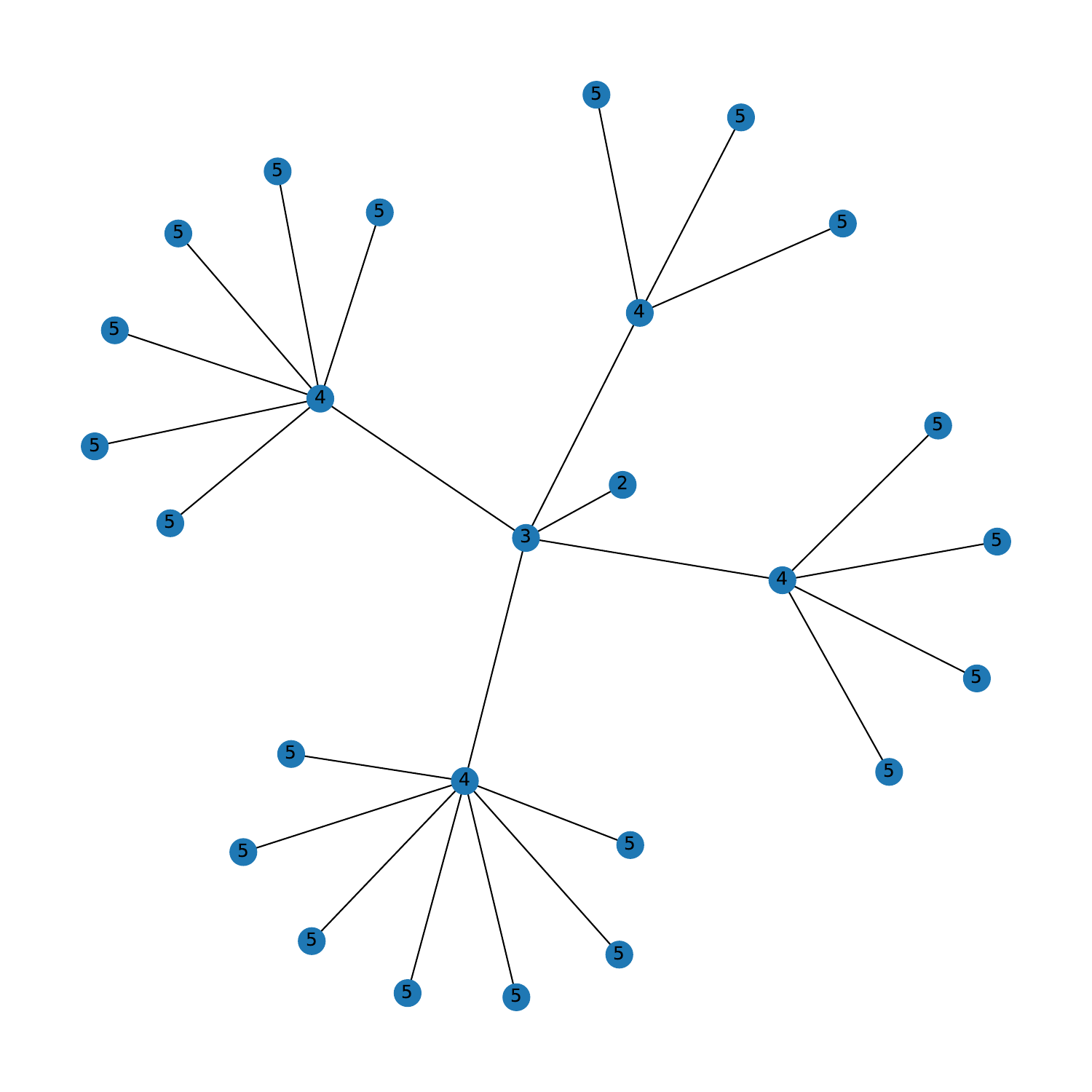}
        \caption*{\footnotesize{Figure 17:} $S^2 \times S^1$ Depth 3}
    \end{subfigure} 
    \begin{subfigure}{0.47\textwidth}
        \centering
        \includegraphics[width=0.98\textwidth]{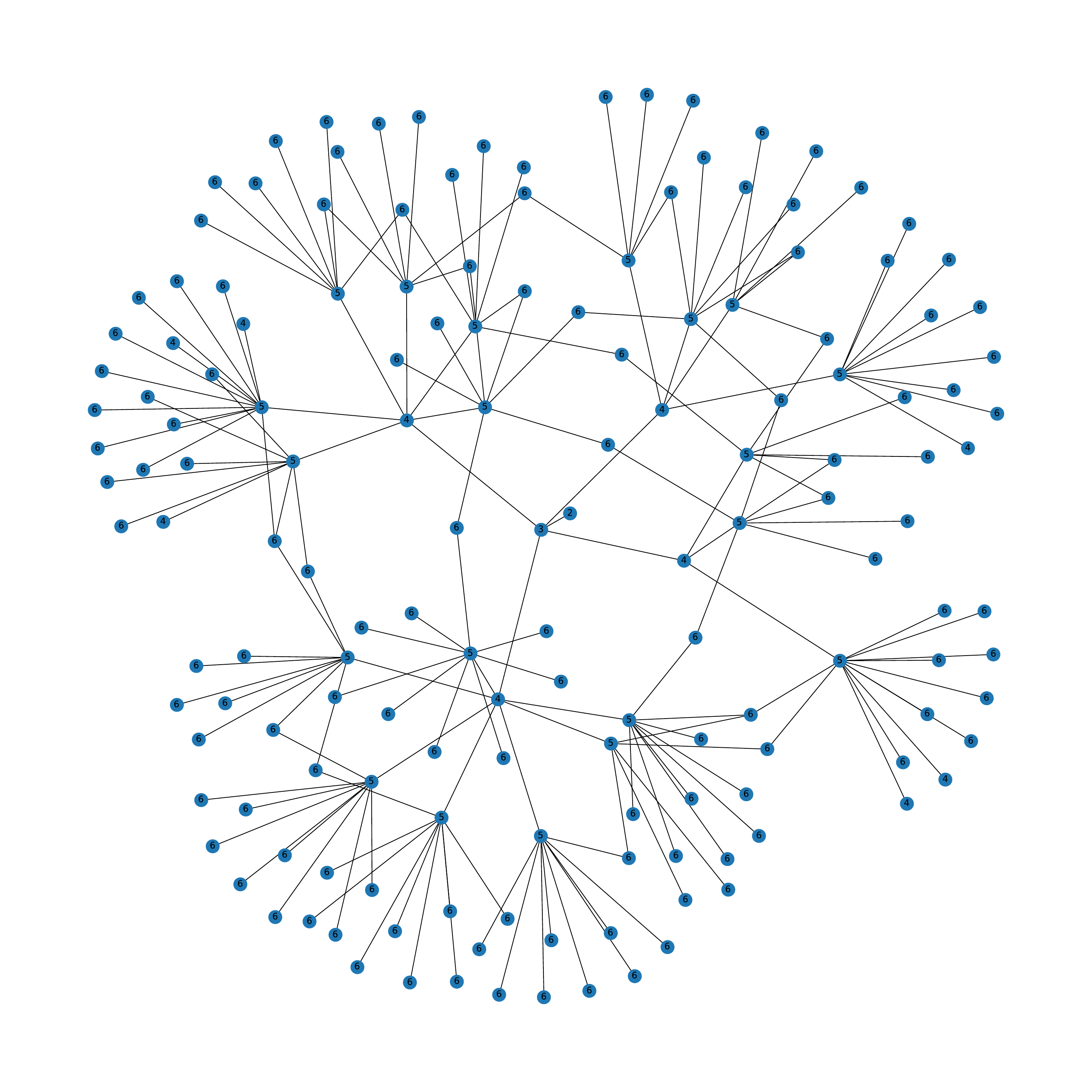}
        \caption*{\footnotesize{Figure 18:} $S^2 \times S^1$ Depth 4}
    \end{subfigure}\\
    \begin{subfigure}{0.47\textwidth}
        \centering
        \includegraphics[width=0.98\textwidth]{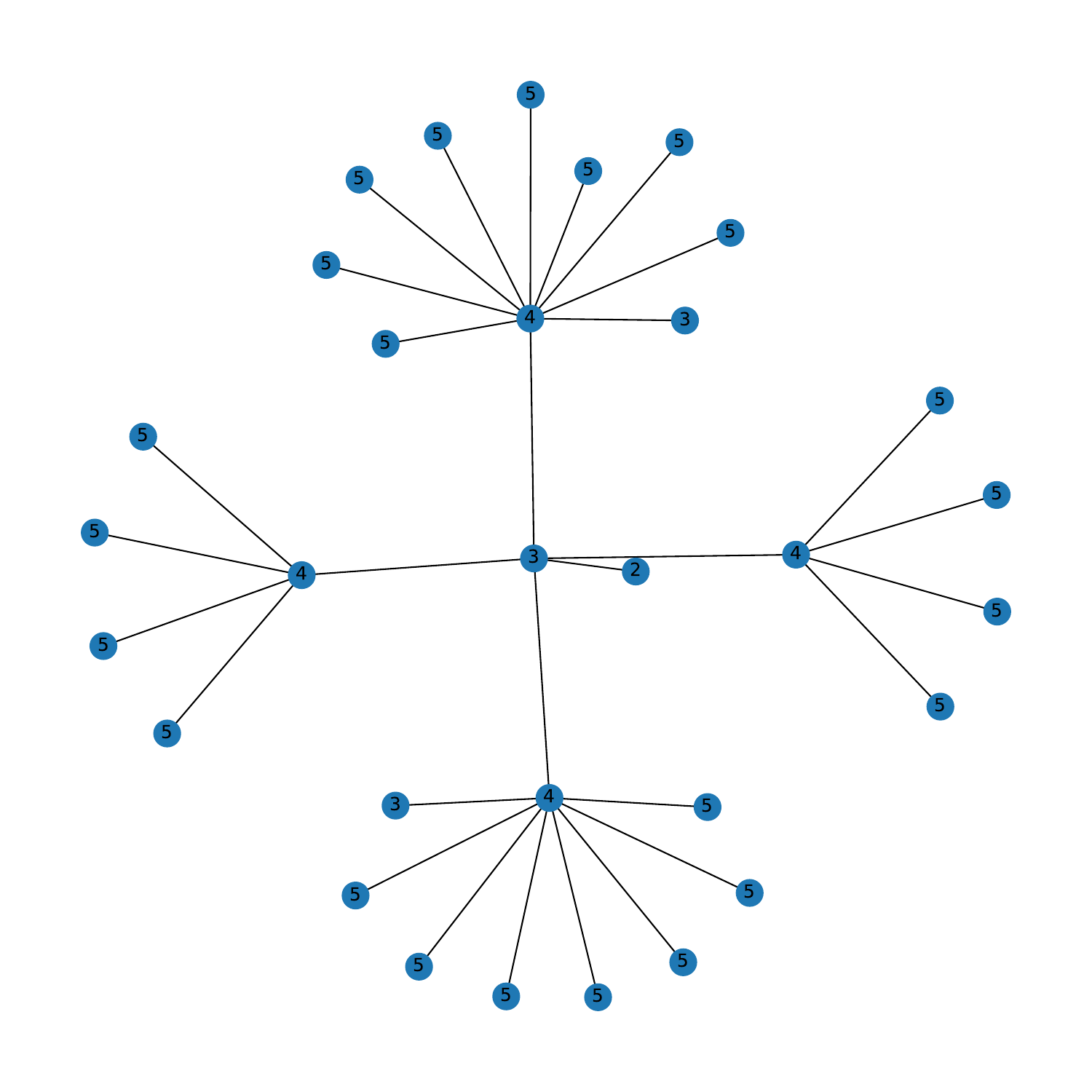}
        \caption*{\footnotesize{Figure 19:} $\mathbb{RP}^3$ Depth 3}
    \end{subfigure} 
    \begin{subfigure}{0.47\textwidth}
        \centering
        \includegraphics[width=0.98\textwidth]{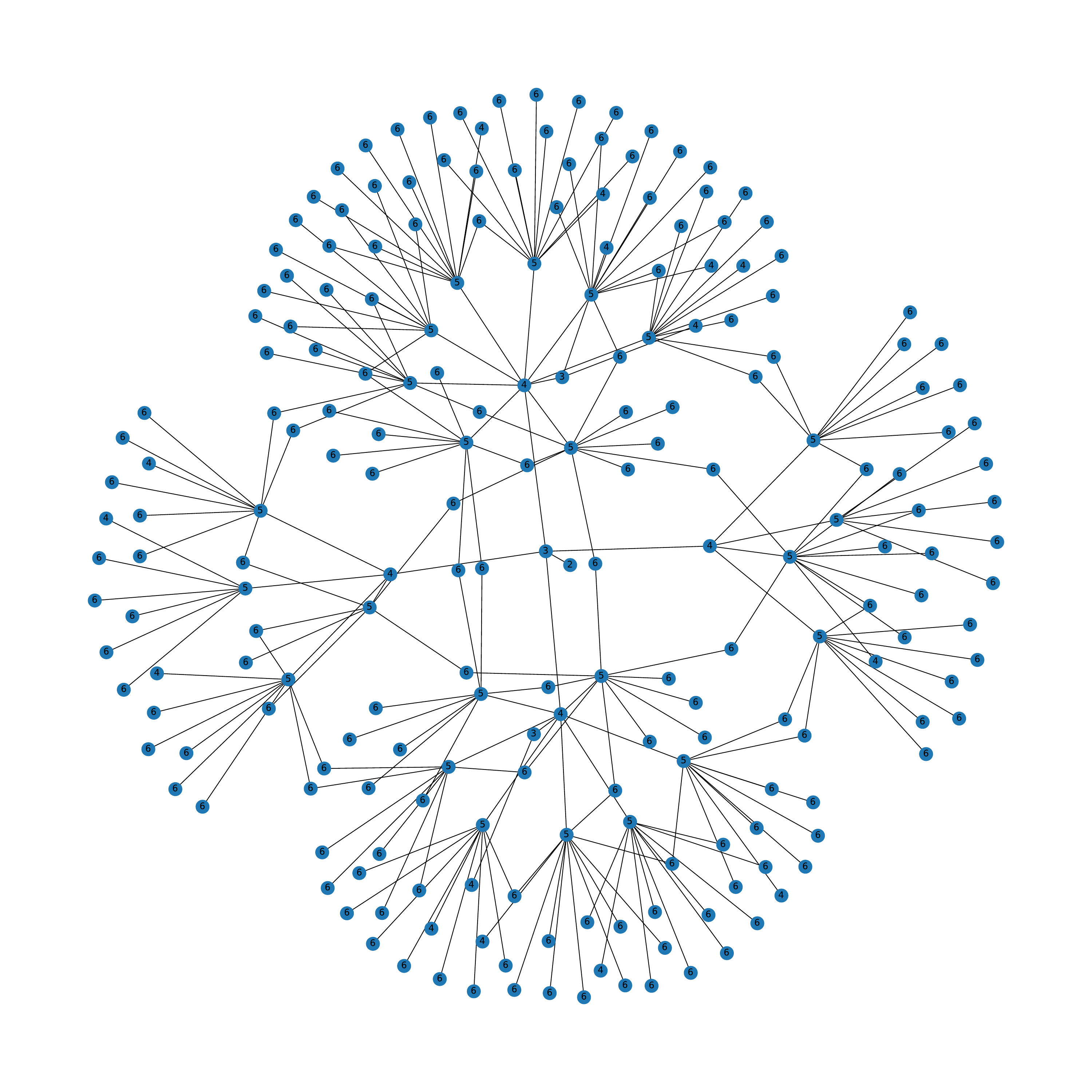}
        \caption*{\footnotesize{Figure 20:} $\mathbb{RP}^3$ Depth 4}
    \end{subfigure}
    \label{extra_PGs1}
\end{figure}

\begin{figure}[b]
    \centering
    \begin{subfigure}{0.47\textwidth}
        \centering
        \includegraphics[width=0.98\textwidth]{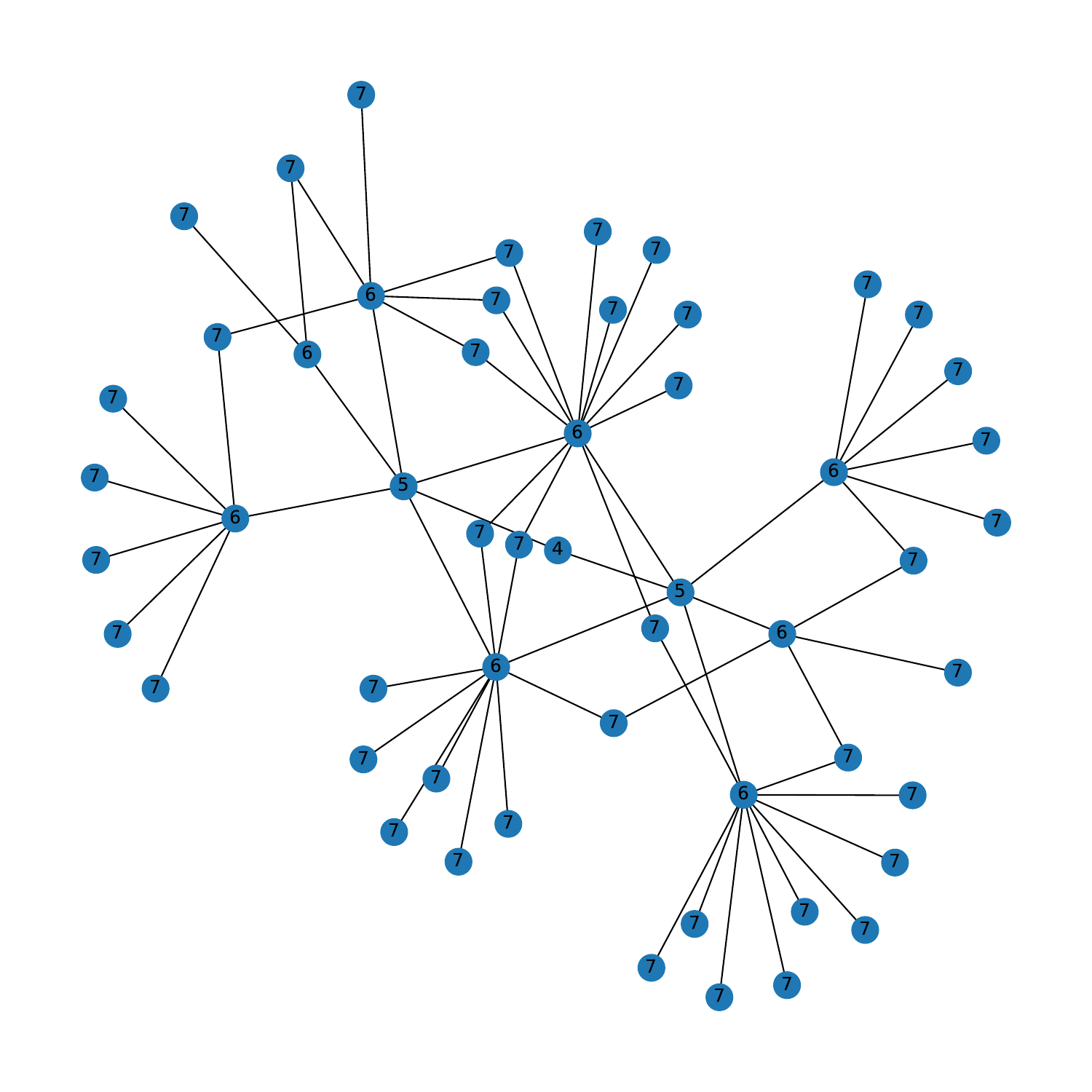}
        \caption*{\footnotesize{Figure 21:} $L(7,1)$ Depth 3}
    \end{subfigure} 
    \begin{subfigure}{0.47\textwidth}
        \centering
        \includegraphics[width=0.98\textwidth]{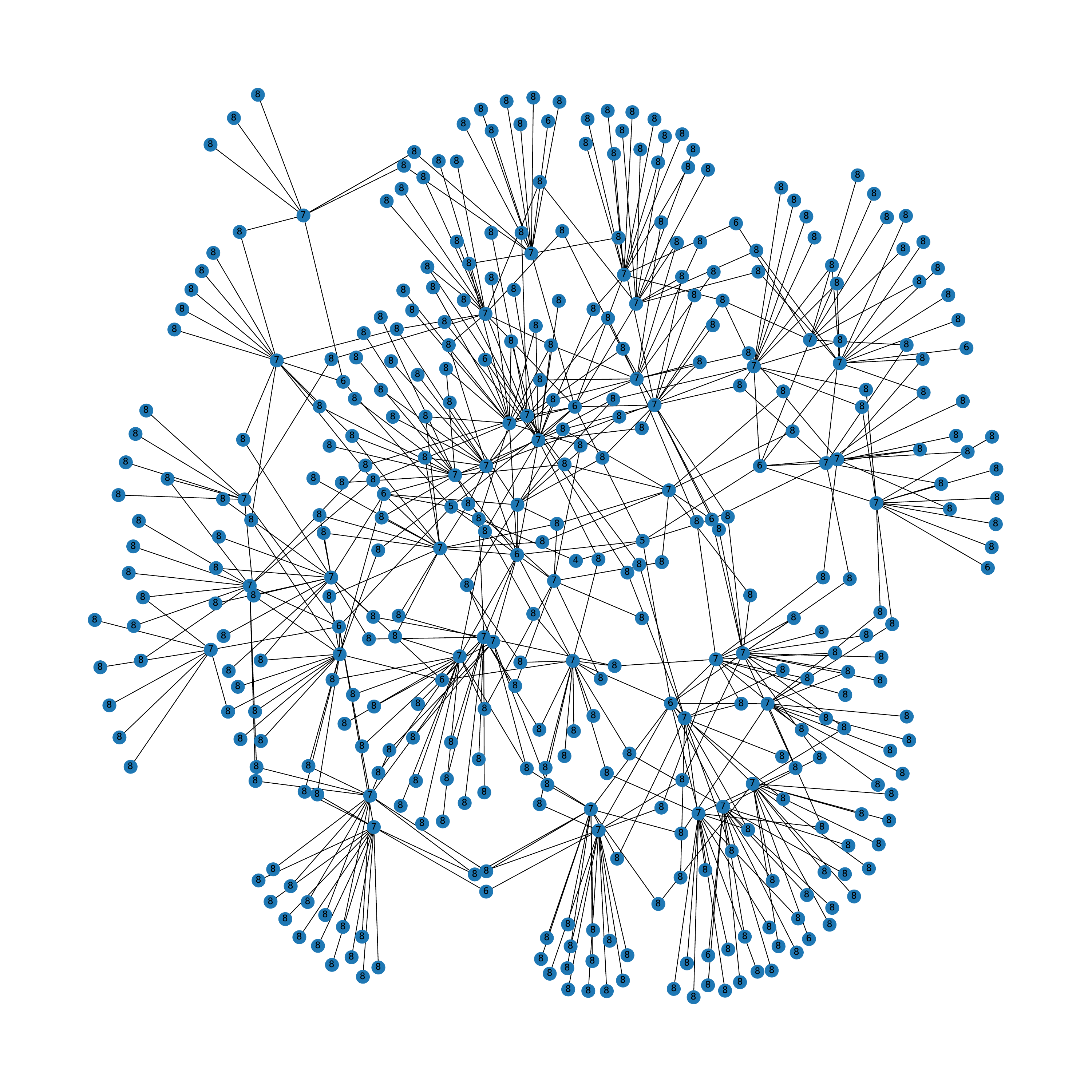}
        \caption*{\footnotesize{Figure 22:} $L(7,1)$ Depth 4}
    \end{subfigure}
        \begin{subfigure}{0.47\textwidth}
        \centering
        \includegraphics[width=0.98\textwidth]{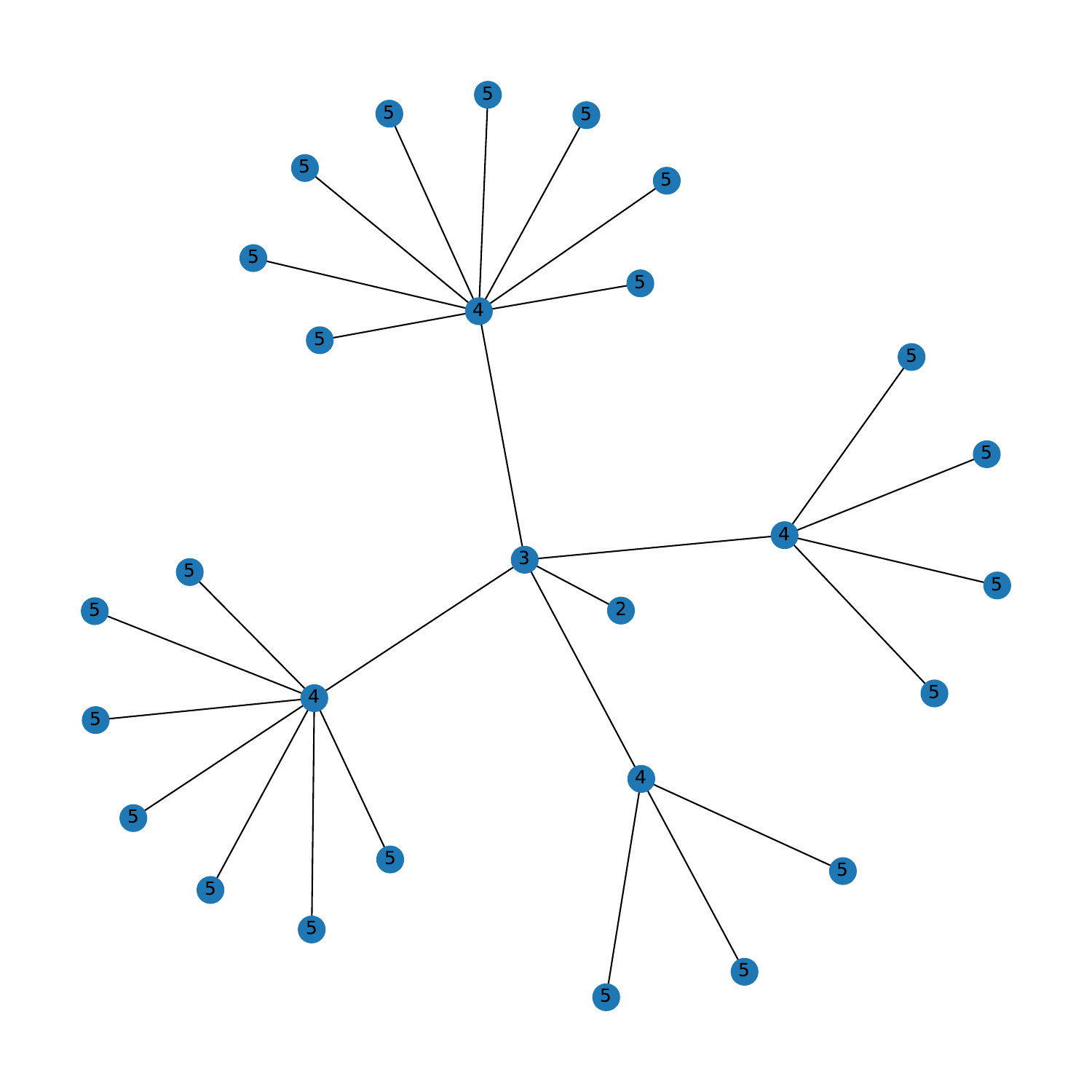}
        \caption*{\footnotesize{Figure 23:} $L(7,2)$ Depth 3}
    \end{subfigure} 
    \begin{subfigure}{0.47\textwidth}
        \centering
        \includegraphics[width=0.98\textwidth]{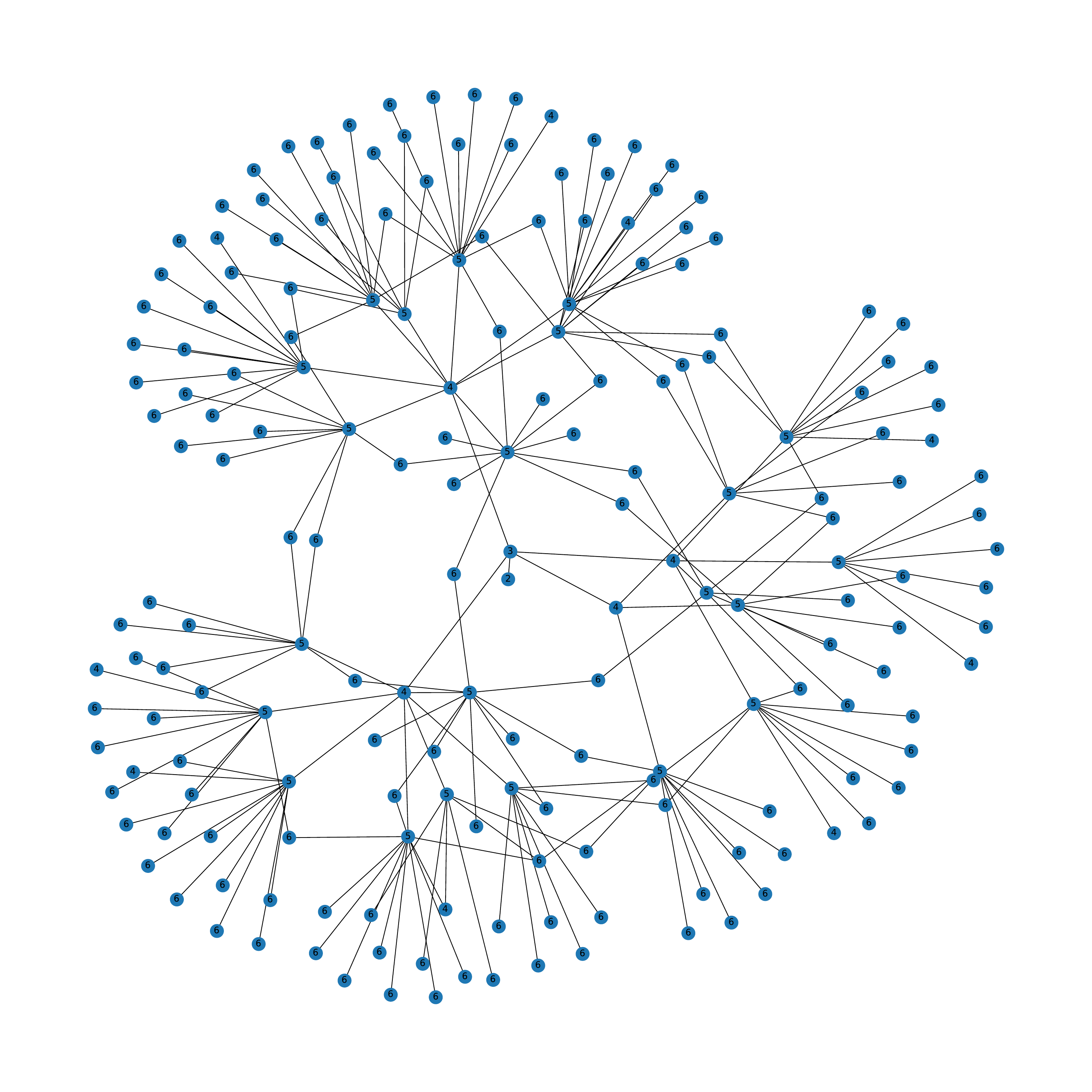}
        \caption*{\footnotesize{Figure 24:} $L(7,2)$ Depth 4}
    \end{subfigure}
    \label{extra_PGs2}
\end{figure}

\begin{figure}[b]
    \centering
    \begin{subfigure}{0.47\textwidth}
        \centering
        \includegraphics[width=0.98\textwidth]{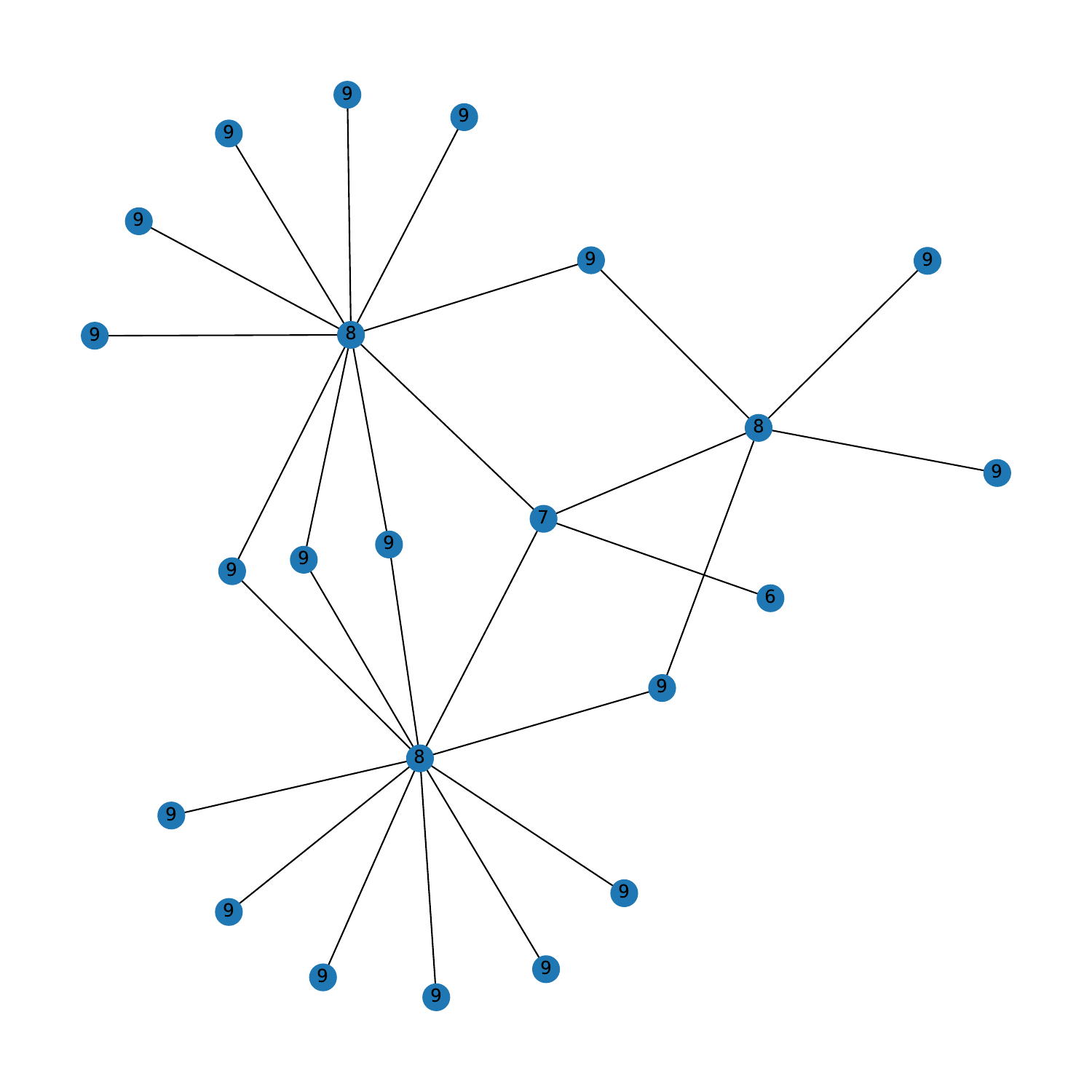}
        \caption*{\footnotesize{Figure 25:} $T^3$ Depth 3}
    \end{subfigure} 
    \begin{subfigure}{0.47\textwidth}
        \centering
        \includegraphics[width=0.98\textwidth]{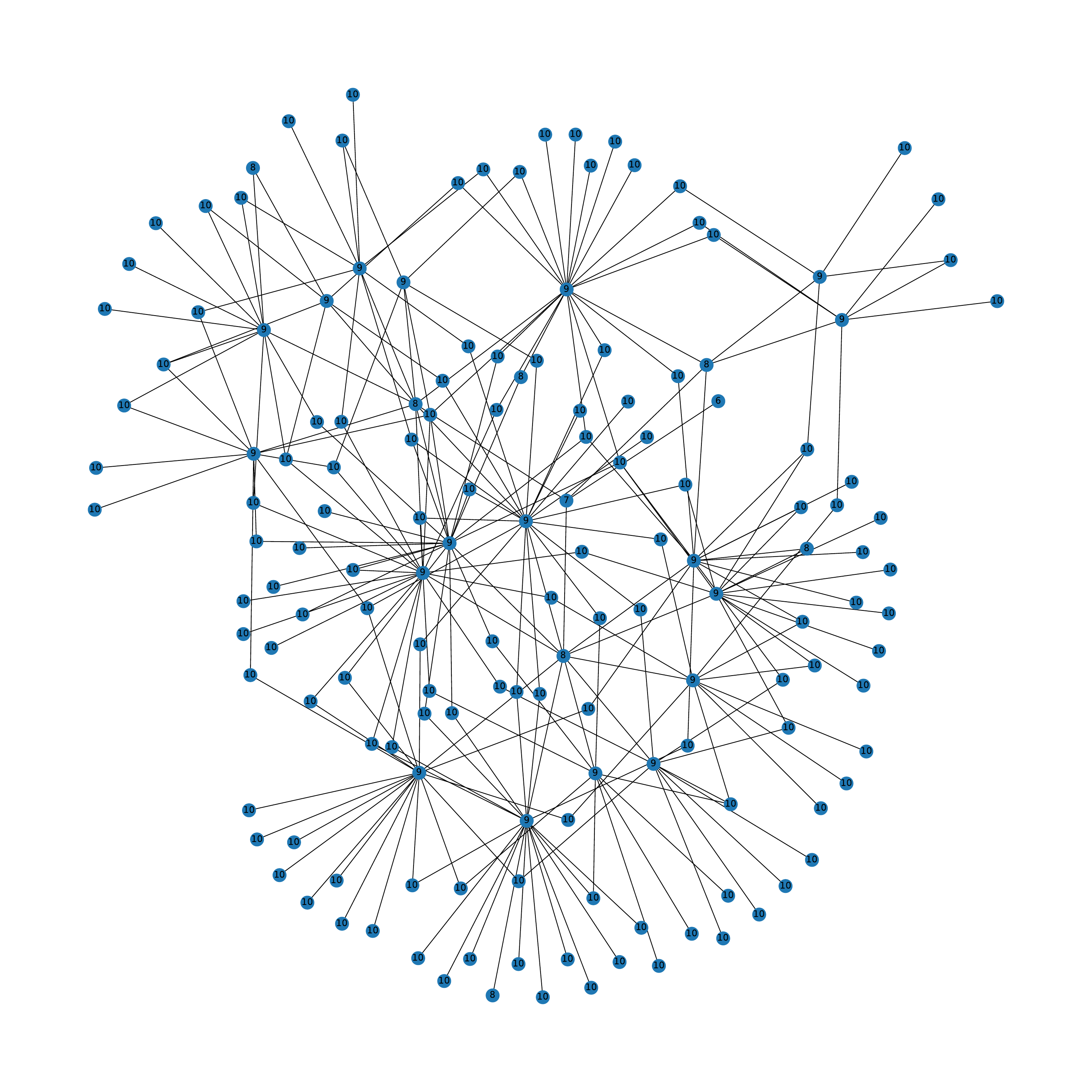}
        \caption*{\footnotesize{Figure 26:} $T^3$ Depth 4}
    \end{subfigure}
    \begin{subfigure}{0.47\textwidth}
        \centering
        \includegraphics[width=0.98\textwidth]{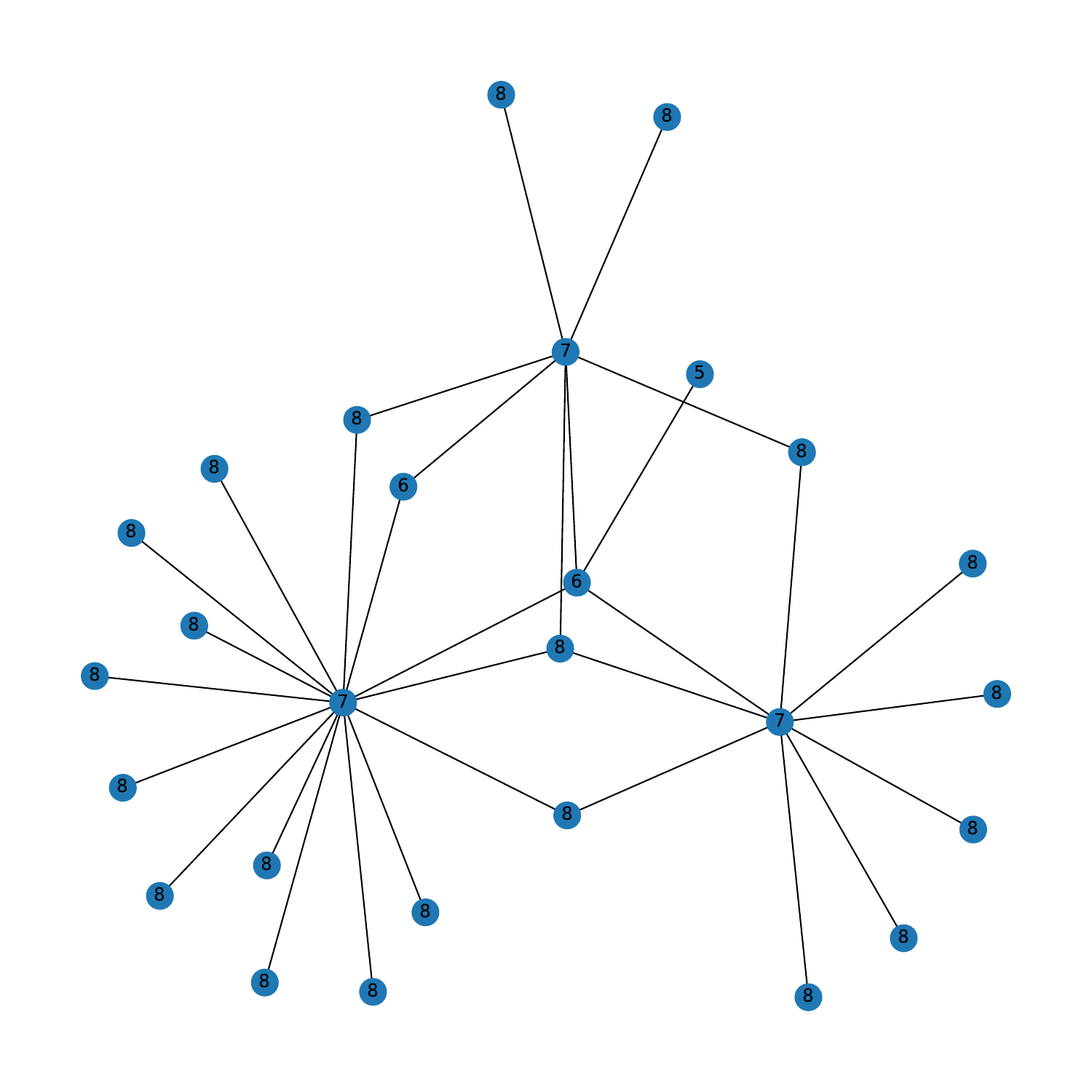}
        \caption*{\footnotesize{Figure 27:} PHS Depth 3}
    \end{subfigure} 
    \begin{subfigure}{0.47\textwidth}
        \centering
        \includegraphics[width=0.98\textwidth]{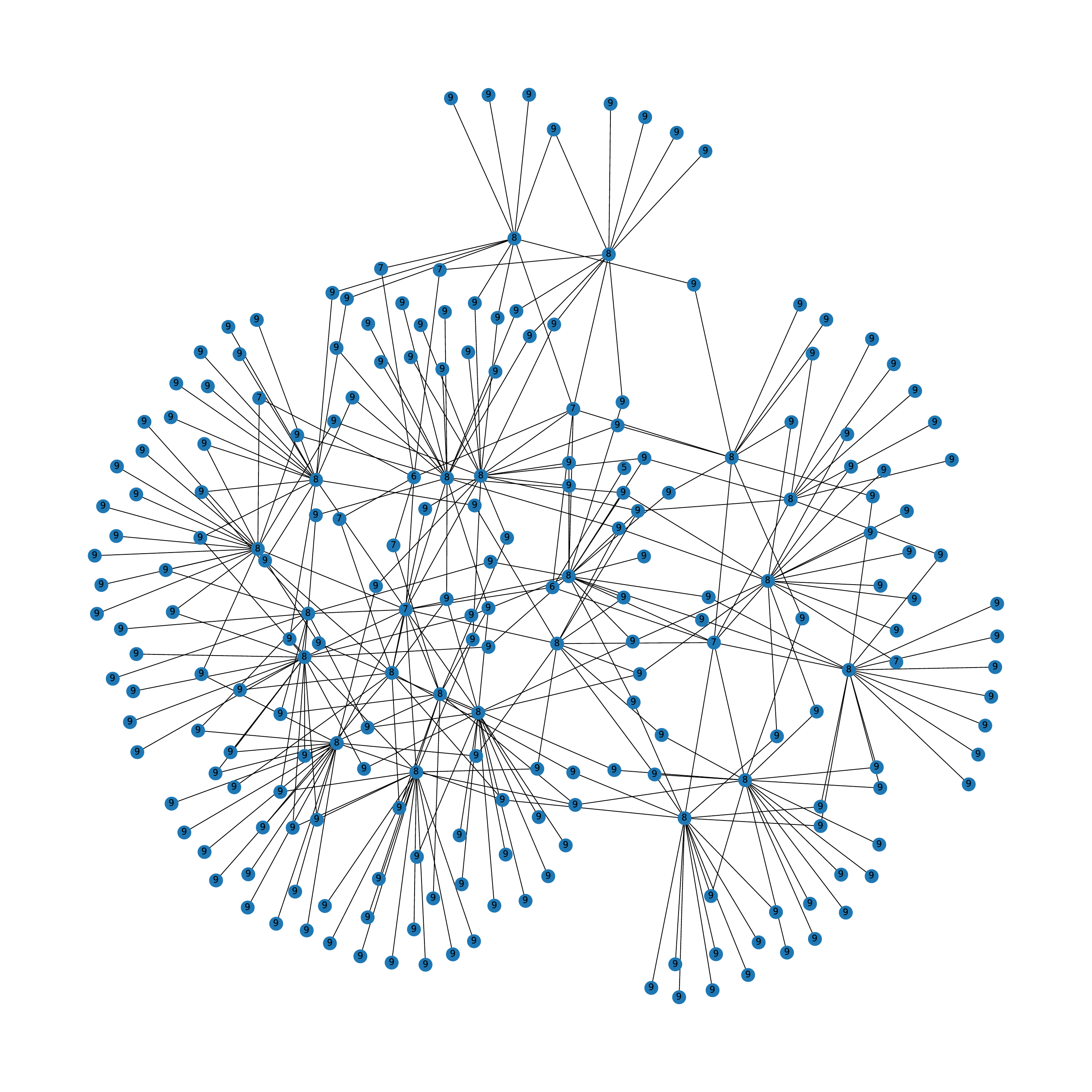}
        \caption*{\footnotesize{Figure 28:} PHS Depth 4}
    \end{subfigure}
    \label{extra_PGs3}
\end{figure}

\begin{figure}[b]
    \centering
    \begin{subfigure}{0.47\textwidth}
        \centering
        \includegraphics[width=0.98\textwidth]{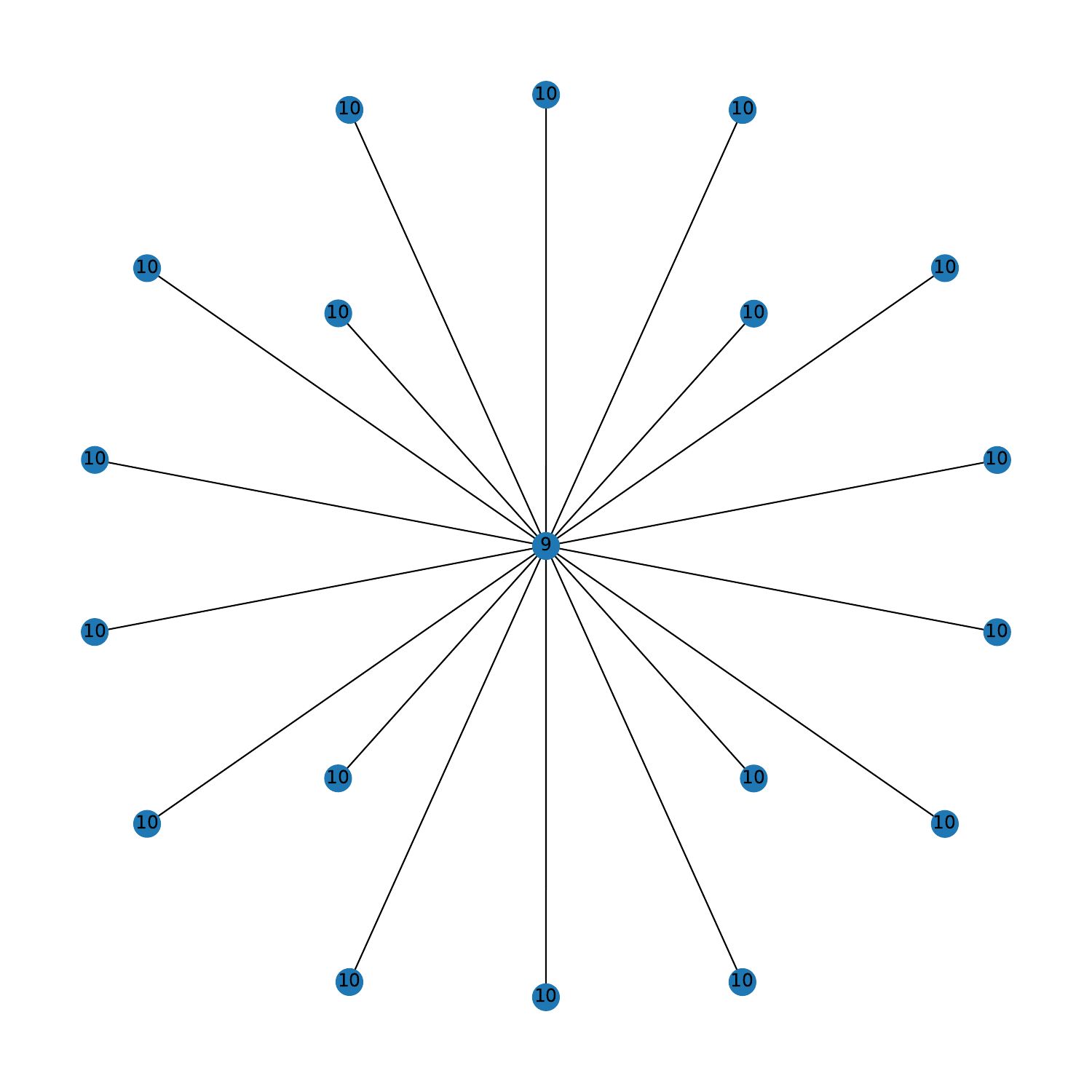}
        \caption*{\footnotesize{Figure 29:} $H_{SC}$ Depth 1}
    \end{subfigure} 
    \begin{subfigure}{0.47\textwidth}
        \centering
        \includegraphics[width=0.98\textwidth]{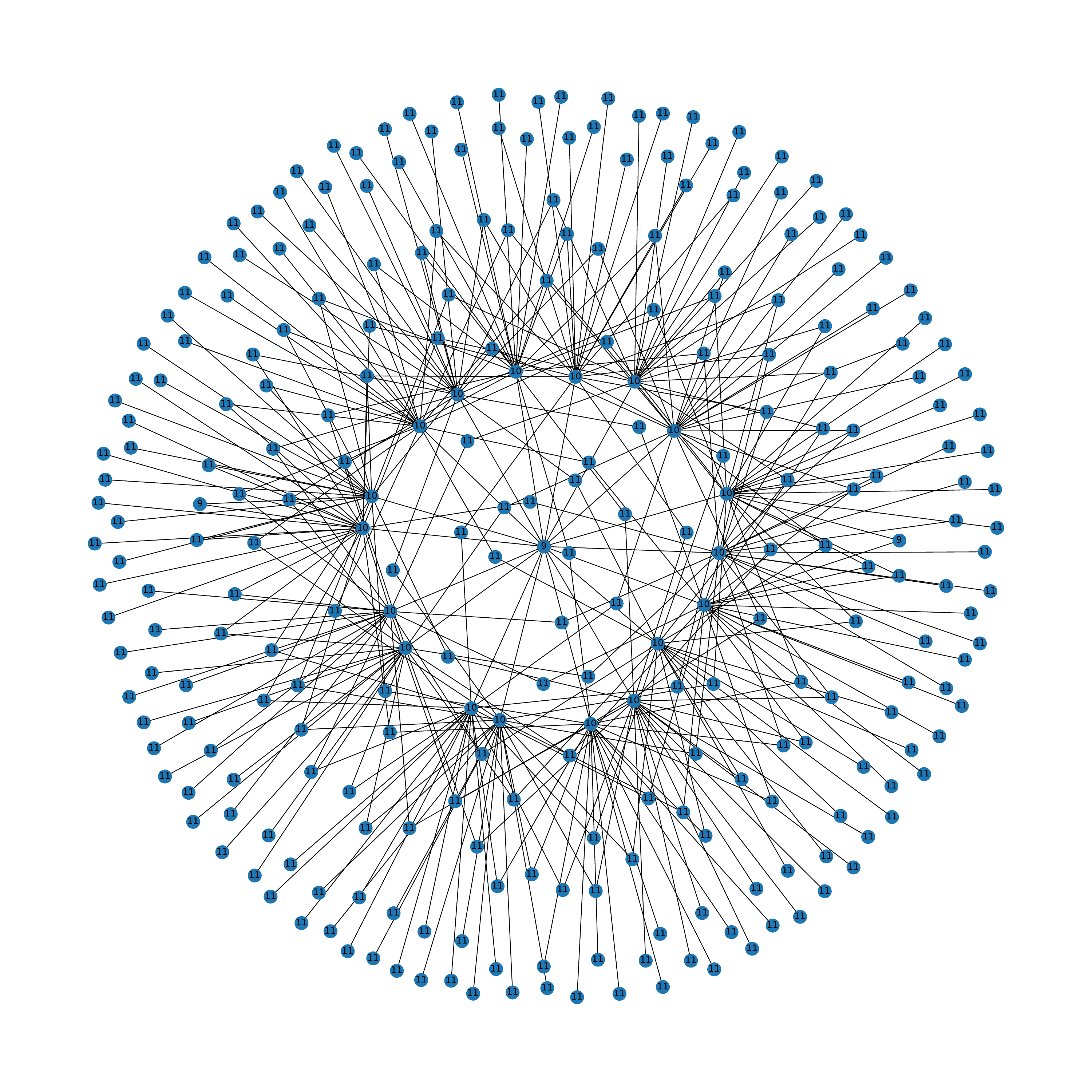}
        \caption*{\footnotesize{Figure 30:} $H_{SC}$ Depth 2}
    \end{subfigure}
    \label{extra_PGs4}
\end{figure}

\begin{figure}[!t]
    \centering    
    \includegraphics[width=0.7\textwidth]{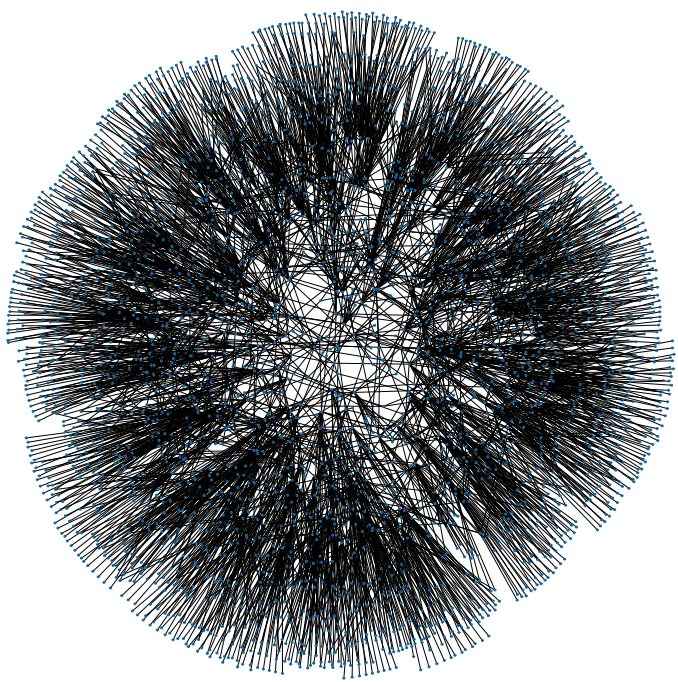}
    \caption*{\footnotesize{Figure 31:} $S^3$ depth 5}
    \label{PGS3_5}
\end{figure}

\end{document}